\documentclass[reqno, 11pt]{amsart}
\usepackage[usenames,dvipsnames]{color}
\usepackage{hyperref}
\hypersetup{
    colorlinks=true,       			% false: boxed links; true: colored links
    linkcolor=blue,          			% color of internal links
    citecolor=blue,		 		% color of links to bibliography
    filecolor=blue,      				% color of file links
    urlcolor=blue           			% color of external links
}

\usepackage{amsmath, amssymb, hyperref, empheq}
\usepackage[perpage,multiple]{footmisc}

\usepackage{tikz}
\usetikzlibrary{shapes,backgrounds}

\usepackage[all,cmtip,knot]{xy}
\xyoption{arc}
\usepackage{mathtools}

\newtheorem* {theorem}{Theorem}
\newtheorem* {proposition}{Proposition}
\newtheorem* {corollary}{Corollary}
\newtheorem* {lemma}{Lemma}

\theoremstyle{definition}
\newtheorem* {definition}{Definition}
\newtheorem*{rem}{Remark}

\numberwithin {equation}{section}

\newenvironment{pf}{\paragraph{{\sc Proof}}}{\qed\par\medskip}

\newcommand {\IC}{\mathbb{C}}
\newcommand {\IN}{\mathbb{N}}

\newcommand {\IZ}{\mathbb{Z}}

\newcommand {\bfI}{\mathbf I}

\newcommand {\A}{\mathcal A}

\newcommand {\C}{\mathcal C}
\newcommand {\D}{\mathcal D}

\newcommand {\F}{\mathcal F}
\newcommand {\G}{\mathcal G}
\renewcommand {\H}{\mathcal H}
\newcommand {\I}{\mathcal I}

\newcommand {\K}{\mathcal K}

\newcommand {\N}{\mathcal N}
\renewcommand {\O}{\mathcal O}

\newcommand {\R}{\mathcal R}

\newcommand {\V}{\mathcal V}

\renewcommand {\a}{\mathfrak a}
\renewcommand {\b}{\mathfrak b}
\renewcommand {\c}{\mathfrak c}

\newcommand {\g}{\mathfrak{g}}
\newcommand {\h}{\mathfrak h}
\renewcommand {\ll}{\mathfrak l}
\newcommand {\m}{\mathfrak m}
\newcommand {\n}{\mathfrak n}

\renewcommand {\r}{\mathfrak r}
\newcommand {\z}{\mathfrak z}

\renewcommand {\SS}{\mathfrak S}

\newcommand {\sfA}{\mathsf{A}}

\newcommand {\ie}{{\it i.e., }}
\newcommand {\eg}{{\it e.g.}, }

\newcommand {\wrt}{with respect to }

\newcommand {\ol}{\overline}
\newcommand {\wt}{\widetilde}
\newcommand {\wh}{\widehat}

\newcommand {\ns}{nested set }

\newcommand {\supp}{\operatorname{supp}}

\newcommand {\DCP}{De Concini--Procesi }
\renewcommand {\DJ}{Drinfeld--Jimbo }

\newcommand {\nEK}{Etingof--Kazhdan }

\newcommand {\KM}{Kac--Moody }
\newcommand {\AKM}{affine Kac--Moody }

\newcommand {\EEK}{^{\scriptscriptstyle{\operatorname{EK}}}}

\newcommand{\Lg}{\mathfrak{g}}
\newcommand{\Lh}{\mathfrak{h}}

\newcommand{\Lp}{\mathfrak{p}}

\newcommand {\gD}{\mathfrak{g}_{D}}
\newcommand {\gDp}{\mathfrak{g}_{D,+}}
\newcommand {\gDm}{\mathfrak{g}_{D,-}}
\newcommand {\gDpm}{\mathfrak{g}_{D,\pm}}

\newcommand {\iD}{i_D}

\newcommand {\gp}{\mathfrak{g}_+}
\newcommand {\gm}{\mathfrak{g}_-}
\newcommand{\gpm}{\mathfrak{g}_{\pm}}
\newcommand{\gmp}{\mathfrak{g}_{\mp}}

\newcommand{\Lpp}{\mathfrak{p_+}}
\newcommand{\Lpm}{\mathfrak{p_-}}
\newcommand{\Lppm}{\mathfrak{p}_{\pm}}

\newcommand{\Lmp}{\mathfrak{m_+}}
\newcommand{\Lmm}{\mathfrak{m_-}}
\newcommand{\Lmpm}{\mathfrak{m_{\pm}}}

\newcommand{\ctp}{\hat{\otimes}}
\newcommand{\ten}{\otimes}

\def\iip#1#2{\langle{#1},{#2}\rangle}

\def\slantfrac#1#2{\kern.2em\raise.2em\hbox{$#1$}\kern-.2em\left/\lower.4em\hbox{$#2$}\kern-.2em\right.}
\def\backslantfrac#1#2{\kern.2em\lower.4em\hbox{$#1$}\kern-.3em\left\backslash\raise.4em\hbox{$#2$}\right.}

\DeclareMathOperator{\Hom}{Hom}
\DeclareMathOperator{\End}{End}
\DeclareMathOperator{\Ker}{Ker}
\DeclareMathOperator{\colim}{colim}
\DeclareMathOperator{\Ind}{Ind}

\DeclareMathOperator{\id}{id}

\DeclareMathOperator{\rank}{rank}
\DeclareMathOperator{\corank}{corank}
\DeclareMathOperator{\eRep}{Rep^{\operatorname{eq}}}

\renewcommand {\sl}[1]{\mathfrak{sl}_{#1}}

\newcommand {\Uhg}{U_{\hbar}\g}

\renewcommand{\1}{\mathbf{1}}
\newcommand{\Mp}{M_+}
\newcommand{\Mm}{M_-}
\newcommand{\Mpd}{\Mp^*}
\newcommand{\Lm}{L_-}
\newcommand{\Np}{N_+}
\newcommand{\Npd}{\Np^*}

\newcommand{\um}{1_-}
\newcommand{\up}{1_+}
\newcommand{\coup}{\epsilon_+}
\newcommand{\cou}{\epsilon}

\newcommand{\ipd}{i_+^{\vee}}
\newcommand{\im}{i_-}
\newcommand{\ips}{i_+^*}
\newcommand{\ip}{i_+}

\newcommand{\EK}{U_{\hbar}^{\EEK}}
\renewcommand{\DJ}{U_{\hbar}}
\def\qEK#1{\EK#1}
\def\Ue#1{U#1}
\def\UE#1{U(#1)}
\def\Ueh#1{U#1[[\hbar]]}
\newcommand{\oEK}{\EEK}

\newcommand{\DY}{\mathcal{YD}}
\newcommand{\TDY}[2]{\DY_{#2}(#1)}
\newcommand{\DC}[2]{\D_{#2}(#1)}

\newcommand{\oKZ}{^{\operatorname{KZ}}}

\newcommand{\Oint}{\O^{\operatorname{int}}}

\newcommand {\odots}[1]{#1\cdots #1}

\newcommand {\aand}{\qquad\text{and}\qquad}

\newcommand {\qc}{quasi--Coxeter }

\newcommand {\qcqtqba}{quasi--Coxeter quasitriangular quasibialgebra }

\renewcommand {\t}{t}

\newcommand {\ND}{\N_D}

\newcommand {\UhgD}{U_{\hbar}\gD}

\newcommand {\ul}[1]{\underline{#1}}

\newcommand {\nperp}{\perp\negthickspace\negthickspace\negthickspace\diagup}

\newcommand {\Db}{\mathfrak{D}}
\newcommand{\sfMns}{\mathsf{Mns}}

\DeclareMathOperator{\Rep}{Rep}
\DeclareMathOperator{\obj}{Obj}
\DeclareMathOperator{\vect}{Vect}
\DeclareMathOperator{\ad}{ad}
\newcommand {\sfk}{\mathsf{k}}

\newcommand {\kvect}{\vect_{\sfk[[\hbar]]}}

\newcommand{\Csd}[1]{\mathsf{SD}(#1)}%{\sfC(#1)}
\newcommand{\Mns}[1]{\sfMns(#1)}

\newcommand{\Omit}[1]{}
\newcommand{\ind}{\operatorname{Ind}}
\newcommand{\PLBA}{\underline{PLBA} }
\newcommand{\PLBAD}{\underline{PLBAD} }
\newcommand{\LBA}{\underline{LBA} }
\newcommand{\HA}{\underline{HA} }
\newcommand{\PHA}{\underline{PHA} }
\newcommand{\PHAD}{\underline{PHAD} }
\newcommand{\PROP}{{\sc Prop} }
\newcommand{\PROPs}{{\sc Prop}s }
\newcommand{\hLBA}{S({\LBA}) }
\newcommand{\hPLBA}{S({\PLBA}) }
\newcommand{\hPLBAD}{S({\PLBAD}) }
\newcommand{\Gammah}{\Gamma_{\hbar}}
\newcommand{\Lmh}{\Lm^{\hbar}}

\newcommand{\Npdh}{(N_+^{\hbar})^*}

\newcommand{\sfEnd}[1]{\mathsf{End}(#1)}
\newcommand{\Nat}[2]{\mathsf{Nat}(#1,#2)}
\newcommand{\sfNat}{\mathsf{Nat}}
\newcommand{\sfAd}[1]{\mathsf{Ad}(#1)}
\newcommand{\sfAut}[1]{\mathsf{Aut}(#1)}

\newcommand{\Cat}{\mathsf{Cat}}
\newcommand{\qC}[1]{\mathsf{qC}(#1)}

\newcommand{\sfV}{\mathsf{V}}
\newcommand{\bfJ}{\mathbf{J}}

\newcommand {\abelian}{additive }%{abelian }
\newcommand {\exact}{additive }%{exact }
\newcommand {\ID}{I(D)}
\newcommand {\MD}{M(D)}
\newcommand {\Drin}{\mathsf{DCat^\ten}}

\newcommand{\return}{{\color{white}.}\linebreak\vspace{-0.2in}}

\title[Quasi--Coxeter categories and a relative EK functor]{Quasi--Coxeter categories and a relative Etingof--Kazhdan quantization functor}
\author[A. Appel]{Andrea Appel}
\address{Department of Mathematics,
Northeastern University,
360 Huntington Avenue,
Boston, MA 02115}
\email{appel.an@husky.neu.edu}
\author[V. Toledano Laredo]{Valerio Toledano Laredo}
\address{Department of Mathematics,
Northeastern University,
360 Huntington Avenue,
Boston MA 02115}
\email{V.ToledanoLaredo@neu.edu}
\thanks{Both authors were supported in part through the NSF grant DMS--0854792}
\date{\today}

\begin{document}
%\dedicatory{Preliminary version}
\maketitle

\begin{abstract}
Let $\g$ be a symmetrizable \KM algebra and $\DJ{\g}$ its quantized
enveloping algebra. The quantum Weyl group operators of $\DJ{\g}$
and the universal $R$--matrices of its Levi subalgebras endow $\DJ
{\g}$ with a natural \qcqtqba structure which underlies the action of
the braid group of $\g$ and Artin's braid groups on the tensor product
of integrable, category $\O$ modules. We show that this structure can
be transferred to the universal enveloping algebra $U\g[[\hbar]]$. The
proof relies on a modification of the Etingof--Kazhdan quantization
functor, and yields an isomorphism between (appropriate completions
of) $\DJ{\g}$ and $U\g[[\hbar]]$ preserving a given chain of Levi
subalgebras. We carry it out in the more general context of chains
of Manin triples, and obtain in particular a relative version of the
Etingof--Kazhdan functor with input a split pair of Lie bialgebras.
Along the way, we develop the notion of quasi--Coxeter categories,
which are to generalized braid groups what braided tensor categories
are to Artin's braid groups. This leads to their succint description as a
2--functor from a 2--category whose morphisms are \DCP associahedra.
These results will be used in the sequel to this paper to give a monodromic
description of the quantum Weyl group operators of an affine \KM algebra,
extending the one obtained by the second author for a semisimple Lie algebra.
\end{abstract}

\tableofcontents

\section{Introduction}
%===============

\subsection{}
%--------------

This is the first of a series of three papers the aim of which is to extend
the description of the monodromy of the rational Casimir connection in
terms of quantum Weyl group operators given in \cite{vtl-3,vtl-4,vtl-6}
to the case of an affine \KM algebra $\g$.

The method we follow is close to that of \cite{vtl-4}, and relies on the notion
of a \qcqtqba (qCqtqba), which is informally a bialgebra carrying actions
of a given generalized braid group and Artin's braid groups on the tensor
products of its modules. A cohomological rigidity result, proved in the second
paper of this series \cite{ATL2}, shows that there is at most one such structure
with prescribed local monodromies on the classical enveloping algebra
$\Ueh{\g}$. It follows that the generalized braid group actions arising from
quantum Weyl groups and the monodromy of the Casimir connection \cite
{ATL3} are equivalent, provided the quasi--Coxeter quasitriangular quasibialgebra
structure responsible for the former can be transferred from $\DJ{\g}$ to
$\Ueh{\g}$. This result is the purpose of the present article.

\subsection{}
%--------------

Its proof differs substantially from that given in \cite{vtl-4}. Indeed, for
a semisimple Lie algebra $\g$, the transfer of structure ultimately rests
on the vanishing of the first and second Hochschild cohomology groups
of $\Ueh{\g}$, and in particular on the fact that $\DJ{\g}$ and $\Ueh{\g}$
are isomorphic as algebras, a fact which does not hold for \AKM algebras.
Rather than the cohomological methods of \cite{vtl-4}, we use instead
the Etingof--Kazhdan (EK) quantization functor \cite{ek-1,ek-2,ek-6},
which yields a canonical isomorphism
\[\Psi\EEK:\wh{\DJ{\g}}\stackrel{\sim}{\longrightarrow}\wh{\Ueh{\g}}\]
between the completions of $\DJ{\g}$ and $\Ueh{\g}$ \wrt category
$\O$.

Surprisingly perhaps, and despite its functorial construction, the
isomorphism $\Psi\EEK$ does not preserve the inclusions of Levi
subalgebras
\[\DJ{\g_D}\subseteq \DJ{\g}\qquad\text{and}\qquad\Ueh{\g_D}\subseteq\Ueh{\g}\]
determined by a subdiagram $D$ of the Dynkin diagram of $\g$,
something which is required by the transfer of structure. The bulk
of this paper is therefore devoted to modifying $\Psi\EEK$ so as to
make it compatible with such inclusions.

\subsection{} % Etingof-Kazhdan functor and isomorphism
%--------------

To outline our construction, which works more generally for an
inclusion $(\g_D,\g_{D,-},\g_{D,+})\subset(\g,\g_-,\g_+)$ of Manin
triples over a field $\sfk$ of characteristic zero, recall first that the
main steps of the EK construction are as follows.
\begin{enumerate}
% Drinfeld category
\item One considers the Drinfeld category $\D_\Phi(\g)$ of (deformation)
equicontinuous $\g$--modules, with associativity constraints given by a
fixed Lie associator $\Phi$ over $\sfk$. This category can be thought of
as a topological analogue of category $\O$ when $\g$ is the Manin triple
associated to a Kac--Moody algebra. It can equivalently be described as
the category of Drinfeld--Yetter modules over the Lie bialgebra $\g_-$.
% Tensor functor F (==> End(F) is a topological bialgebra which obviously acts on any F(V))
\item One constructs a tensor functor $F$ from $\D_\Phi(\g)$ to
the category $\kvect$ of topologically free $\sfk[[\hbar]]$--modules.
The algebra of endomorphisms $H=\sfEnd{F}$ is then a topological
bialgebra, \ie it is endowed with a coproduct $\Delta$ mapping $H$
to a completion of $H\otimes H$. 
%,namely $\sfEnd{F\boxtimes F}$.
% The Hopf algebra U_h g_-\subset End(F)
\item Inside $H$, one constructs a subalgebra $\DJ{\g_-}$ such that
$\Delta(\DJ{\g_-})\subset\DJ{\g_-}\otimes\DJ{\g_-}$, and which is a
quantization of $U\g_-$. The quantum group $\DJ{\g}$ is then defined
as the quantum double of $\DJ{\g_-}$.
%It is a dense subbialgebra of $H$.
% The lift of F
\item By construction, $U\g_-$ acts and coacts on any $F(V)$, $V\in
\D_\Phi(\g)$, so that the functor $F$ lifts to $\wt{F}:\D_\Phi(\g)\to\Rep
(\DJ{\g})$ where, by definition, the latter is the category of Drinfeld--Yetter
modules over $\DJ(\g_-)$.
% F is an equivalence of categories
\item Finally, one proves that $\wt{F}$ is an equivalence of categories.
\end{enumerate}

Since $F$ is isomorphic to the forgetful functor $f:\D_\Phi(\g)\to\kvect$
as abelian functors, we obtain the following diagram
\[\xymatrix{
%\ar@{}[ddrr]|{\stackrel{\alpha}{\leftarrow}}
\D_\Phi(\g)\ar@{=}[rr]\ar[dd]_{f}&&\D_\Phi(\g)\ar@{->}[rr]^{\wt{F}}\ar[dd]^{F}&&\Rep(\Uhg)\ar[dd]^{f_\hbar}\\
&&\ar@{=>}[ll]&&\ar@{=}[ll]*+{\phantom{ciao}}\\
\kvect\ar@{=}[rr]&&\kvect\ar@{=}[rr]&&\kvect}\]
where $f_\hbar:\Rep(\DJ{\g})\to\kvect$ is the forgetful functor. The EK
isomorphism $\Psi\EEK$ is then given by the identifications
\[\wh{U\g[[\hbar]]}:=\sfEnd{f}\cong\sfEnd{F}=\sfEnd{f_\hbar\circ\wt{F}}\cong\sfEnd{f_\hbar}=:\wh{\Uhg}\]

\subsection{} % what we do
%--------------

Overlaying the above diagrams for an inclusion $\iD:\gD\hookrightarrow\g$
of Manin triples shows that constructing an isomorphism $\wh{\Uhg}\stackrel
{\sim}{\longrightarrow}\wh{U\g[[\hbar]]}$ compatible with $\iD$ may be achieved
by filling in the diagram
\[\xymatrix@C=0.1in@R=0.3in{
			&\D_\Phi(\g)\ar[dl]_{i_D^*}\ar@{=}[rr]\ar'[d][dd]^{f}&&\D_\Phi(\g)\ar@{-->}[dl]_{\Gamma}\ar@{->}[rr]^(.55){\wt{F}}\ar'[d][dd]^{F}&&\Rep(\Uhg)\ar[dl]_{i_{D,\hbar}^*}\ar[dd]^(.72){f_\hbar}\\
\D_\Phi(\gD)\ar@{=}[rr]\ar[dr]^{f_D}&&\D_\Phi(\gD)\ar@{->}[rr]^(.65){\wt{F}_D}\ar[dr]^{F_D}&&\Rep(\UhgD)\ar[dr]^{f_{D,\hbar}}&\\
			&\kvect\ar@{=}[rr]&&\kvect\ar@{=}[rr]&&\kvect}\]
where $f_D,f_{D,\hbar}$ are forgetful functors, $F_D$ the EK functor
for $\gD$, and $i_{D,\hbar}:\DJ{\gD}\to\DJ{\g}$ is the inclusion derived
from the functoriality of the quantization. %of \g_-

To do so, we first construct a {\it relative} fiber functor, that is a (tensor)
functor $\Gamma$ on $\D_\Phi(\g)$ whose target category is $\D_\Phi
(\gD)$ rather than $\kvect$, and which is isomorphic as abelian functor
to the restriction $\iD^*$. We then show the existence of a natural transformation
between the composition $\wt{F}_D\circ \Gamma$ and $i_{D,\hbar}^*
\circ\wt{F}$. Our constructions do not immediately yield a commutative
diagram, \ie the two factorizations $F\cong F_D\circ\Gamma$ deduced
from $f=f_D\circ\iD^*$ and $f_h=f_{D,\hbar}\circ i_{D,\hbar}^*$ do not
coincide, but this can easily be adjusted by using a different identification
$F\cong f$, which amounts to modifying the original EK isomorphism.

\subsection{}
%--------------

The construction of the functor $\Gamma$ is very much inspired by
\cite{ek-1}. The principle adopted by Etingof and Kazhdan is the following.
In a $\sfk$-linear monoidal category $\C$, a coalgebra structure on an
object $C\in\obj(\C)$ induces a tensor structure on the Yoneda functor 
\begin{equation*}
h_C=\Hom_{\C}(C, -):\C\to\vect_{\sfk}
\end{equation*}
If $\C$ is braided and $C_1,C_2$ are coalgebra objects in $\C$, then so
is $C_1\ten C_2$, and there is therefore a canonical tensor structure on
$h_{C_1\ten C_2}$.

If $\g$ is finite--dimensional, the polarization $\Ue{\g}\simeq\Mm\ten\Mp$,
where $M_\pm$ are the Verma modules $\ind_{\g_\mp}^\g\sfk$, realizes
$U\g$ as the tensor product of two coalgebra objects in $\DC{\Ueh{\g}}
{\Phi}$. This yields a tensor structure on the forgetful functor
\[h_{U\g}:\DC{U\g}{\Phi}[[\hbar]]\to\kvect\]
%and a definition of $\qEK{\g}$ as the Hopf algebra $\End(h_{M_-\otimes M_+})$ \cite{ek-1}.

Our starting point is to apply the same principle to the (abelian) restriction
functor $i_D^*:\DC{U{\g}}{\Phi}\to\DC{\Ueh{\g_D}}{\Phi}$. We therefore
factorize $\Ue{\g}$ as a tensor product of two coalgebra objects $\Lm,\Np$
in the braided monoidal category of $(\g,\gD)$--bimodules, with associator
$(\Phi \cdot \Phi_D^{-1})$, where $\Phi_D^{-1}$ acts on the right. Just as the
modules $\Mm,\Mp$ are related to the decomposition $\g=\gm\oplus\gp$,
$\Lm$ and $\Np$ are related to the asymmetric decomposition
$$\g=\Lmm\oplus\Lpp$$
where $\Lmm=\gm\cap\gD^{\perp}$ and $\Lpp=\gD\oplus\Lmp$. This factorization
induces a tensor structure on the functor $\Gamma=h_{\Lm\ten\Np}$, canonically
isomorphic to $i_D^*$ through the right $\gD$--action on $\Np$. As
in \cite[Part II]{ek-1}, this tensor structure can also be defined in the
infinite--dimensional case.

\subsection{}
%--------------

To construct a natural transformation making the following diagram
commute
\[\xymatrix{
\D_\Phi(\g)\ar[rr]^{\wt{F}}\ar[dd]_{\Gamma}	&&\Rep(\DJ{\g})\ar@{=>}[ddll]\ar[dd]^{(\iD)^*_\hbar}\\														
&&& \\
\D_\Phi(\g_D)\ar[rr]_{\wt{F}_D}			&&\Rep(\DJ{\gD})
}\]
we remark, as suggested to us by P. Etingof, that a quantum analogue
{$\Gamma_\hbar$ of} $\Gamma$ can be similarly {defined} using a
quantum version $\Lm^{\hbar},\Np^{\hbar}$
of the modules $\Lm,\Np$. The functor $\Gamma_\hbar=\Hom_{\qEK\g}(L_-^\hbar\otimes
N_+^\hbar,-)$ is naturally isomorphic to $(\iD)_\hbar^*$ as tensor functor,
since there is no associator involved on this side. Moreover, {an identification}
$$\wt{F}_D\circ\Gamma\simeq\Gamma_\hbar\circ\wt{F}$$
{is readily obtained, provided one establishes isomorphisms of $(\qEK{\g},\qEK\g_D)$--bimodules
\[\wt{F}_D\circ \wt{F}(\Lm)\simeq\Lm^{\hbar}\aand
\wt{F}_D\circ \wt{F}(\Np)\simeq\Np^{\hbar}\]}

\subsection{}
%--------------

While for $M_{\pm}$ it is easy to {construct} an isomorphism between $\wt{F}
(M_{\pm})$ and {the} quantum counterparts of {$M_\pm$}, the proof for $\Lm,
\Np$ is more involved. It relies on {a description of} the quantization functor $F
\oEK$ in terms of {\sc Prop} categories (cf. \cite{ek-2, EG}) and the realization
of $\Lm,\Np$ as universal objects in a suitable colored {\sc Prop} describing the
inclusion of bialgebras $\gD\subset\g$. This yields in particular a relative extension
of the EK functor with input a pair of Lie bialgebras $\a,\b$ which is {\it split}, \ie
endowed with maps  $\a\leftrightarrows^i_p\b$ such that $p\circ i=\id$.

\subsection{}
%--------------

Given that we work throughout with completions of algebras obtained as
endomorphisms of fiber functors, the transfer of structure from $\DJ{\g}$
to $\Ueh{\g}$ is more conveniently phrased in terms of categories. Part of
this paper is therefore devoted to rephrasing the definition of \qcqtqba in
categorical terms. This yields the notion of a {\it quasi--Coxeter category},
which is to a generalized braid group $B$ what a braided tensor category
is to Artin's braid groups, and of a {\it \qc tensor category}. Interestingly perhaps,
both notions be concisely rephrased in terms of a 2--functor from a combinatorially
defined 2--category $\qC{D}$ to the 2--categories $\Cat,\Cat^\ten$ of categories
and tensor categories respectively. The objects of $\qC{D}$ are the subdiagrams
of the Dynkin diagram $D$ of $B$ and, for two subdiagrams $D'\subseteq
D''$, $\Hom_{\qC{D}}(D'',D')$ is the fundamental 1--groupoid of the \DCP
associahedron for the quotient diagram $D''/D'$ \cite{DCP2,vtl-4}.

%---------------------------
%OUTLINE
%---------------------------

\subsection{Outline of the paper}
%---------------------------------------

We begin in Section \ref{s:diagrams} by reviewing a number of combinatorial
notions which will be used in later sections. In Section \ref{s:qcqtqba} we define
\qc (tensor) categories. In Section \ref{s:ek}, we review the construction of the
\nEK quantization functor and the isomorphism $\Psi\oEK$ following \cite{ek-1,ek-6}.
In Section \ref{s:Gamma}, we modify this construction by using generalized
Verma modules $L_-,N_+$, and obtain a relative fiber functor $\Gamma:\D
_\Phi(\g)\to\D_{\Phi_D}(\gD)$. In Section \ref{s:Gammah}, we define the quantum
generalized Verma modules $L_-^\hbar$ and $N_+^\hbar$. Using suitably defined
{\sc Props} we then show, in Section \ref{s:props} that these are isomorphic to
the EK quantization of their classical counterparts. In Section
\ref{s:mchain}, we use these results to show that, for any given chain of Manin
triples ending in a given $\g$, there exists a quantization of $U\g$ which is
compatible with each inclusion and independent, up to isomorphism, of the
choice of the given chain. Finally, in Section \ref{s:qcstructure}, we apply
these results to the case of a \KM algebra $\g$ and obtain the desired tranport
of its \qcqtqba structure to the completion of $\Ueh{\g}$ \wrt category $\O$,
integrable modules.

\subsection{Acknowledgments}
%-------------------------------------

We are very grateful to Pavel Etingof for his interest in the present work and for many
enlightening discussions.

%Section \ref{s:ek}
%Section \ref{s:qcqtqba}
%Section \ref{s:Gamma}
%Section \ref{s:tensor-structure}
%Section \ref{s:Gammah}
%Section \ref{s:props}
%Section \ref{s:qcstructure}

%--------------------------------------------------------------------
% COMPLETION OF QCQTQBA
%---------------------------------------------------------------------

\section{Diagrams and nested sets}\label{s:diagrams}
%------------------------------------------

We review in this section a number of combinatorial notions associated to a diagram
$D$, in particular the definition of nested sets on $D$ and of the \DCP associahedron
of $D$ following \cite{DCP2} and \cite[Section 2]{vtl-4}.

\subsection{Nested sets on diagrams} 
%---------------------------------------------

%We review some basic definitions from \cite[Section 2.1]{vtl-4} and \cite{DCP1,DCP2}.
By a \emph{diagram} we shall mean a nonempty undirected graph $D$ with no multiple
edges or loops. We denote the set of vertices of $D$ by $\mathsf{V}(D)$ and set $|D|=
|\mathsf{V}(D)|$. A \emph{subdiagram} $B\subseteq D$ is a full subgraph of $D$, that
is, a graph consisting of a subset $\sfV(B)$ of vertices of $D$, together with all edges
of $D$ joining any two elements of $\sfV(B)$. We will often abusively identify such a
$B$ with its set of vertices and write $i\in B$ to mean $i\in\sfV(B)$. We denote by
$\Csd{D}$ the set of subdiagrams of $D$.

The union $B_1\cup B_2$ of two subdiagrams $B_1,B_2\subset D$ is the subdiagram
having $\sfV(B_1)\cup\sfV(B_2)$ as its set of vertices. Two subdiagrams $B_1,B_2\subset
D$ are \emph{orthogonal} if $\sfV(B_1)\cap\sfV(B_2)=\emptyset$ and no two vertices $i\in
B_1$, $j\in B_2$ are joined by an edge in $D$. $B_1$ and $B_2$ are \emph{compatible}
if either one contains the other or they are orthogonal.
 
\begin{definition}
A \emph{nested set} on a diagram $D$ is a collection $\H$ of pairwise compatible,
connected subdiagrams of $D$ which contains the connected components $D_1,
\dots, D_r$ of $D$.
\end{definition}

\subsection{The \DCP associahedron}
%----------------------------------------------

Let $\N_D$ be the partially ordered set of nested sets on $D$, ordered by reverse inclusion.
$\N_D$ has a unique maximal element ${\bf 1}=\{D_i\}$ and its minimal elements are the
\emph{maximal nested sets}. We denote the set of maximal nested sets on $D$ by $\Mns
{D}$. Every nested set $\H$ on $D$ is uniquely determined by  a collection $\{\H_i\}_{i=1}
^r$ of nested sets on the connected components of $D$. We therefore obtain canonical
identifications
\[\N_D=\prod_{i=1}^r \N_{D_i}\qquad\text{and}\qquad\Mns{D}=\prod_{i=1}^r\Mns{D_i}\]

The \emph{De Concini--Procesi associahedron} $\A_D$ is the regular CW--complex
whose poset of (nonempty) faces is $\N_D$. It easily follows from the definition that 
\[\A_D=\prod_{i=1}^r \A_{D_i}\]
It can be realized as a convex polytope of dimension $|D|-r$. For any $\H\in\N_D$,
we denote by $\dim(\H)$ the dimension of the corresponding face in $\A_D$.\\

\subsection{The rank function of $\N_D$}
%------------------------------------------------

For any nested set $\H$ on $D$ and $B\in\H$, we set
$i_{\H}(B)=\bigcup_{i=1}^m B_i$
where the $B_i$'s are the maximal elements of $\H$ properly contained in $B$.

\begin{definition}
Set $\ul{\alpha}_{\H}^B=B\setminus i_{\H}(B)$. We denote by
\[
n(B;\H)=|\ul{\alpha}_{\H}^B|\qquad\text{and}\qquad n(\H)=\sum_{B\in\H} \big(n(B;\H)-1\big)
\]
An element $B\in\H$ is called \emph{unsaturated} if $n(B;\H)>1$. 
\end{definition}

\begin{proposition}\return
\begin{itemize}
\item[(i)] For any nested set $\H\in\N_D$,
\[n(\H)=|D|-|\H|=\dim(\H)\]
\item[(ii)] If $\H$ is a maximal nested set if and only if $n(B;\H)=1$ for any $B\in\H$.
\item[(iii)] Any maximal nested set is of cardinality $|D|$.
\end{itemize}
\end{proposition}

For any $\F\in\Mns{D}$, $B\in\F$, $i_{\F}(B)$ denotes the maximal element in $\F$ properly contained in $B$ and $\ul{\alpha}_{\F}^B=B\setminus i_{\F}(B)$ consists of one vertex, denoted $\alpha_{\F}^B$. 

For any $\F\in\Mns{D}$, $B\in\F$, we denote by $\F_{B}\in\Mns{B}$ the maximal nested set induced by $\F$ on $B$.

\subsection{Quotient diagrams}
%----------------------------------------

Let $B\subsetneq D$ a proper subdiagram with connected components $B_1,\dots, B_m$.

\begin{definition}%QUOT
The set of vertices of the quotient diagram $D/B$ is $\sfV(D)\setminus\sfV(B)$.
Two vertices $i\neq j$ of $D/B$ are linked by an edge if and only if the following
holds in $D$
\[i\not\perp j\qquad\text{or}\qquad i,j\not\perp B_i\qquad\text{for some }i=1,\dots, m\]
\end{definition}

For any connected subdiagram $C\subseteq D$ not contained in $B$,
we denote by $\ol{C}\subseteq D/B$ the connected subdiagram with
vertex set $V(C)\setminus V(B)$. %We have the following

\subsection{Compatible subdiagrams of $D/B$}
%---------------------------------------------------------

\begin{lemma}\label{le:quotient compatibility}
Let $C_1,C_2\nsubseteq B$ be two connected subdiagrams of $D$
which are compatible. Then
\begin{enumerate}
\item $\ol{C}_1,\ol{C}_2$ are compatible unless $C_1\perp C_2$ and
$C_1,C_2\nperp B_i$ for some $i$.
\item If $C_1$ is compatible with every $B_i$, then $\ol{C}_1$ and $
\ol{C}_2$ are compatible.
\end{enumerate}
In particular, if $\F$ is a \ns on $D$ containing each $B_i$, then $\ol
{\F}=\{\ol{C}\}$, where $C$ runs over the elements of $\F$ such that
$C\nsubseteq B$, is a \ns on $D/B$.
\end{lemma}

Let now $A$ be a connected subdiagram of $D/B$ and denote by
$\wt{A}\subseteq D$ the connected sudbdiagram with vertex set
\begin{equation*}
V(\wt{A})=V(A)\bigcup_{i:B_i\nperp V(A)} V(B_i)
\end{equation*}
Clearly, $A_1\subseteq A_2$ or $A_1\perp A_2$ imply $\wt{A_1}
\subseteq\wt{A_2}$ and $\wt{A_1}\perp\wt{A_2}$ respectively, so
the lifting map $A\rightarrow\wt{A}$ preserves compatibility.\\

\subsection{Nested sets on quotients}
%---------------------------------------------

For any connected subdiagrams $A\subseteq D/B$ and $C\subseteq
D$, we have
\begin{equation*}
\ol{\wt{A}}=A
\aand\wt{\ol{C}}=C\bigcup_{i:B_i\nperp C} B_i
\end{equation*}
In particular, $\wt{\ol{C}}=C$ if, and only if, $C$ is compatible with $B_1,
\ldots,B_m$ and not contained in $B$. The applications $C\rightarrow
\ol{C}$ and $A\rightarrow\wt{A}$ therefore yield a bijection between the
connected subdiagrams of $D$ which are either orthogonal to or strictly
contain each $B_i$ and the connected subdiagrams of $D/B$. This bijection preserves compatibility and
therefore induces an embedding $\N_{D/B}\hookrightarrow\ND$. This
yields an embedding
\begin{equation*}
\N_{D/B}\times\N_{B}=\N_{D/B}\times\big(\N_{B_1}\odots{\times}\N_{B_m}\big)\hookrightarrow\ND
\end{equation*}
with image the poset of nested sets on $D$ containing each $B_i$. Similarly, for any $B\subseteq B'\subseteq B''$, we obtain a map
\[\cup: \N_{B''/B'}\times\N_{B'/B}\hookrightarrow\N_{B''/B}\]
The map $\cup$ restricts to maximal nested sets. For any $B\subset B'$, we denote by $\Mns{B',B}$ the collection of maximal nested sets on $B'/B$. Therefore, for any $B\subset B' \subset B''$, we obtain an embedding
\[\cup:\Mns{B'',B'}\times\Mns{B',B}\to\Mns{B'',B}\]
such that, for any $\F\in\Mns{B'',B'},\G\in\Mns{B',B}$,
\[(\F\cup\G)_{B'/B}=\G\]

\subsection{Elementary and equivalent pairs}
%-----------------------------------------------------------

\begin{definition}%ELPAIR
An ordered pair $(\G,\F)$ in $\Mns{D}$ is called \emph{elementary} if $\G$ and $\F$ differ by one element. A sequence $\H_1,\dots,\H_m$ in $\Mns{D}$ is called \emph{elementary} if $|\H_{i+1}\setminus\H_i|=1$ for any $i=1,2,\dots, m-1$.
\end{definition}
\vspace{0.05in}
\begin{definition}%SUPP
The \emph{support} $\supp(\F,\G)$ of an elementary pair in $\Mns{D}$ is the unique unsaturated element of $\F\cap\G$. The \emph{central support} $\z\supp(\F,\G)$ is the union of the maximal elements of $\F\cap\G$ properly contained in $\supp(\F,\G)$. Thus
\[\z\supp(\F,\G)=\supp(\F,\G)\setminus\alpha_{\F\cap\G}^{\supp(\F,\G)}\]
\end{definition}
\vspace{0.05in}
\begin{definition}%EQUIPAIR
Two elementary pairs $(\F,\G),(\F',\G')$ in $\Mns{D}$ are \emph{equivalent} if
\[
\begin{array}{c}
\supp(\F,\G)=\supp(\F',\G')\\[1.1ex]
\alpha^{\supp(\F,\G)}_{\F}=\alpha_{\F'}^{\supp(\F',\G')}\qquad\alpha^{\supp(\F,\G)}_{\G}=\alpha_{\G'}^{\supp(\F',\G')}
\end{array}
\]
\end{definition}
\vspace{0.05in}

%------------------------------------------------------------
% D--CATEGORY
%------------------------------------------------------------

\section{Quasi--Coxeter categories}\label{s:qcqtqba}
%===========================

The goal of this section is to rephrase the notion of \qcqtqba defined in \cite{vtl-4}
in terms of terms of categories of representations.

\Omit{The goal of this section is to adapt the definition of a quasi--Coxeter quasitriangular
quasibialgebra, introduced in \cite[Section 3]{vtl-4}, to completion of an algebra with
respect to a subcategory of its modules (\eg category $\O$), compatible with restriction.
This description amounts to the definition of a categorical notion of quasi--Coxeter
quasitriangular quasibialgebras.}

\subsection{Algebras arising from fiber functors}\label{ss:endiso}
%-----------------------------------------------------------

We shall repeatedly need the following elementary

\begin{lemma}
Consider the following situation
\begin{equation*}
\xymatrix{\C \ar@/_/[drr]_{F}="F" \ar[rr]^{H}& & \D \ar[d]^{G} \ar@{=>}"F"_{\alpha} \\ & & \A }
\end{equation*}
where $\A,\C,\D$ are additive $\sfk$--linear categories, $F,G,H$ functors, and $\alpha$ is an
invertible transformation. If $H$ is an equivalence of categories, the map $\mathsf{End}(G)
\longrightarrow\mathsf{End}(F)$ given by
\[\{g_W\}\mapsto\{{\sf Ad}(\alpha_V^{-1})(g_{H(V)})\}\]
is an algebra isomorphism.
\end{lemma}

\subsection{$D$--categories} 
%----------------------------------

Recall \cite[Section 3]{vtl-4} that, given a diagram $D$, a $D$--algebra is a pair
$(A,\{A_B\}_{B\in\Csd{D}})$, where $A$ is an associative algebra and $\{A_B\}_
{B\in\Csd{D}}$ is a collection of subalgebras indexed by $\Csd{D}$ and satisfying 
\begin{equation*}
A_B\subseteq A_{B'}\quad\mbox{ if }B\subseteq B'
\qquad\text{and}\qquad
[A_B,A_{B'}]=0\quad\mbox{ if }B\perp B'
\end{equation*}

The following rephrases the notion of $D$--algebras in terms of their category of
representations.

\begin{definition}
A \emph{$D$--category}
\[\C=(\{\C_{B}\}, \{F_{BB'}\})\]
is the datum of 
\begin{itemize}
\item a collection of $\sfk$--linear \abelian categories $\{\C_{B}\}_{B\subseteq{D}}$
\item for any pair of subdiagrams $B\subseteq B'$, an \exact $\sfk$--linear functor
$F_{BB'}:\C_{B'}\to\C_{B}$\footnote{When $B=\emptyset$ we will omit the index $B$.}
\item for any $B\subset B'$, $B'\perp B''$, $B',B''\subset B'''$, a homomorphism of 
$\sfk$--algebras
\[
\eta:\sfEnd{F_{BB'}}\to\sfEnd{F_{(B\cup B'')B'''}}
\]
\end{itemize}
satisfying the following properties
\begin{itemize}
\item For any $B\subseteq D$, $F_{BB}=\id_{\C_B}$.
\item For any $B\subseteq B'\subseteq B''$, $F_{BB'}\circ F_{B'B''}=F_{BB''}$.
\item For any $B\subset B'$, $B'\perp B''$, $B',B''\subset B'''$, the following diagram
of algebra homomorphisms commutes:
%\[
%\xymatrix{
%\sfEnd{F_{BB'}}\ar[r]^{\id\ten\id_{F_{B'B'''}}} \ar[d]^{\id_{F_{B(B\cup B'')}}\ten\eta} & \sfEnd{F_{BB'}}\ten\sfEnd{F_{B'B'''}}\ar[d]^{\circ}\\
%\sfEnd{F_{B(B\cup B'')}}\ten\sfEnd{F_{(B\cup B'')B'''}} \ar[r]^{\circ} & \sfEnd{F_{BB'''}}
%}
%\]
\[
\xymatrix@C=.1cm@R=.3cm{
& \sfEnd{F_{BB'}}\ar[rd]^{\id\ten\id_{F_{B'B'''}}} \ar[dl]_{\id_{F_{B(B\cup B'')}}\ten\eta} &\\
 \sfEnd{F_{BB'}}\ten\sfEnd{F_{B'B'''}}\ar[dr]_{\circ} & & \sfEnd{F_{B(B\cup B'')}}\ten\sfEnd{F_{(B\cup B'')B'''}} \ar[ld]^{\circ}\\
&  \sfEnd{F_{BB'''}} & 
}
\]

\end{itemize}
\end{definition}

\begin{rem}
It may seem more natural to replace the equality of functors $F_{BB'}\circ F_{B'B''}
=F_{BB''}$ by the existence of invertible natural transformations $\alpha_{BB''}^{B'}:
F_{BB'}\circ
F_{B'B''}\Rightarrow F_{BB''}$ for any $B\subseteq B'$ satisfying the associativity
constraints $\alpha_{BB'''}^{B'}\circ F_{BB'}(\alpha_{B'B'''}^{B''})=\alpha_{BB'''}^
{B''}\circ(\alpha_{BB''}^{B'})_{F_{B''B'''}}$ for any $B\subseteq B'\subseteq B''
\subseteq B'''$. A simple coherence argument shows however that this leads
to a notion of $D$--category which is equivalent to the one given above.
\end{rem}

\begin{rem} We will usually think of $\C_{\emptyset}$ as a base category and at
the functors $F$ as forgetful functors. Then the family of algebras $\sfEnd{F_B}$
defines, through the morphisms $\alpha$, a structure of $D$--algebra on $\sfEnd
{F_D}$. Conversely, every $D$--algebra $A$ admits such a description setting
$\C_B=\Rep A_B$ for $B\neq\emptyset$ and $\C_\emptyset=\vect_k$, $F_{BB'}
=i_{B'B}^*$, where $i_{B'B}:A_B\subset A_{B'}$ is the inclusion.
\end{rem}

\begin{rem}
The conditions satisfied by the maps $\eta$ imply that, given $B=\bigsqcup_{j=1}^rB_j$, with $B_j\in\Csd{D}$ pairwise orthogonal, the images in $\sfEnd{F_B}$ of the maps
\[\sfEnd{F_{B_j}}\to\sfEnd{F_{B_j}F_{B_jB}}=\sfEnd{F_B}\]
pairwise commute. This condition rephrases for the endomorphism 
algebras the $D$--algebra axiom 
\[[A_{B'},A_{B''}]=0 \qquad\forall\quad B'\perp B''\]
that is equivalent to the condition, for any $B\supset B',B''$, 
\[A_{B'}\subset A^{B''}_B\]
\end{rem}

\begin{rem} The above definition of $D$--category may be rephrased as follows. Let $\ID$ be the
category whose objects are subdiagrams $B\subseteq D$ and morphisms $B'\to
B$ the inclusions $B\subset B'$. Then a $D$--category is a functor \[\C: \ID\to\Cat\]
\end{rem}

%-------------------------------
% D-MORPHISMS
%-------------------------------

\subsection{Strict morphisms of $D$--categories}
%------------------------------------------------------------

The interpretation of $D$--categories in terms of $\ID$ suggests that a morphism of
$D$--categories $\C,\C'$ is one of the corresponding functors
\[
\xymatrix{\ID \ar@/^15pt/@<2pt>[r]^{\C}="C" \ar@/_15pt/@<-2pt>[r]_{\C'}="C'" \ar@{=>}"C";"C'"& \Cat}
\]
This yields the following definition. For simplicity, we assume that $\C_{\emptyset}=
\C'_{\emptyset}$.

\begin{definition}
A \emph{strict morphism} of $D$--categories $\C,\C'$ is the datum of 
\begin{itemize}
\item for any $B\subseteq D$, a functor $H_{B}:\C_B\to\C'_B$ 
\item for any $B\subseteq B'$, a natural transformation
\begin{equation}\label{eq:B B'}
\xymatrix@C=0.6in@R=0.6in{
\C_{B'} \ar[r]^{H_{B'}} \ar[d]_{F_{BB'}} & \C'_{B'}\ar[d]^{F'_{BB'}} \ar@{=>}[dl]_{\gamma_{BB'}} \\ 
\C_{B}\ar[r]_{H_{B}} & \C'_{B}}
\end{equation}
\end{itemize}
such that 
\begin{itemize}
\item $H_\emptyset=\id$
\item $\gamma_{BB}=\id_{H_B}$
\item For any $B\subseteq B'\subseteq B''$,
\[\gamma_{BB''}=\gamma_{BB'}\circ\gamma_{B'B''}\] 
where $\circ$ is the composition of natural transformations defined by
\begin{equation}\label{eq:B B' B''}
\xymatrix{\C_{B''} \ar[r]\ar[d]& \C'_{B''}\ar@{=>}[dl]\ar[d] \ar[d]\\ \C_{B'}\ar[r]\ar[d] & \C'_{B'}\ar@{=>}[dl]\ar[d]\\ \C_{B}\ar[r] & \C'_{B}}
\end{equation}
\end{itemize}
\end{definition}

\noindent The diagram \eqref{eq:B B'}, with $B=\emptyset$, induces an algebra homomorphism
$\sfEnd{F'_{B'}} \to \sfEnd{F_{B'}}$ which, by \eqref{eq:B B' B''} is compatible with the maps $\sfEnd
{F_B}\to\sfEnd{F_{B'}}$ and $\sfEnd{F'_B}\to\sfEnd{F'_{B'}}$ for any $B\subset B'$. As pointed out
in \cite[3.3]{vtl-4}, this condition is too restrictive and will be weakened in the next paragraph.

%-------------------
%MORPHISM
%-------------------

\subsection{Morphisms of $D$--categories}
%-----------------------------------------------------

\begin{definition}
A \emph{morphism} of $D$--categories $\C,\C'$, with $\C_{\emptyset}=\C'_{\emptyset}$, is the datum
of
\begin{itemize}
\item for any $B\subseteq D$ a functor $H_{B}:\C_B\to\C'_B$
\item for any $B\subseteq B'$ and $\F\in\Mns{B,B'}$, a natural transformation
\begin{equation*}
\xymatrix@C=0.6in@R=0.6in{
\C_{B'} \ar[r]^{H_{B'}} \ar[d]_{F_{BB'}} & \C'_{B'}\ar[d]^{F'_{BB'}} \ar@{=>}[dl]_{\gamma^\F_{BB'}} \\ 
\C_{B}\ar[r]_{H_{B}} & \C'_{B}}
\end{equation*}
\end{itemize}
such that
\begin{itemize}
\item $H_\emptyset=\id$
\item $\gamma^\F_{BB}=\id_{H_B}$
\item for any $B\subseteq B'\subseteq B''$, $\F\in\Mns{B,B'}$, $\G\in\Mns{B',B''}$,
\[\gamma_{BB'}^{\F}\circ\gamma_{B'B''}^{\G}=\gamma_{BB''}^{\F\vee\G}\]
\end{itemize}
\end{definition}

\begin{rem} For any $\F\in\Mns{B'}$, the natural transformation $\gamma^{\F}_{B'}$
induces an algebra homomorphism $\Psi^{\F}_{B'}:\sfEnd{F'_{B'}}\to\sfEnd{F_{B'}}$
such that the following diagram commutes for any $B\in\F$
\[
\xymatrix{\sfEnd{F'_{B'}} \ar[r]^{\Psi_{B'}^{\F}}& \sfEnd{F_{B'}}\\ \sfEnd{F'_B}\ar[r]^{\Psi_B^{\F_B}}\ar[u]& \sfEnd{F_B}\ar[u]}
\]
In particular, the collection of homomorphisms $\{\Psi^{\F}_D\}$ defines a morphism of
$D$--algebras $\sfEnd{F'_D}\to\sfEnd{F_D}$ in the sense of \cite[3.4]{vtl-4}.
\end{rem}

\begin{rem} The above definition may be rephrased as follows. Let $\MD$ be the category
with objects the subdiagrams $B\subseteq D$ and
morphisms $\Hom(B',B)=\Mns{B',B}$, with composition given by union. There is a
forgetful functor $\MD\to\ID$ which is the identity on objects and maps $\F\in\Mns
{B',B}$ to the inclusion $B\subseteq B'$. Given two $D$--categories $\C,\C': \ID\to\mathsf
{Cat}$ a morphism $\C\to\C'$ as defined above coincides with a morphism of the functors
$\MD\to\Cat$ given by the composition
\[\xymatrix{\MD\ar[r] & \ID \ar@<+3pt>[r]^{\C} \ar@<-3pt>[r]_{\C'} & \Cat}\]
\end{rem}

% labelling, and braid groups, v symbol

\subsection{Quasi--Coxeter categories}
%-----------------------------------------------

\begin{definition}%LABELLING 
A {\it labelling} of the diagram $D$ is the assignment of an integer $m_{ij}\in
\{2,3,\dots, \infty\}$ to any pair $i,j$ of distinct vertices of $D$ such that
\[m_{ij}=m_{ji}\qquad m_{ij}=2\]
if and only if $i\perp j$.
\end{definition}

Let $D$ be a labeled diagram.

\begin{definition}% BRAID GROUP 
The \emph{Artin braid group} $B_D$ is the group generated by elements $S_i$ 
labeled by the vertices $i\in D$ with relations
\[\underbrace{S_iS_j\cdots}_{m_{ij}}=\underbrace{S_jS_i\cdots}_{m_{ij}}\]
for any $i\neq j$ such that $m_{ij}<\infty$. We shall also refer to $B_D$ as the 
braid group corresponding to $D$.
\end{definition}

\begin{definition}
A \emph{quasi--Coxeter category of type $D$}
\[\C=(\{\C_B\},\{F_{BB'}\},\{\Phi_{\F\G}\}, \{S_i\})\] is the datum of
\begin{itemize}
\item a $D$--category $\C=(\{\C_B\},\{F_{BB'}\})$
\item for any elementary pair $(\F,\G)$ in $\Mns{B,B'}$, a natural transformation
\[\Phi_{\F\G}\in\sfAut{F_{BB'}}\]
\item for any vertex $i\in \mathsf{V}(D)$, an element
\[S_{i}\in\sfAut{F_i}\]
\end{itemize}
satisfying the following conditions
\begin{itemize}
\item {\bf Orientation.} For any elementary pair $(\F,\G)$,
\begin{equation*}
\Phi_{\G\F}=\Phi_{\F\G}^{-1}
\end{equation*}
\item {\bf Coherence.} For any elementary sequences $\H_1,\dots,\H_m$ and $\K_1,\dots, \K_l$ in $\Mns{B,B'}$ such that $\H_1=\K_1$ and $\H_m=\K_l$,
\begin{equation*}
\Phi_{\H_{m-1}\H_{m}}\cdots \Phi_{\H_1\H_2}=\Phi_{\K_{l-1}\K_l}\cdots\Phi_{\K_1\K_2}
\end{equation*}
\item {\bf Factorization.} The assignment
\[\Phi:\Mns{B,B'}^2\to{\sf Aut}(F_{B'B}) \]
is compatible with the embedding
\[\cup:\Mns{B,B'}\times\Mns{B',B''}\to\Mns{B,B''}\]
for any $B''\subset B'\subset B$, \ie the diagram
\[
\xymatrix@C=.8in{
\Mns{B,B'}^2\times\Mns{B',B''}^2 \ar[d]_{\cup}\ar[r]^{\Phi\times\Phi} & \sfAut{F_{B''B'}}\times\sfAut{F_{B'B}} \ar[d]^{\circ}\\
\Mns{B,B''}^2\ar[r]^{\Phi} & \sfAut{F_{B''B}}
}
\]
is commutative.
\item {\bf Braid relations.} For any pairs $i,j$ of distinct vertices of $B$, such that $2< m_{ij}<\infty$, and elementary pair $(\F,\G)$ in $\Mns{B}$ such that $i\in\F,j\in\G$, the following relations hold in $\sfEnd{F_B}$
\begin{equation*}
\sfAd{\Phi_{\G\F}}(S_i)\cdot S_j \cdots = S_j\cdot\sfAd{\Phi_{\G\F}}(S_i)\cdots
\end{equation*}
where, by abuse of notation, we denote by $S_i$ its image in $\sfEnd{F_B}$
and the number of factors in each side equals $m_{ij}$.
\end{itemize}
\end{definition}

The elements $S_i$ will be commonly referred at as \emph{local monodromies}.

\begin{rem} It is clear that the factorization property implies the support and forgetful properties
as stated in \cite[Def. 3.12]{vtl-4}.
\begin{itemize}
\item {\bf Support.} For any elementary pair $(\F,\G)$ in $\Mns{B,B'}$, let $S=\supp(\F,\G), Z=\z\supp(\F,\G)\subseteq D$ and 
\[\wt{\F}=\F|^{\supp(\F,\G)}_{\z\supp(\F,\G)}\qquad \wt{\G}=\G|^{\supp(\F,\G)}_{\z\supp(\F,\G)}\] 
Then
\begin{equation*}
\Phi_{\F\G}=\id_{BZ}\circ\Phi_{\wt{F}\wt{\G}}\circ\id_{B'S}
\end{equation*}
where the expression above denotes the composition of natural transformations
\begin{equation*}
\xymatrix{
\C_{B'} \ar@/_15pt/[dddd]_{F_{BB'}}="F"  \ar@/^15pt/[dddd]^{F_{BB'}}="G" \ar@{=>}"F";"G"^{\Phi_{\G\F}}  &  &\C_{B'} \ar[d]^{F_{SB'}}  
\\
 &&\C_{S} \ar@/_15pt/[dd]_{F_{ZS}}="tF" \ar@/^15pt/[dd]^{F_{ZS}}="tG" \ar@{=>}"tF";"tG"^{\Phi_{\wt{\G}\wt{\F}}} \\
 &=&\\
  &&\C_Z \ar[d]^{F_{BZ}}\\
 \C_{B} &&\C_{B}
 }
\end{equation*}
\item {\bf Forgetfulness.} For any equivalent elementary pairs $(\F,\G),(\F',\G')$ in $\Mns{B,B'}$
\begin{equation*}
\Phi_{\F\G}=\Phi_{\F'\G'}
\end{equation*}
\end{itemize}
\end{rem}

\begin{rem}
To rephrase the above definition, consider the 2--category $\qC{D}$ obtained by
adding to $\MD$ a unique 2--isomorphism $\varphi^{BB'}_{\G\F}:\F\to\G$ for any
pair of 1--morphisms $\F,\G\in\Mns{B',B}$, with the compositions
\[\varphi^{BB'}_{\H\G}\circ\varphi^{BB'}_{\G\F}=\varphi^{BB'}_{\H\F}
\qquad\text{and}\qquad
\varphi^{BB'}_{\F_2\G_2}\circ\varphi^{B'B''}_{\F_1\G_1}=\varphi^{BB''}_{\F_2\cup\F_1\,\G_2\cup\G_1}\]
where $\F,\G,\H\in\Mns{B',B}$, $B\subset B'\subseteq B''$ and $\F_1,\G_1\in\Mns{B'',B'}$, $\F_2,\G_2
\in\Mns{B',B}$. There is a unique functor $\qC{D}\to\ID$ extending $\MD\to\ID$, and
a quasi--Coxeter category is the same as a 2--functor $\qC{D}\to\Cat$ fitting in a diagram
\[\xymatrix{\qC{D}\ar[rr]^(.53){}="A" \ar[dr] & &  \Cat\\ & \ID \ar@{<=}"A" \ar[ur] & }\] 

Note that, for any $B\subset B'$, the category $\Hom_{\qC{D}}(B',B)$ is the
1--groupoid of the \DCP associahedron on $B'/B$ \cite{vtl-4}.
\end{rem}

\subsection{Morphisms of quasi--Coxeter categories}
%--------------

\begin{definition}
A {\it morphism of quasi--Coxeter categories} $\C,\C'$ of type $D$ is a morphism
$(H,\gamma)$ of the underlying $D$--categories such that
\begin{itemize}
\item For any $i\in B$, the corresponding morphism $\Psi_i:\sfEnd{F_i'}\to\sfEnd{F_i}$ satisfies
\[\Psi_i(S_i')=S_i\]
\item For any elementary pair $(\F,\G)$ in $\Mns{B,B'}$, 
\[
H_{B}(\Phi_{\F\G})\circ\gamma^{\F}_{BB'}\circ (\Phi'_{\G\F})_{H_{B'}}=\gamma^{\G}_{BB'}
\] 
in $\Nat{F'_{BB'}\circ H_{B'}}{H_B\circ F_{BB'}}$, as in the diagram %DA AGGIUSTARE
\[
\xymatrix{
 & & & \C'_{B'} \ar@{-->}@/_15pt/[dd]_{}="F" \ar@/^15pt/[dd]^{}="G" \ar@{-->}"F";"G"^{\Phi'_{\F\G}}  \ar@{<-}[dlll]_{H_{B'}}\\
\C_{B'} \ar@/_15pt/[dd]_{}="F'" \ar@/^15pt/[dd]^{}="G'" \ar@{=>}"F'";"G'"^{\Phi_{\F\G}} & & \ar@/^10pt/@<5pt>@{=>}[dl]^{\gamma^{\G}} \ar@/_10pt/@{-->}@<-5pt>[dl]_{\gamma^{\F}} & \\
& & & \C'_{B}\ar@{<-}[dlll]^{H_{B}}\\
\C_{B} & & & 
}
\]
\end{itemize}
\end{definition}

\begin{rem}
Note that the above condition can be alternatively stated in terms of morphisms $\Psi_F$
as the identity
\[
\Psi_{\G}\circ\sfAd{\Phi_{\G\F}}=\sfAd{\Phi'_{\G\F}}\circ\Psi_{\F}
\]
\end{rem}

\subsection{Strict $D$--monoidal categories}
%------------------------------------------------------

\begin{definition}
A \emph{strict $D$--monoidal category} $\C=(\{\C_B\}, \{F_{BB'}\}, \{J_{BB'}\}\})$
is a $D$--category $\C=(\{\C_B\}, \{F_{BB'}\}\})$ where
\begin{itemize}
\item for any $B\subseteq D$, $(\C_B,\ten_B)$ is a strict monoidal category
\item for any $B\subseteq B'$, the functor $F_{BB'}$ is endowed with a tensor
structure $J_{BB'}$
\end{itemize}
with the additional condition that, for every $B\subseteq B'\subseteq B''$,
$J_{BB'}\circ J_{B'B''}=J_{BB''}$.
\end{definition}

\begin{rem} The tensor structure $J^B$ induces on $\sfEnd{F_B}$ a coproduct
$\Delta_B: \sfEnd{F_B}\to\sfEnd{F_B^2}$, where $F_B^2:=\ten\circ(F_B\boxtimes F_B)$,
given by
\[\{g_V\}_{V\in\C_{B}}\mapsto\{\Delta_B(g)_{VW}:=\sfAd{J^B_{VW})(g_{V\ten W}}\}_{V,W\in\C_B}\]
Moreover, for any $B\subseteq B'$, $\sfEnd{F_{B}}$ is a subbialgebra of $\sfEnd{F_{B'}}$,
\ie the following diagram is commutative
\begin{equation*}
\xymatrix{
\sfEnd{F_B}\ar[d]\ar[r]^{\Delta_B} & \sfEnd{F_{B}^2} \ar[d]\\ 
\sfEnd{F_{B'}} \ar[r]_{\Delta_{B'}} & \sfEnd{F_{B'}^2}}
\end{equation*}
\end{rem}

\begin{rem} Note that a strict $D$--monoidal category can be thought of as functor
\[\C:\ID\to\Cat_0^{\ten}\] 
where $\Cat_0^{\ten}$ denotes the $2$--category of strict monoidal category, with
monoidal functors and gauge transformations.
\end{rem}

\begin{definition}
A morphism of strict $D$--monoidal categories is a natural transformation of the
corresponding $2$--functors $\MD\to\Cat_0^{\ten}$, obtained by composition with
$\MD\to\ID$.
\end{definition}

\subsection{$D$--monoidal categories} % (TRUE) D--quasibialgebra
%-----------------------------------------------

\begin{definition}\label{def:D-qba}
A \emph{$D$--monoidal category}
\[\C=(\{(\C_B,\ten_B,\Phi_B)\}, \{F_{BB'}\}, \{J^{\F}_{BB'}\})\]
is the datum of
\begin{itemize}
\item A $D$--category $(\{(\C_B\}, \{F_{BB'}\})$ such that each $(\C_B,\ten_B,\Phi_B)$
is a tensor category, with $\C_{\emptyset}$ a strict tensor category, \ie $\Phi_{\emptyset}
=\id$.
\item for any pair $B\subseteq B'$ and $\F\in\Mns{B,B'}$, a tensor structure $J^{BB'}_{\F}$
on the functor $F_{BB'}:\C_{B'}\to\C_{B}$
\end{itemize}
with the additional condition that, for any $B\subseteq B'\subseteq B''$, $\F\in\Mns{B'',B'}$,
$\G\in\Mns{B',B}$,
\[J_{BB'}^\G\circ J_{B'B''}^\F=J_{BB''}^{\F\cup\G}\]
\end{definition}

\begin{rem} The usual comparison with the algebra of endomorphisms leads to a collection
of bialgebras $(\sfEnd{F_B}, \Delta_{\F}, \varepsilon)$ endowed with multiple coproducts,
indexed by $\Mns{B}$.
\end{rem}

\begin{rem} A $D$--monoidal category can be thought of as a functor $\MD\to\Cat^\ten$
fitting in a diagram
\[
\xymatrix{\MD \ar[d]\ar[r]& \Cat^{\ten}\ar[d] \ar@{=>}[dl]\\ \ID \ar[r] & \Cat}
\]
Accordingly, a morphism of monoidal $D$--categories is one of the corresponding functors.
\[
\xymatrix{\MD\ar@<2pt>@/^15pt/[r]^{\C}="C" \ar@<-2pt>@/_15pt/[r]_{\C'}="C'" \ar@{=>}"C";"C'" & \Cat^{\ten}}
\] 
\end{rem}

\subsection{Fibered $D$--monoidal categories}
%---------------------------------------------------------

We shall often be concerned with $D$--monoidal categories such that
the underlying categories $(\C_B,\otimes_B)$ are strict, and the functors $F
_{BB'}:(\C_{B'},\otimes_{B'})\to(\C_B,\otimes_B)$ are tensor functors. This
may be described in terms
of the category $\MD$ as follows. Let $\Drin$ be the 2--category of Drinfeld
categories, that is strict tensor categories $(\C,\otimes)$ endowed with an
additional associativity constraint $\Phi$ making $(\C,\otimes,\Phi)$ a monoidal
category. There is a canonical forgetful 2--functor $\Drin\to\Cat^\ten_0$.\\

We shall say that a $D$--monoidal category {\it fibers over} a strict
$D$--monoidal category if the corresponding functor $\MD\to\Cat^\ten$
maps into $\Drin$ and fits in a commutative diagram
\[
\xymatrix{\MD \ar[d]\ar[r]& \Drin\ar[d] \ar@{=>}[dl]\\ \ID \ar[r] & \Cat^\ten_0}
\]
In this case, the coproduct $\Delta_\F$ on a bialgebra $\sfEnd{F_B}$ is the
twist of a reference coassociative coproduct $\Delta_0$ on $\sfEnd{F_D}$
such that $\Delta_0:\sfEnd{F_B}\to\sfEnd{F_B^2}$.

\subsection{Braided $D$--monoidal categories} % (TRUE) D-qtqba
%----------------------------------------------------------

\begin{definition}
A \emph{braided $D$--monoidal category}
\[\C=(\{(\C_B,\ten_B,\Phi_B,\beta_B)\}, \{(F_{BB'}, J_{BB'}^{\F}\})\] 
is the datum of
\begin{itemize}
\item a $D$--monoidal category $(\{(\C_B,\ten_B,\Phi_B)\},\{(F_{BB'}, J_{BB'}^{\F}\})$
\item for every $B\subseteq D$, a commutativity constraint $\beta_B$ in $\C_B$,
defining a braiding in $(\C_B,\ten_B,\Phi_B)$.
\end{itemize}
\end{definition}

\begin{rem} Note that the tensor functors $(F_{BB'},J_{BB'}^\F):\C_{B'}\to\C_B$ are
\emph{not} assumed to map the commutativity constraint $\beta_{B'}$ to $\beta_B$.
\end{rem}

\begin{definition}
A morphism of braided $D$--monoidal categories from $\C$ to $\C'$ is a morphism of
the underlying $D$--monoidal categories such that the functors $H_B:\C_B\to\C_B'$
are braided tensor functors.
\end{definition}

\begin{rem} The fact that $H_B$ are braided tensor functors automatically implies that
\[
\Psi_{\F}^{\ten 2}((R_B)_{J_{\F}})=(R'_{B})_{J'_{\F}}
\]
in analogy with \cite{vtl-4}, where $R_B=(12)\circ\beta_B$, and we are assuming that
$\C_\emptyset=\C'_\emptyset$ is a symmetric strict tensor category.
\end{rem}

\subsection{Quasi--Coxeter braided monoidal categories}
%----------------------------------------------------------------------

\begin{definition}\label{def:qc-cat}
A \emph{quasi--Coxeter braided monoidal category of type $D$}
\[\C=(\{(\C_B,\ten_B, \Phi_B,\beta_B)\}, \{(F_{BB'}, J^{\F}_{BB'})\},\{\Phi_{\F\G}\},\{S_i\})\]
is the datum of
\begin{itemize}
\item a quasi--Coxeter category of type $D$, 
\[\C=(\{\C_B\},\{F_{BB'}\}, \{\Phi_{\F\G}\}, \{S_i\})\]
\item a braided $D$--monoidal category 
\[\C=(\{(\C_B,\ten_B,\Phi_B,\beta_B)\}, \{(F_{BB'}, J^{BB'}_{\F})\})\]
\end{itemize}
satisfying the following conditions
\begin{itemize}
\item for any $B\subseteq B'$, and $\G,\F\in\Mns{B,B'}$, the natural transformation
$\Phi_{\F\G}\in\sfAut{F_{BB'}}$ determines an isomorphism of tensor functors
$(F_{BB'},J_{BB'}^\G)\to(F_{BB'},J_{BB'}^\F)$, that is for any $V,W\in\C_{B'}$,
\[(\Phi_{\G\F})_{V\ten W}\circ (J_{BB'}^{\F})_{V,W}=(J^{\G}_{BB'})_{V,W}\circ((\Phi_{\G\F})_V\ten(\Phi_{\G\F})_W)\]
\item for any $i\in D$, the following holds:
\begin{equation*}
\Delta_i(S_i)=(R_i)^{21}_{J_i}\cdot(S_i\ten S_i)
\end{equation*}
\end{itemize}
A morphism of quasi--Coxeter braided monoidal categories of type $D$ is a morphism
of the underlying quasi--Coxeter categories and braided $D$--monoidal categories.
\end{definition}

\begin{rem}
A quasi--Coxeter braided monoidal category of type $D$ determines a 2--functor
$\qC{D}\to\Cat^\ten$ fitting in a diagram
\[
\xymatrix{\qC{D} \ar[d]\ar[r]& \Cat^{\ten}\ar[d] \ar@{=>}[dl]\\ \ID \ar[r] & \Cat}
\]
Note however that this functor does not entirely capture the \qc braided monoidal
category since it does not encode the commutativity constraints $\beta_B$ and
automorphisms $S_i$.
\end{rem}

%TO CHECK
%==========
%\begin{rem}
%The element $S_i$ turns the pair $(\sfEnd{F_i}, S_i)$ in a topological half-ribbon
%algebra \cite{tingley}.
%\end{rem}

%----------------------------
%ETINGOF-KAZHDAN
%----------------------------
\section{Etingof-Kazhdan quantization}\label{s:ek}

We review in this section the results obtained in \cite{ek-1,ek-6}. More specifically, we follow the quantization
of Lie bialgebras given in \cite[Part II]{ek-1} and the case of generalized Kac--Moody algebras from \cite{ek-6}.

%---------------------------------------------------------------%
%TOPOLOGY
%---------------------------------------------------------------%

\subsection{Topological vector spaces}

The use of topological vector spaces is needed in order to deal with convergence
issues related to duals of infinite dimensional vector spaces and tensor product of such spaces.\\

Let $\sfk$ be a field of characteristic zero with the discrete topology and $V$ a topological vector space over $\sfk$. 
The topology on $V$ is {\em linear} if open subspaces in $V$ form a basis of neighborhoods of zero. 
Let $V$ be endowed with a linear topology and $p_V$ the natural map
\[p_V:V\to\lim(V/U)\]
where the limit is taken over the open subspaces $U\subseteq V$. Then $V$ is called \emph{separated} 
if $p_V$ is injective and \emph{complete} if $p_V$ is surjective. Throughout this section, we shall call 
\emph{topological vector space} a linear, complete, separated topological space.\\

If $U$ is an open subspace of a topological vector space $V$, then the quotient $V/U$ is discrete. It is then 
possible, given two topological vector spaces $V$ and $W$, to define the topological tensor product as
\[V\ctp W:= \lim V/V' \ten W/W'\]
where the limit is take over open subspaces of $V$ and $W$. We then denote by $\Hom_{\sfk}(V,W)$ the topological vector space
of continuous linear operators from $V$ to $W$ equipped with the weak topology. Namely, a basis of neighborhoods of zero in
$\Hom_{\sfk}(V,W)$ is given by the collection of sets
\[Y(v_1,\dots,v_n,W_1,\dots,W_n):=\{f\in \Hom_{\sfk}(V,W)\;|\: f(v_i)\in W_i, i=1,\dots, n\}\]
for any $n\in\IN, v_i\in V$ and $W_i$ open subspace in $W$ for all $i=1,\dots, n$.In particular, if $W=\sfk$ with 
the discrete topology, the space $V^*=\Hom_{\sfk}(V,\sfk)$ has a basis of neighborhoods of zero given by 
orthogonal complements of finite--dimensional subspaces in $V$. When $V$ is finite--dimensional, $V^*$ 
coincides with the linear dual and the weak topology coincides with the discrete topology. The canonical
map $V\to {V^*}^*$ is a linear isomorphism, when $V$ is discrete, and it is not topological in general.\\

The space of formal power series in $\hbar$ with coefficients in a topological vector space $V$, $V[[\hbar]]=V\ctp\sfk[[\hbar]]$, 
is also a complete topological space with a natural structure of a topological $\sfk[[\hbar]]$-module. A topological 
$\sfk[[\hbar]]$-module is complete if it is isomorphic to $V[[\hbar]]$ for some complete $V$. The additive category of 
complete $\sfk[[\hbar]]$-module, denoted $\A$, where morphisms are continuous $\sfk[[\hbar]]$-linear maps, has a natural
symmetric monoidal structure. Namely, the tensor product on $\A$ is defined to be the quotient of the tensor product $V\ctp W$ by the image of the operator $\hbar\ten1-1\ten\hbar$. This tensor product will be still denoted by $\ctp$. There is an extension of scalar functor from the category of topological spaces to $\A$, mapping $V$ to $V[[\hbar]]$. This functor respects the tensor product, \ie $(V\ctp W)[[\hbar]]$ is naturally isomorphic to $V[[\hbar]]\ctp W[[\hbar]]$. 

%-----------------------------
%EQUICONTINUOUS 
%----------------------------

\subsection{Equicontinuous modules}

Fix a topological Lie algebra $\g$.

\begin{definition}
Let $V$ be a topological vector space. We say that $V$ is an \emph{equicontinuous} $\g$-module if:
\begin{itemize}
\item the map $\pi_V:\g\to\End_{\sfk} V$ is a continuous homomorphism of topological Lie algebras;
\item $\{\pi_V(g)\}_{g\in\g}$ is an equicontinuous family of linear operators,\ie for any open subspace $U\subseteq V$, there exists $U'$
such that $\pi_V(g)U'\subset U$ for all $g\in\g$.
\end{itemize}
\end{definition}

Clearly, a topological vector space with a \emph{trivial} $\g$-module structure is an equicontinuous $\g$-module. Moreover, given
 equicontinuous $\g$-modules $V,W, U$, the tensor product $V\ctp W$ has a natural structure of equicontinuous 
 $\g$-module and $(V\ctp W)\ctp U$ is naturally identified with $V\ctp (W\ctp U)$. The category of equicontinuous $\g$-modules is then 
 a symmetric monoidal category, with braiding defined by permutation of components. We denote this category by $\eRep\g$.
 
 %---------------------------------------------------------------%
%LIE BIALGEBRAS, SUBALGEBRAS, MANIN TRIPLES
%---------------------------------------------------------------%

\subsection{Lie bialgebras and Manin triples}\label{ss:lbamt}

A Manin triple is the data of a Lie algebra $\g$ with 
\begin{itemize}
\item a nondegenerate invariant inner product $\iip{}{}$;
\item isotropic Lie subalgebras $\gpm\subset\g$;
\end{itemize}
such that 
\begin{itemize}
\item $\Lg=\gp\oplus\gm$ as vector space;
\item the inner product defines an isomorphism $\gp\simeq\gm^*$;
\item the commutator of $\Lg$ is continuous with respect to the topology 
obtained by putting the discrete and the weak topology on $\gm,\gp$ 
respectively.
\end{itemize}

Under these assumptions, the commutator on $\gp\simeq\gm^*$ induces 
a cobracket on $\gm$, satisfying the cocycle condition \cite{drin-2}. 
Therefore, $\gm$ is canonically endowed with a Lie bialgebra structure.
Notice that, in absolute generality, $\gp$ is only a topological Lie bialgebra, 
\ie $\delta(\gp)\subset\gp\hat{\ten}\gp$. The inner product also gives rise
to an isomorphism of vector spaces $\gm\simeq\gm^{**}\simeq\gp^*$,
where the latter is the continuous dual, though this isomorphism does
not respect the topology. Conversely, every Lie bialgebra $\a$ 
defines a Manin triple $(\a\oplus\a^*, \a, \a^* )$.

%CUTTED GRADED MANIN TRIPLE
\Omit{
It is known that, if $\gm$ is $\IN$-graded of finite type (i.e. with $\dim(\gm)_n<\infty$ for all $n\in\IN$), then we can formulate the definition of \emph{graded} Manin triple where the isomorphism $\gp\simeq\gm^*$ is instead an isomorphism with the \emph{graded} dual of $\gm$, i.e. $\gm^*=\bigoplus_{n\in\IN}(\gm)_n^*$. Under this assumption, $\gp$ and then $\Lg$ have natural structures of Lie bialgebras.\\ 
} 

% CONVERGENCE LEMMA
\Omit{
Let then $(\g,\gp,\gm)$ be a Manin triple. Let $\{a_i, i\in I\}$ be a basis for $\gm$ and $\{b_i, i\in I\}\subset\gp$ be the linear functionals on $\gm$ defined by $\iip{b^i}{a_j}=\delta_{ij}$.

\begin{lemma}
Let $V$ be a topological vector space with a continuous homomorphism $\g\to\End_{\sfk}V$. 
Then for all $v\in V$ and open subspace $U\subset V$, one has $b^i v\in U$ for all but finitely many $i\in I$.
\end{lemma}

\begin{pf}
A basis of neighborhood of zero in $\gp$ is given by sets
\[Y_j=\{f\in\gp\;|\; \iip{f}{a_j}=0\}\]
Let $\{i_s\in I\}$ be any sequence of distinct elements. Since for any $Y_j$ there exists $i_{s_0}\in I$ such that 
$b^{is}\in Y_j$ for $s>s_0$, then $\lim_{s\to\infty}b^{i_s}=0$. By continuity of the map $\g\to\End M$, 
\[\lim_{s\to\infty}b^{i_s}v=0 \qquad\mbox{for any}\quad v\in V\]
This means that $b^i v\in U$ for almost all $i\in I$.
\end{pf}
}

%-------------------------------------
% VERMA MODULES
%-------------------------------------

\subsection{Verma modules}

In \cite{ek-1}, Etingof and Kazhdan constructed two main examples of equicontinuous $\g$-modules in the case when 
$\g$ belongs to a Manin triple $(\g,\gp,\gm)$. The modules $M_{\pm}$, defined as
\[\Mp=\Ind_{\gm}^{\g}\sfk\qquad \Mm=\Ind_{\gp}^{\g}\sfk\]
are freely generated over $\UE{\gpm}$ by a vector $1_{\pm}$ such that $\gmp1_{\pm}$. Therefore, they are naturally identified, 
as vector spaces, to $\UE{\gpm}$ via $x1_{\pm}\to x$. The modules $\Mm$ and $\Mpd$, with appropriate topologies, 
are equicontinuous $\g$-modules.\\

The module $\Mm$ is an equicontinuous $\g$-module with respect to the discrete topology. The topology on $\Mp$ comes, 
instead, from the identification of vector spaces 
\[\Mp\simeq \UE{\gp}=\bigcup_{n\geq0} \UE{\gp}_n\] 
where $\UE{\gp}_n$ is the set of elements of degree at most $n$. The topology on $U(\gp)_n$ is defined through
the linear isomorphism
\[\xi_n:\bigoplus_{j=0}^n S^j\gp\to \UE{\gp}_n\]
where $S^j\gp$ is considered as a topological subspace of $(\gm^{\ten j})^*$, embedded with the weak topology.
Finally, $U(\gp)$ is equipped with the topology of the colimit. Namely, a set $U\subseteq U(\gp)$ is open if and
only if $U\cap U(\gp)_n$ is open for all $n$. With respect to the topology just described, the action of $\g$ on $\Mp$ 
is continuous.\\

Consider now the vector space of continuous linear functionals on $\Mp$
\[\Mpd=\Hom_{\sfk}(\Mp,\sfk)\simeq\colim\Hom_{\sfk}(U(\gp)_n, \sfk)\]
It is natural to put the discrete topology on $U(\gp)_n^*$, since, as a vector space, 
\[U(\gp)_n^*\simeq\bigoplus_{j=0}^n S^j\gp^*\simeq\bigoplus_{j=0}^n S^j\gm\simeq U(\gm)_n\]
We then consider on $\Mpd$ the topology of the limit. This defines, in particular, a filtration by subspaces 
$(\Mpd)_n$ satisfying
\[0\to(\Mpd)_n\to\Mpd\to(\UE{\gp}_n)^*\to 0\]
and such that $\Mpd=\lim \Mpd/(\Mpd)_n$. The topology of the limit on $\Mpd$ is, in general, stronger than the weak topology of 
the dual. Since the action of $\g$ on $\Mp$ is continuous, $\Mpd$ has a natural structure of $\g$--module. In particular, this is an 
equicontinuous $\g$--action.

%---------------------------------------------------------------------
% CASIMIR ELEMENT AND DRINFELD CATEGORY
%---------------------------------------------------------------------

\subsection{Drinfeld category}

The natural embedding
\[\gm\ten\gm^*\subset\End_{\sfk}(\gm)\]
induces a topology on $\gm\ten\gm^*$ by restriction of the weak topology in 
$\End_{\sfk}({\gm})$. With respect to this topology, the image of $\gm\ten\gm^*$ is dense in 
$\End_{\sfk}(\gm)$ and the topological completion $\gm\ten\gm^*$ is identified with
$\End_{\sfk}(\gm)$. Under this identification, the identity operator defines an element $r\in\gm\ctp\gm^*$.\\

Given two equicontinuous $\g$--modules $V, W$, the map
\[\pi_V\ten\pi_W:\gm\ten\gm^*\to\End_{\sfk}(V\ctp W)\]
naturally extends to a continuous map $\gm\ctp\gm^*\to\End_{\sfk}(V\ctp W)$. Therefore, the Casimir operator
\[\Omega=r+r^{\operatorname{op}}\in\gm\ctp\gm^*\oplus\gm^*\ctp\gm\]
defines a continuous endomorphism of $V\ctp W$, $\Omega_{VW}=(\pi_V\ten\pi_W)(\Omega)$, commuting with 
the action of $\g$.\\

Following \cite{drin-3}, it is possible to define a structure of braided monoidal category on the category of
deformed equicontinuous $\g$--module, depending on the choice of a Lie associator $\Phi$, the bifunctor
$\ctp$ and the Casimi operator $\Omega$. The commutativity constraint is explicitly defined by the formula
\[\beta_{VW}=(12)\circ e^{\frac{\hbar}{2}\Omega_{VW}}\in\Hom_{\g}(V\ctp W, W\ctp V)[[\hbar]]\]
We denote this braided tensor category braided tensor category $\DC{U{\g}}{\Phi}$. The category of equicontinuous 
$\g$--modules is equivalent to the category of Yetter-Drinfeld module over $\gm$, $\DY({\gm})$. The equivalence holds
at the level of tensor structure induced by the choice of an associator $\Phi$,
\[\DC{U{\g}}{\Phi}\simeq\TDY{\Ueh{\gm}}{\Phi}\]

%---------------------------------------------------
%PROPERTIES OF VERMA MODULES
%----------------------------------------------------

\subsection{Verma modules}
%------------------------------------------------------------------

The modules $M_{\pm}$ are identified, as vector spaces, with the enveloping universal algebras $U\gpm$. 
Their comultiplications induce the $U\g$--intertwiners $i_{\pm}:M_{\pm}\to M_{\pm}\ctp M_{\pm}$, mapping 
the vectors $1_{\pm}$ to the $\gmp$-invariant vectors $1_{\pm}\ten 1_{\pm}$.\\

For any $f,g\in\Mpd$, consider the linear functional $\Mp\to\sfk$ defined by $v\mapsto (f\ten g)(\ip(v))$. This f
unctional defines a map $\ips:\Mpd\ten\Mpd\to\Mpd$, that is continuous and extends to a morphism in 
$\Rep{\g}[[\hbar]]$, $\ips:\Mpd\ctp\Mpd\to\Mpd$. The pairs $(\Mm,\im)$ and $(\Mpd, \ips)$ form, respectively, 
a coalgebra and an algebra object in $\DC{U{\g}}{\Phi}$.\\

For any $V\in\DC{U{\g}}{\Phi}$, the vector space $\Hom_{\g}(\Mm,\Mpd\ctp V)$ is naturally 
isomorphic to $V$, as topological vector space, through the isomorphism $f\mapsto (1_+\ten 1)f(1_-)$.

%---------------------------------------------------
%EKQUANTIZATION
%----------------------------------------------------

\subsection{The fiber functor and the EK quantization}
%------------------------------------------------------------------
We will now recall the main results from \cite{ek-1,ek-2}. Where no confusion is possible, 
we will abusively denote $\ctp$ by $\ten$. Let then $F$ be the functor
\[F: \DC{U{\g}}{\Phi}\to\A\qquad F(V)=\Hom_{\DC{U{\g}}{\Phi}}(\Mm,\Mpd\ten V)\]
There is a natural transformation 
\[J\in\Nat{\ten\circ(F\boxtimes F)}{F\circ\ten}\]
defined, for any $v\in F(V), w\in F(W)$, by
\[J_{VW}(v\ten w)= (\ipd\ten1\ten1)A^{-1}\beta_{23}^{-1}A(v\ten w)\im\]
where $A$ is defined as a morphism
\[(V_1\ten V_2)\ten (V_3\ten V_4)\to V_1\ten( (V_2\ten V_3)\ten V_4)\]
by the action of $(1\ten\Phi_{2,3,4})\Phi_{1,2,34}$.
\begin{theorem}
The natural transformation $J$ is invertible and defines a tensor structure on the functor $F$. 
\end{theorem}
The tensor functor $(F,J)$ is called \emph{fiber functor}. The algebra of endomorphisms 
of $F$ is therefore naturally endowed with a topological bialgebra structure, as described 
in the previous section.\footnote{By \emph{topological} bialgebra we do not mean 
\emph{topological over $\sfk[[\hbar]]$}. We are instead referring to the fact that the algebra 
$\sfEnd{F}$ has a natural comultiplication $\Delta:\sfEnd{F}\to\sfEnd{F^2}$, where $\sfEnd{F^2}$ 
can be interpreted as an appropriate completion of $\sfEnd{F}^{\ten 2}$.}\\
 
The object $F(\Mm)\in\A$ has a natural structure of Hopf algebra, defined by the multiplication
\[m: F(\Mm)\ten F(\Mm)\to F(\Mm)\qquad m(x,y)=(\ipd\ten1)\Phi^{-1}(1\ten y)x\]
and the comultiplication
\[\Delta:F(\Mm)\to F(\Mm)\ten F(\Mm)\qquad \Delta(x)=J^{-1}(F(\im)(x))\]
The algebra $F(\Mm)$ is naturally isomorphic as a vector space with $\Mm[[\hbar]]\simeq 
U\gm[[\hbar]]$ and 
\begin{theorem}
The algebra $\qEK{\gm}=F(\Mm)$ is a quantization of the algebra $U\gm$.
\end{theorem}

In \cite{ek-2}, it is shown that this construction defines a functor 
\[Q\oEK:{\sc LBA}(\sfk)\to{\sc QUE}(K)\]
where ${\sc LBA}(\sfk)$ denotes the category of Lie bialgebras over $\sfk$ and {\sc QUE}(K) 
denotes the category of quantum universal enveloping algebras over $K=\sfk[[\hbar]]$. Another 
important result in \cite{ek-2} states the invertibility of the functor $Q\oEK$.\\

The map
\[m_-:\qEK{\gm}\to\sfEnd{F}\qquad m_-(x)_V(v)=(\ipd\ten 1)\Phi^{-1}(1\ten v)x\]
where $V\in\TDY{\Ueh{\gm}}{\Phi}$ and $v\in F(V)$, is, indeed, an inclusion of Hopf algebras. 
The map $m_-$ defines an action of $\qEK{\gm}$ on $F(V)$. Moreover, the map
\[F(V)\to F(\Mm)\ten F(V)\qquad v\mapsto R_{J}(1\ten v)\]
where $R_{J}$ denotes the twisted $R$--matrix, defines a coaction of $\qEK{\gm}$ on $F(V)$ 
compatible with the action, therefore
\begin{theorem}
The fiber functor $F:\TDY{\Ueh{\gm}}{\Phi}\to\A$ lifts to an equivalence of braided tensor categories
\[\wt{F}:\TDY{\Ueh{\gm}}{\Phi}\to\TDY{\qEK{\gm}}{}\]
\end{theorem}

%----------------------
% KAC--MOODY
%----------------------

\subsection{Generalized Kac-Moody algebras}
%------------------------------------------------------------

Denote by $\sfk$ a field of characteristic zero. We recall definitions from \cite{K} and \cite{ek-6}. 
Let $\sfA=(a_{ij})_{i,j\in \bfI}$ be an $n\times n$ symmetrizable matrix with entries in $\sfk$, i.e. there 
exists a (fixed) collection of nonzero numbers $\{d_i\}_{i\in \bfI}$ such that $d_ia_{ij}=d_ja_{ji}$ for all 
$i,j\in \bfI$. Let $(\Lh, \Pi,\Pi^{\vee})$ be a realization of $\sfA$. It means that $\Lh$ is a vector space 
of dimension $2n-\rank(\sfA)$, $\Pi=\{\alpha_1,\dots,\alpha_n\}\subset\Lh^*$ and 
$\Pi^{\vee}=\{h_1,\dots,h_n\}\subset\Lh$ are linerly independent, and $(\alpha_i,h_j)=a_{ji}$.

\begin{definition}
The Lie algebra $\tilde{\Lg}=\tilde{\Lg}(\sfA)$ is generated by $\Lh, \{e_i,f_i\}_{i\in \bfI}$ with 
defining relations
\[[h,h']=0\quad h,h'\in\Lh;\qquad [h,e_i]=(\alpha_i,h)e_i\]\[ [h,f_i]=-(\alpha_i,h)f_i;
\qquad [e_i,f_j]=\delta_{ij}h_i\]
\end{definition}

There exists a unique maximal ideal $\r$ in $\tilde{\Lg}$ that intersect $\Lh$ trivially. Let 
$\Lg:=\tilde{\Lg}/\r$. The algebra $\Lg$ is called \emph{generalized Kac-Moody algebra}. 
The Lie algebra $\Lg$ is graded by principal gradation $\deg(e_i)=1,\deg(f_i)=-1,\deg(\Lh)=0$, 
and the homogenous component are all finite-dimentional.\\ 

Let us now choose a non--degenerate bilinear symmetric form on $\Lh$ such that 
$\iip{h}{h_i}=d_i^{-1}(\alpha_i, h)$. Following \cite{K}, there exists a unique extension 
of the form $\iip{}{}$ to an invariant symmetric bilinear form on $\tilde{\Lg}$. For this 
extension, one gets $\iip{e_i}{f_j}=\delta_{ij}d_i^{-1}$. The kernel of this form on $\wt{\g}$ 
is $\r$, therefore it descends to a non--degenerate bilinear form on $\g$.\\

Let $\n_{\pm},\b_{\pm}$ be the nilpotent and the Borel subalgebras of $\g$, \ie $\n_{\pm}$ are 
generated by $\{e_i\}$, $\{f_i\}$, respectively, and $\b_{\pm}:=\n_{\pm}\oplus\h$. Since 
$[\n_{\pm}, \h]\subset\n_{\pm}$, we get Lie algebra maps $\bar{\;}:\b_{\pm}\to\h$ and we 
can consider the embeddings of Lie subalgebras $\eta_{\pm}:\b_{\pm}\to\g\oplus\h$ given 
by
$$\eta_{\pm}(x)=(x,\pm\bar{x})$$
Define the inner product on $\g\oplus\h$ by $\iip{}{}_{\g\oplus\h}=\iip{}{}_{\g}-\iip{}{}_{\h}$. 

\begin{proposition}
The triple $(\g\oplus\h,\b_+,\b_-)$ with inner product $\iip{-}{-}_{\g\oplus\h}$ and embeddings 
$\eta_{\pm}$ is a graded Manin triple.
\end{proposition}

Under the embeddings $\eta_{\pm}$, the Lie subalgebras $\b_{\pm}$ are isotropic with 
respect to $\iip{}{}_{\g\oplus\h}$. Since $\iip{}{}_{\g}$ and $\iip{}{}_{\h}$ are invariant 
symmetric non--degenerate bilinear form, so is $\iip{}{}_{\g\oplus\h}$.\\

The proposition implies that $\g\oplus\h,\b_+,\b_-$ are Lie bialgebras. Moreover, 
$\b_+^*\simeq\b_-^{cop}$ as Lie bialgebras (where $\b_+^*:=\bigoplus (\b_+)_n^*$ 
denotes the \emph{restricted} dual space and by $cop$ we mean the opposite cocommutator). 
The Lie bialgebra structures on $\b_{\pm}$ are then described by the following formulas:
\[\delta(h)=0,\quad h\in\h\subset\b_{\pm};\]\[\delta(e_i)=\frac{d_i}{2}(e_i\ten h_i-h_i\ten e_i)=
\frac{d_i}{2}e_i\wedge h_i;\quad\delta(f_i)=\frac{d_i}{2}f_i\wedge h_i\]
The Lie subalgebra $\{(0,h)\;|\; h\in\h\}$ is therefore an ideal and a coideal in $\g\oplus\h$. 
Thus, the quotient $\g=(\g\oplus\h)/\h$ is also a Lie bialgebra with Lie subbialgebras 
$\b_{\pm}$ and the same cocommutator formulas.

%--------------------------------------------------------------------
% EK QUANTIZATION OF KACMOODY ALGEBRAS
%--------------------------------------------------------------------

\subsection{Quantization of Kac--Moody algebras and category $\O$}\label{ss:isoek}
%--------------------------------------------------------------------

In \cite{ek-6}, Etingof and Kazhdan proved that, for any symmetrizable irreducible Kac-Moody 
algebra $\g$, the quantization $\qEK{\g}$ is isomorphic with the Drinfeld--Jimbo quantum group 
$\DJ{\g}$.\\

In particular, they construct an isomorphism of Hopf algebras $\DJ{\b_+}\simeq\qEK{\b_+}$, 
inducing the identity on $\Ueh{\h}$, where $\b_+$ is the Borel subalgebra and $\h$ is the 
Cartan subalgebra of $\g$.
Thanks to the compatibility with the doubling operations
\[\D\qEK{\b_+}\simeq\qEK{\Db\b_+}\]
proved by Enriquez and Geer in \cite{EG}, the isomorphism for the Borel subalgebra induces  
an isomorphism $\DJ{\g}\simeq\qEK{\g}$.\\

Recall that the category $\O$ for $\g$, denoted $\O_{\g}$ is defined to be the category of all 
$\h$--diagonalizable $\g$--modules $V$, whose set of weights $\mathsf{P}(V)$ belong to a 
union of finitely many cones 
\[\mathsf{D}(\lambda_s)=\lambda_s+\sum_{i}\IZ_{\geq0}\alpha_i\qquad \lambda_s\in\h^*, s=1,...,r\]
and the weight subspaces are finite--dimensional. We denote by $\O_{\g}[[\hbar]]$ the category of  
deformation $\g$--representations,\ie representations of $\g$ on topologically free $\sfk[[\hbar]]$--modules 
with the above properties (with weights in $\h^*[[\hbar]]$).\\

In a similar way, one defines the category $\O_{\DJ{\g}}$: it is the category of $\DJ{\g}$--modules 
which are topologically free over $\sfk[[\hbar]]$ and satisfy the same conditions as in the classical 
case.\\ 

The morphism of Lie bialgebras
\[\Db\b_+\to\g\simeq\Db\b_+/(\h\simeq\h^*)\]
gives rise to a pullback functors
\[\O_{\g}\to\TDY{U\b_+}{}\qquad \O_{\g,\Phi}[[\hbar]]\to\TDY{\Ueh{\b_+}}{\Phi}\]
where $\O_{\g,\Phi}$ denotes the category $\O_{\g}$ with the tensor structure of the Drinfeld category. 
Similarly, the morphism of Hopf algebras
\[D\qEK{\b_+}\to\qEK{\g}\simeq\DJ{\g}\]
gives rise to a pullback functor
\[\O_{\DJ{\g}}\to\TDY{\qEK{\b_+}}{}\]
\begin{theorem}
The equivalence $\wt{F}$ reduces to an equivalence of braided tensor categories
\[\wt{F}_{\O}:\O_{\g,\Phi}[[\hbar]]\to\O_{\DJ{\g}}\]
which is isomorphic to the identity functor at the level of $\h$--graded $\sfk[[\hbar]]$--modules. 
\end{theorem}

%----------------------------------
% THE ISOMORPHISM EK
%----------------------------------

\subsection{The isomorphism $\Psi\oEK$}
In \cite{ek-6}, Etingof--Kazhdan showed that the equivalence $\wt{F}$ induces an isomorphism of algebras
\begin{equation*}
\Psi\oEK: \wh{\Ueh{\g}}\to\wh{\DJ{\g}}
\end{equation*}
where 
\begin{equation*}
\wh{U\g}=\lim U_{\beta}\qquad U_{\beta}=U\g/I_{\beta}, \beta\in\IN^I
\end{equation*}
$I_{\beta}$ being the left ideal generated by elements of weight less or equal $\beta$ (analogously 
for $\wh{\DJ{\g}}$, cf.~\cite[Sec. 4]{ek-6}). 
\begin{proposition}
The isomorphism $\Psi\oEK$ coincides with the isomorphism induced by the equivalence $\wt{F}_{\O}$, 
as explained in Section \ref{ss:endiso}. 
\end{proposition}
\begin{pf}
The identification of the two isomorphism is constructed in the following way:
\begin{itemize}
\item[(a)] First, we show that there is a canonical map
\begin{equation*}
\sfEnd{{\sf f}_{\O}}\to{\sf C}_{\End(\wh{U})}(\End_{\g}(\wh{U}))
\end{equation*}
\item[(b)] There is a \emph{canonical}  multiplication in $\wh{U}$, so that
	\begin{itemize}
	\item[(i)] There is a canonical map \[{\sf C}_{\End(\wh{U})}(\End_{\g}(\wh{U}))\to\wh{U}\]
	\item[(ii)] For every $V\in\O$ the action of $U\g$ lifts to an action of $\wh{U}$
	\begin{equation*}
	\xymatrix{U\g\ar[rr]\ar[dr]&&\End(V)\\&\wh{U}\ar[ur]&}
	\end{equation*}
	\end{itemize}
\item[(c)] It defines a map $\wh{U}\to\sfEnd{{\sf f}_{\O}}$ and we have an isomorphism of algebras \[\wh{U}\simeq\sfEnd{{\sf f}_{\O}}\]
\end{itemize}
\end{pf}

If $\g$ is a semisimple Lie algebra, the equivalence of categories $\wt{F}$ leads to an isomorphism of 
algebras
\[\Ueh{(\Db\b_+)}\simeq\D\qEK{\b_ +} \Longrightarrow \Ueh{\g}\simeq\DJ{\g}\]
which is the identity modulo $h$. Toledano Laredo proved in \cite[Prop. 3.5]{vtl-4} that such an 
isomorphism cannot be compatible with all the isomorphisms
\[\Ueh{\sl{2}^{\alpha_i}}\simeq\DJ{\sl{2}^{\alpha_i}}\qquad\forall i\]
where $\{\alpha_i\}$ are the simple roots of $\g$. This amounts to a simple proof that the 
isomorphism $\Psi\oEK$ cannot be, in general, an isomorphism of $D$--algebras.

%---------------------------------------------------
% GENERALIZED EK QUANTIZATION
%---------------------------------------------------

\section{A relative Etingof--Kazhdan functor}\label{s:Gamma}
%=================================

\subsection{}
%--------------

In this section, we consider a split inclusion of Manin triples
\[i_D: (\gD,\gDp,\gDm)\hookrightarrow(\g,\gp,\gm)\]
We then define a relative version of the Verma modules $M_\pm$, and
use them to prove the following
\begin{theorem}\label{th:Gamma}
There is a tensor functor
\[\Gamma:\DC{U{\g}}{\Phi}\to\DC{U\gD}{\Phi_D}\]
canonically isomorphic, as abelian functor, to the restriction functor $i_D^*$. 
\end{theorem}

%SPLIT MANIN PAIR
%---------------------------

\subsection{Split inclusions of Manin triples}\label{ss:split}
%-----------------------------------------------------

\begin{definition}
An embedding of Manin triples
\[i:(\gD,\gDm,\gDp)\longrightarrow(\g,\gm,\gp)\]
is a Lie algebra homomorphism $i:\gD\to\g$ preserving inner products,
and such that $i(\g_{D,\pm})\subset\g_\pm$.
\end{definition}

%It is easy to see that such an $i$ is injective.
Denote the restriction of $i$ to $\g_{D,\pm}$ by $i_\pm$. $i_\pm$
give rise to maps $p_\pm=i_\mp^*:\g_\pm\to\g_{D,\pm}$, defined
via the identifications $\g_\pm\simeq\g_\mp^*$ and $\g_{D,\pm}
\simeq\g_{D,\mp}^*$ by
\[\iip{p_\pm(x)}{y}_{D}=\iip{x}{i_\mp(y)}\] 
for any $x\in\g_{\pm}$ and $y\in\g_{D,\mp}$. These map satisfy
$p_\pm\circ i_\pm=\id_{\g_{D,\pm}}$ since, for any $x\in\g_{D,\pm}$,
$y\in\g_{D,\mp}$,
\[\iip{p_\pm\circ i_\pm(x)-x}{y}_{D}=\iip{i_\pm(x)}{i_\mp(y)}-\iip{x}{y}_D=0\]
This yields in particular a a direct sum decomposition $\g_\pm=i(\g_{D,\pm})
\oplus\Lmpm$, where
\[\Lmpm=\Ker(p_\pm)=\g_\pm\cap i(\gD)^\perp\]

\begin{definition}
The embedding $i:(\gD,\gDm,\gDp)\longrightarrow(\g,\gm,\gp)$ is called
{\it split} if the subspaces $\Lmpm\subset\gpm$ are Lie subalgebras.
\end{definition}

\subsection{Split pairs of Lie bialgebras}
%------------------------------------------------

For later use, we reformulate the above notion in terms of bialgebras
via the double construction.

\begin{definition}\label{def:splitlba}
A \emph{split pair} of Lie bialgebras is the data of
\begin{itemize}
\item Lie bialgebras $(\a, [,]_{\a},\delta_{\a})$ and $(\b, [,]_{\b}, \delta_{\b})$.
\item Lie bialgebra morphisms $i:\a\to\b$ and $p:\b\to\a$ such that $p\circ i=\id_{\a}$.
\end{itemize}
\end{definition}

\begin{proposition}\label{pr:one to one}
There is a one--to--one correspondence between split inclusions of
Manin triples and split pairs of Lie bialgebras. Specifically,
\begin{enumerate}
\item If $i:(\gD,\gDm,\gDp)\longrightarrow(\g,\gm,\gp)$ is a split inclusion
of Manin triples, then $(\g_{D,-},\g_-,i_-,i_+^*)$ is a split pair of Lie bialgebras.
\item Conversely, if $(\a,\b,i,p)$ is a split pair of Lie bialgebras, then 
$i\oplus p^*:(\Db\a, \a, \a^*)\longrightarrow(\Db\b, \b, \b^*)$ is a split
inclusion of Manin triples.
\end{enumerate}
\end{proposition}

\subsection{Proof of (i) of Proposition \ref{pr:one to one}}\label{ss:(i) of one to one}
%---------------------------------------------------------------------

Given a split inclusion \[i=i_-\oplus i_+:(\gD,\gDm,\gDp)\longrightarrow(\g,\gm,\gp)\]
we need to show that $i_-$ and $i_+^*$ are Lie bialgebra morphisms.
By assumption, $i_-$ is a morphism of Lie algebras, and $i_+^*$ one
of coalgebras. Since $i_-=(i_-^*)^*$, it suffices to show that $p_\pm=i
_\mp^*$ preserve Lie brackets. 

We claim to this end that $\m_\pm$ are ideals in $\g_\pm$. Since
$[\m_\pm,\m_\pm]\subseteq\m_\pm$ by assumption, this amounts
to showing that $[i(\g_{D,\pm}),\m_\pm]\subseteq\m_\pm$. This
follows from the fact that $[i(\g_{D,\pm}),\m_\pm]\subseteq\g_\pm$,
and from
\begin{multline*}
\iip{[i(\g_{D,\pm}),\m_\pm]}{i(\g_{D,\mp})}\\
=
\iip{\m_\pm}{[i(\g_{D,\pm}),i(\g_{D,\mp})]}\subset
\iip{\m_\pm}{i(\g_{D,\pm})}+\iip{\m_\pm}{i(\g_{D,\mp})}
\end{multline*}
where the first term is zero since $\g_\pm$ is isotropic, and the second
one is zero by definition of $\m_\pm$.

Let now $X_1,X_2\in\g_\pm$, and write $X_j=i_\pm(x_j)+y_j$, where
$x_j\in\g_{D,\pm}$ and $y_j\in\m_\pm$. Since $\m_\pm=\Ker(p_\pm)$
and $p_\pm\circ i_\pm=\id$, we have $[p_\pm (X_1),p_\pm (X_2)]=
[x_1,x_2]$, while
\[p_\pm [X_1,X_2]=
p_\pm\left(i_\pm[x_1,x_2]+[i_\pm x_1,y_2]+[y_1,i_\pm x_2]+[y_1,y_2]\right)=
[x_1,x_2]\]
where the last equality follows from the fact that $\m_\pm$ is an ideal.

\subsection{Proof of (ii) of Proposition \ref{pr:one to one}}
%----------------------------------------------------------------------

The bracket on $\Db\a$ is defined by
\[[a,\phi]=\ad^*(a)(\phi)-\ad^*(\phi)(a)=-\iip{\phi}{[a,-]_{\a}}+\iip{\phi\ten\id}{\delta_{\a}(a)}\]
for any $a\in\a$, $\phi\in\a^*$. Analogously for $\Db\b$. Therefore, the equalities
\begin{align*}
\iip{p^*(\phi)\ten\id}{\delta_{\b}(i(a))}
=&\iip{\phi\ten\id}{(p\ten\id)(i\ten i)\delta_{\a}(a)}\\
=&\iip{\phi\ten\id}{(\id\ten i)\delta_{\a}(a)}=i(\iip{\phi\ten\id}{\delta_{\a}(a)})
\end{align*}
and
\begin{align*}
\iip{p^*(\phi)}{[i(a), b]_{\b}}=\iip{\phi}{p([i(a), b]_{\b})}=\iip{\phi}{[a, p(b)]_{\a}}
\end{align*}
for all $a\in\a$ and $b\in\b$, imply that the map $i\oplus p^*$ is a Lie algebra
map. It also respects the inner product, since for any $a\in\a,\phi\in\a^*$,
\[\iip{p^*(\phi)}{i(a)}=\iip{\phi}{p\circ i(a)}=\iip{\phi}{a}\]
Finally, $\m_-=\Ker(p)$ and $\m_+=\Ker i^*$ are clearly subalgebras.

%PARABOLIC SUBALGEBRAS
%----------------------------------------

\subsection{Parabolic Lie subalgebras}\label{s:parabolic} 
%-----------------------------------------------

Let
\[i_D=\im\oplus\ip:(\gD,\gDm,\gDp)\to(\g,\gm,\gp)\]
be a split embedding of Manin triples. We henceforth identify $\g_D$ as a
Lie subalgebra of $\g$ with its induced inner product, and $\g_{D,\pm}$ as
subalgebras of $\g_\pm$ noting that, by Proposition \ref{pr:one to one},
$\g_{D,-}$ is a sub Lie bialgebra of $\g_-$.

The following summarizes the properties of the subspaces $\Lmpm=
\gpm\cap\gD^{\perp}$ and $\Lppm=\Lmpm\oplus\gD$.
\begin{proposition}\hfill\break\vspace{-0.4cm}
\begin{itemize}
\item[(i)] $\Lmpm$ is an ideal in $\gpm$, so that $\gpm=\Lmpm\rtimes\gDpm$.\\
\item[(ii)] $[\gD,\m_\pm]\subset\Lmpm$, so that $\Lp_{\pm}=\m_\pm\rtimes\gD$
are Lie subalgebras of $\g$.\\
\item[(iii)] $\delta(\m_-)\subset\m_-\ten\g_{D,-}+\g_{D,-}\ten\m_-$, so that $\m_-
\subseteq\g_-$ is a coideal.\\
\end{itemize}
\end{proposition}

\begin{pf} (i) was proved in \ref{ss:(i) of one to one}. (ii) Since
\[\iip{[\g_D,\m_\pm]}{\g_D}=\iip{\m_\pm}{[\g_D,\g_D]}=0\]
we have $[\g_D,\m_\pm]\subset\gD^{\perp}=\m_-\oplus\m_+$. Moreover,
\[\iip{[\g_D,\m_\pm]}{\m_\pm}=\iip{\g_D}{[\m_\pm,\m_\pm]}=\iip{\g_D}{\m_\pm}=0\]
since $\m_\pm$ is a subalgebra, and it follows that $[\g_D,\m_\pm]\subset\m_\pm$.
(iii) is clear since $\Lmm$ is the kernel of a Lie coalgebra map.
\end{pf}

\begin{rem}
If the inclusion $i_D$ is compatible with a finite type $\IN$--grading,
then $\Lmp\subset\gp$ is a coideal. Moreover, $\Lp_{\pm}$ are Lie subbialgebras
of $\g$ such that the projection $\Lp_\pm\to\gD$ is a morphism of bialgebras.
Namely, a finite type $\IN$--grading allows to define a Lie bialgebra structure on
$\g,\gp$. We then get a vector space decomposition $\gpm=\Lmpm\oplus\gDpm$
and a Lie bialgebra map $\gpm\to\gDpm$. It is also possible to define the
Lie subalgebras
\[\Lppm=\Lmpm\oplus\gD\subset\g\]
If we assume the existence of a compatible grading on $\g$ and $\gD$,\ie
preserved by $i_D$, then the natural maps 
\[\Lppm\subset\Lg\qquad\Lppm\to\gD\] 
are morphisms of Lie bialgebras.
\end{rem}

\subsection{The relative Verma Modules}\label{ss:rel Verma}
%-------------------------------------------------

\begin{definition}
Given a split embedding of Manin triples $\gD\subset\g$, and the corresponding
decomposition $\g=\Lmm\oplus\Lpp$, let $\Lm,\Np$ be the \emph{relative} Verma
modules defined by
\begin{equation*}
\Lm=\ind_{\Lpp}^{\g}\sfk\aand\Np=\ind_{\Lmm}^{\g}\sfk
\end{equation*}
\end{definition}

\begin{proposition}
The $\g$--modules $\Lm$ and $\Npd$ are equicontinuous.
\end{proposition}

\noindent
The description of the appropriate topologies on $\Lm$ and $\Npd$, and
the proof of their equicontinuity will be carried out in \ref{ss:Lm}--\ref{ss:end N+*}.

\subsection{Equicontinuity of $\Lm$}\label{ss:Lm}
%--------------------------------------------

As vector spaces, 
\[\Lm\simeq U\Lmm\subset U\gm\]
so it is natural to equip $\Lm$ with the discrete topology. The set of operators 
$\{\pi_{\Lm}(x)\}_{x\in\g}$ is then an equicontinuous family, and the continuity
of $\pi_{\Lm}$ reduces to checking that, for every element $v\in \Lm$, the set
\[Y_v=\{b\in\Lg_+ |\; b.v=0\}\]
is a neighborhood of zero in $\Lg_+$. Since $U\Lmm$ embeds naturally in $U
\Lg_-$ the proof is identical to \cite[Lemma 7.2]{ek-1}. We proceed by induction
on the length of $v=a_{i_1}\dots a_{i_n}\bf{1}_{-}$. If $n=0$, then $v=\bf{1}_-$
and $Y_v=\Lg_+$. If $n>1$, then assume $v=a_jw$, with $w=a_{i_1}\dots a_
{i_{n-1}}\bf{1}_-$ and $Y_w$ open in $\Lg_+$. For every $x\in\Lg_+$
$$x.v=x.(a_jw)=[x,a_j].w+(a_jx).w$$
Call $Z$ the subset of $\Lg_+$
$$Z=\{x\in\Lg_+ |\; [x,a_j]\in Y_w\}$$
Z is open in $\Lg_+$, by continuity of bracket $[,]$, and clearly $Z\cap Y_w\subset
Y_v$.

\subsection{Topology of $\Np$}
%-------------------------------------

As vector spaces,
$$\Np=\text{Ind}_{\Lmm}^{\Lg}{\sfk}\simeq U\Lp_+\simeq\colim U_n\Lp_+$$
where $\{U_n\Lp_+\}$ denotes the standard filtration of $U\Lp_+$, so that
$$U_n\Lp_+\simeq\bigoplus_{m=0}^nS^m\Lp_+=
\bigoplus_{i+j\leq n}\left(S^i\Lg_+\otimes S^j\Lg_{D,-}\right)$$
We turn this isomorphism into an isomorphism of topological vector spaces,
by taking on $S^i\Lg_+$ and $S^j\Lg_{D,-}$ the topologies induced by the
embeddings
$$S^i\Lg_+\hookrightarrow (\Lg_-^{\otimes i})^*
\aand S^j\Lg_{D,-}\hookrightarrow \Lg_{D,-}^{\otimes j}$$ 
With respect to these topologies, $U_m\Lp_+$ is closed inside $U_n\Lp_+$
for $m<n$, and we equip $\Np$ with the direct limit topology. We shall need the
following
\begin{lemma}
For any $x\in\Lg$, the map $\pi_{\Np}(x):\Np\to \Np$ is continuous.
\end{lemma}
\begin{pf}
We need to show that for any neighborhood of the origin $U\subset \Np$, there
exists a neighborhood of zero $U'\subset \Np$ such that $\pi_{\Np}(x)U'\subset
U$. The topology on $\Np$ comes from the decomposition $U\Lpp\simeq 
U\gp\ten U\gDm$, so that an open neighborhood of zero in $\Np$ has the form 
$U\ten U\gDm + U\gp\ten V$, with $U$ open in $U\gp$ and $V$ open in $U\gDm$. 
We apply the same procedure used in \cite[Lemma 7.3]{ek-1} to construct a set 
$U'\ten U\gDm$, with $U'$ open in $U\gp$, such that 
\[\pi_{\Np}(x)(U'\ten U\gDm)\subset U\ten U\gDm\subset U\ten U\gDm + U\gp\ten V\]
Since the topology on $U\gDm$ is discrete, the set $U'\ten U\gDm$ is open 
in $\Np$ and the lemma is proved.
\end{pf}

\subsection{Topology of $\Npd$}
%---------------------------------------

As vector spaces,
$$\Npd\simeq (U\Lp_+)^*\simeq\lim(U_n\Lp_+)^*$$
Define a filtration $\{(\Npd)_n\}$ on $\Npd$ by
$$0\to(\Npd)_n\to(U\Lp_+)^*\to(U_n\Lp_+)^*\to0$$
so that $\Npd\supset (\Npd)_0\supset(\Npd)_1\supset\cdots$, and
we get an isomorphism of vector spaces
\[\Npd\simeq\lim \Npd/(\Npd)_n\]
Finally, we use the isomorphism to endow $\Npd$ with the inverse limit topology.
\begin{lemma}
$\{\pi_{\Npd}(x)\}_{x\in\Lg}$ is an equicontinuous family of operators.
\end{lemma}
\begin{pf}
Since $\Lp_+$ acts on $N_+$ by multiplication, 
\[\Lp_+(\Npd)_n\subset(\Npd)_{n-1}\]
If $x\in\Lmm$ and $x_i\in U\Lp_+$ for $i=1,\dots, n$, then in $U\Lg$,
\[xx_1\cdots x_n=x_1\cdots x_nx-\sum_{i=0}^nx_1\cdots x_{i-1}[x_i,x]x_{i+1}\cdots x_n\]
where $[x_i,x]\in\g$. Iterating shows that $(x.f)(x_1\cdots x_n)=0$ if $f\in (\Npd)_n$,
so that $x(\Npd)_n\subset(\Npd)_n$. Then, for any neighborhood of zero of the form
$U=(\Npd)_n$, it is enough to take $U'=(\Npd)_{n+1}$ to get $\Lg(\Npd)_{n+1}\subset
(\Npd)_n$.
\end{pf}

\subsection{Equicontinuity of $\Npd$} \label{ss:end N+*}
%--------------------------------------------

\begin{lemma}
The map $\pi_{\Npd}: \Lg\to\emph{End}(\Npd)$ is a continuous map.
\end{lemma}
\begin{pf}
Since $\Lg_-$ is discrete, it is enough to check that, for any $f\in\Npd$ and $n\in\IN$,
the subset
$$Y(f,n)=\{b\in\Lg_+|\;b.f\in(\Npd)_n\}$$
is open in $\Lg_+$, i.e.
$$b^i.f\in(\Npd)_n\qquad\text{ for a.a. }i\in I$$
Since $f\in \Npd\simeq\lim \Npd/(\Npd)_n$, we have $f=\{f_n\}$ where $f_n$ is the class
of $f$ modulo $(\Npd)_n$. Therefore $b^i.f\in(\Npd)_n$ iff 
$$(b^i.f)_n=b^i.f_{n+1}=0$$
Now, for any $x_1\cdots x_n\in U_n\Lp_+$, we have
$$b^i.f_{n+1}(x_1\cdots x_n)=-f_{n+1}(b^ix_1\cdots x_n)=0$$
for a.a. $i\in I$ and the lemma is proved (it is enough to exclude the indices corresponding
to the generators involved in the expression of $f_{n+1}$).\\

As a vector spaces, we can identify
$$\Lp_+^*=\Lg_+^*\oplus\Lg_{D,-}^*\simeq \Lg_-\oplus\Lg_{D,+}=\Lp_-$$
We can give as a basis for $\Lp_+$ and $\Lp_-$
$$\Lp_+\supset\{\{b^i\}_{i\in I},\{a_r\}_{r\in \ID}\}\qquad\Lp_-\supset\{\{a_i\}_{i\in I},\{b^r\}_{r\in \ID}\}$$
and obvious relations
$$(b^i,a_j)=\delta_{ij}\qquad (b^i,b^r)=0$$
$$(a_r,a_j)=0\qquad (a_r,b^s)=\delta_{rs}$$
with $i,j\in I$, $r,s\in \ID$.
We can then identify $f_{n+1}$ with an element in $U_{n+1}\Lp_-$. Call $T_{n+1}(f)$
the set of indices of all $a_i$ involved in the expression of $f_{n+1}$. Excluding these
finite set of indices we get the result.
\end{pf}

%PROPERTIES OF VERMA MODULES
\subsection{Coalgebra structure on $L_-,N_+$}
%---------------------------------------------------------

Define $\g$--module maps
\[i_-:\Lm\to\Lm\ctp\Lm\aand i_+:\Np\to\Np\ctp\Np\]
by mapping $\1_\mp$ to $\1_\mp\otimes\1_\mp$. Note that, under the identification
$\Lm\simeq U\Lmm$ and $\Np\simeq U\Lpp$, $i_\mp$ correspond to the coproduct
on $U\Lmm$ and $U\Lpp$ respectively.

Following \cite[Prop. 1.2]{drin-4}, we consider the invertible element $T\in(U\Lg\ctp U
\Lg)[[\hbar]]$ satisfying relations:
\begin{eqnarray*}
S^{\ten 3}(\Phi^{321})\cdot(T\otimes 1)\cdot(\Delta\otimes 1)(T)&=&(1\otimes T)(1\otimes \Delta)(T)\cdot\Phi\\
T\Delta(S(a))&=&(S\otimes S)(\Delta(a))T
\end{eqnarray*}
Let $\Npd$ be as before and $f,g\in\Npd$. Consider the linear functional in $\Hom_{\sfk}(\Np,\sfk)$ 
defined by
\[v\mapsto (f\otimes g)(T\cdot i_+(v))\]
This functional is continuous, so it 
belongs to $\Npd$ and allow us to define the map
\[i_+^{\vee}\in\Hom_{\sfk}(\Npd\otimes\Npd,\Npd)[[\hbar]]\;,\;\; i_+^{\vee}(f\otimes g)(v)=(g\otimes f)(T\cdot i_+(v))\]
This map is continuous and extends to a map from $\Npd\ctp\Npd$ to $\Npd$. 
For any $a\in\g$, we have
\begin{align*}
i_+^{\vee}(a(f\otimes g))(v)&=(f\otimes g)((S\otimes S)(\Delta(a))T\cdot i_+(v))=\\
&=(f\otimes g)(T\Delta(S(a))\cdot i_+(v))=\\
&=i_+^{\vee}(f\otimes g)(S(a).v)=(a.i_+^{\vee}(f\otimes g))(v)
\end{align*}
and then $i_+^{\vee}\in\Hom_{\Lg}(\Npd\ctp\Npd,\Npd)[[\hbar]]$.

The following shows that $L_-$ and $N_+$ are coalgebra objects in the Drinfeld categories
of $\g$--modules and $(\g,\g_D)$--bimodules respectively.

\begin{proposition}\label{pr:coalgebra}
The following relations hold
\begin{itemize}
\item[(i)] $\Phi(i_-\otimes 1)i_-=(1\otimes i_-)i_-$.\\ %in $\Hom_{\g}(\Lm,\Lm^{\otimes 3})$.\\
\item[(ii)] $i_+^{\vee}(1\otimes i_+^{\vee})\Phi=i_+^{\vee}(i_+^{\vee}\otimes 1)S^{\ten 3}(\Phi_{D}^{-1})^{\rho}$
\end{itemize}
where $(-)^{\rho}$ denotes the right $\gD$--action on $\Npd$.
\end{proposition}
\Omit{the right $\gD$--action is given by $\iip{x^{\rho}f}{v}=\iip{f}{vS(x)}$}
\begin{pf}
We begin by showing that
\begin{equation}\label{eq:coalgebra}
\Phi(\1_-^{\otimes 3})=\1_{-}^{\otimes 3}\aand\Phi(\1_{+}^{\otimes 3})=\Phi_{D}(\1_{+}^{\otimes 3})
\end{equation}
To prove the first identity, it is enough to notice that, since $\gp\1_-=0$ and $\Omega=\sum 
(a_i\otimes b^i+b^i\otimes a_i)$, 
$$\Omega_{ij}(\1_{-}^{\otimes 3})=0$$
Then $\Phi(\1_{-}^{\otimes 3})=\1_{-}^{\otimes 3}$. To prove the second one, we notice that
$\Lmm\1_+=0$ and that we can rewrite 
$$\Omega=\sum_{j\in I_D}(a_j\otimes b^j + b^j\otimes a_j)+
\sum_{i\in I\setminus I_D}(a_i\otimes b^i + b^i\otimes a_i)=\Omega_{D}+
\sum_{i\in I\setminus I_D}(a_i\otimes b^i + b^i\otimes a_i)$$
where $\{a_j\}_{j\in I_D}$ is a basis of $\gDm$ and $\{b^j\}_{j\in I_D}$ is the dual basis of $\gDp$. 
Then
\[\Omega_{ij}(\1_{+}^{\otimes 3})=\Omega_{D,ij}(\1_{+}^{\otimes 3})\]
and, since for any element $x\in\gD$, the right and the left $\gD$-action coincide on $\1_{+}$, i.e. 
$x.\1_{+}=\1_{+}.x$, we have
\[\Omega_{ij}(\1_{+}^{\otimes 3})=(\1_{+}^{\otimes 3})\Omega_{D,ij}\]
and consequently $\Phi(\1_{+}^{\otimes 3})=\Phi_{D}(\1_{+}^{\otimes 3})$.

To prove (i), note that since the comultiplication in $U\Lmm$ is coassociative, we
have $(i_-\otimes 1)i_-=(1\otimes i_-)i_-$. We therefore have to show that $\Phi
(i_-\otimes 1)i_-=(1\otimes i_-)i_-$. This is an obvious consequence of \eqref
{eq:coalgebra} and the fact that $\Lmm$ is generated by $\1_-$.

To prove (ii), consider $v\in\Np$,
\begin{align*}
\begin{split}
\ipd&(1\ten\ipd)(\Phi(f\ten g\ten h))(v)=\\
&=(h\ten g\ten f)((S^{\ten 3}(\Phi^{321})\cdot(T\otimes 1)\cdot(\Delta\otimes 1)(T))\cdot(\ip\ten 1)\ip(v))=\\
&=(h\ten g\ten f)((1\otimes T)(1\otimes \Delta)(T)\cdot\Phi(\ip\ten 1)\ip(v))=\\
&=(h\ten g\ten f)((1\otimes T)(1\otimes \Delta)(T)(1\ten \ip)\ip(v)\Phi_D)=\\
&=(S^{\ten3}(\Phi_D)^{\rho}(h\ten g\ten f))((1\otimes T)(1\otimes \Delta)(T)(1\ten \ip)\ip(v))=\\
&=\ipd(\ipd\ten1)(S^{\ten3}(\Phi_D^{321})^{\rho}(f\ten g\ten h))(v)=\\
&=\ipd(\ipd\ten1)S^{\ten 3}(\Phi_D^{-1})^{\rho}(f\ten g\ten h)(v)
\end{split}
\end{align*}
and (ii) is proved.
\end{pf}

\subsection{The fiber functor over $\gD$}
%--------------------------------------------------

To any representation $V[[\hbar]]\in\Rep{U\g}[[\hbar]]$, we can associate the
$\sfk[[\hbar]]$--module
$$\Gamma(V)=\Hom_{\g}(\Lm,\Npd\ctp V)[[\hbar]]$$
where $\Hom_{\g}$ is the set of continuous homomorphisms, equipped with
the weak topology. The right $\g_D$--action on $N_+^*$ endows $\Gamma
(V)$ with the structure of a left $\g_D$--module.

\begin{proposition}\label{pr:isom to forget}
The complete vector space $\Hom_{\g}(\Lm,\Npd\ctp V)$ is isomorphic to $V$
as equicontinous $\gD$--module. The isomorphism is given by
\[\alpha_V: f\mapsto (\1_+\ten1)f(\1_-)\]
for any $f\in\Hom_{\g}(\Lm,\Npd\ctp V)$.
\end{proposition}
\begin{pf}
By Frobenius reciprocity, we get an isomorphism
\[
\Hom_{\g}(\Lm, \Npd\ctp V)\simeq\Hom_{\Lpp}(\sfk, \Npd\ctp V)\simeq\Hom_{\sfk}(\sfk, V)\simeq V\]
given by the map
\[ f\mapsto (\1_+\ten1)f(\1_-) \]
For $f\in \Gamma(V)$ and $x\in U\gD$, $x.f\in \Gamma(V)$ is defined by
\[x.f=(S(x)^{\rho}\otimes\id)\circ f\]
For any $x\in U\gD$, we have
\[\sum_{i,j}x_i^{(1)}f_j\ten x_i^{(2)}v_j=\varepsilon(x)f(1_-)\]
where $\Delta(x)=\sum_i x_i^{(1)}\ten x_i^{(2)}$ and $f(\1_-)=\sum_j f_j\ten v_j$.
Using the identity
\[1\ten x  = \sum_i (S(x_i^{(1)})\ten 1)\cdot\Delta(x_i^{(2)})\]
holding in any Hopf algebra, we obtain
\[(1\ten x)f(\1_-)=\sum_i (S(x_i^{(1)}\varepsilon(x_i^{(2)}))\ten 1) f(\1_-)
= (S(x)\ten 1)f(\1_-)\]
Finally, we have
\begin{align*}
x.\alpha_V(f)&=\iip{\1_+\ten\id}{(1\ten x)f(\1_-)}=\\&=\iip{\1_+\ten\id}{(S(x)\ten 1)f(\1_-)}=\\
&=\iip{\1_+\ten\id}{(S(x)^{\rho}\ten 1)f(\1_-)}=\alpha_V(x.f)
\end{align*}
Therefore, $\Gamma(V)$ is isomorphic to $V[[\hbar]]$ as equicontinuous $\gD$-module.
\end{pf}

\subsection{}
%--------------

For any continuous $\varphi\in\Hom_{\g}(V,V')$, define a map $\Gamma(\varphi):
\Gamma(V)\to \Gamma(V')$ by
\[\Gamma(\varphi):\; f\mapsto(\id\otimes \varphi)\circ f\]
This map is clearly continuous and for all $x\in\gD$
\[\Gamma(\varphi)(x.f)=(S(x)^{\rho}\ten\varphi)\circ f=x.\Gamma(\varphi)(f)\]
then $\Gamma(\varphi)\in\Hom_{\gD}(\Gamma(V),\Gamma(V'))$.

Since the diagram
\[\xymatrix{\Gamma(V)\ar[d]_{\alpha_{V}}\ar[r]^{\Gamma(\varphi)}&\Gamma(V')\ar[d]^{\alpha_{V'}}\\V[[\hbar]]\ar[r]^{\varphi}&V'[[\hbar]]}\]
is commutative for all $\varphi\in\Hom_{\g}(V,V')$, we have a well--defined functor 
\[\Gamma:\eRep{\Ueh{\g}}\to\eRep{\Ueh{\gD}}\]
which is naturally isomorphic to the pullback functor induced by the inclusion $i_{D}
:\gD\hookrightarrow \g$ via the natural transformation
\[\alpha_{V}: \Gamma(V)\simeq i_{D}^*V[[\hbar]]\]

%---------------------------------------------------
% TENSOR STRUCTURE
%----------------------------------------------------

\subsection{Tensor structure on $\Gamma$}\label{s:tensor-structure}
%------------------------------------------------------

Denote the tensor product in the categories $\DC{U{\g}}{\Phi}$, $\DC{U\gD}{\Phi_D}$
by $\ten$, and let $B_{1234}$ and $B'_{1234}$ be the associativity constraints
\[B_{1234}:(V_1\ten V_2)\ten(V_3\ten V_4)\to V_1\ten((V_2\ten V_3)\ten V_4)\] 
%the composition
%\[(V_1\ten V_2)\ten(V_3\ten V_4) \xrightarrow{\Phi_{1,2,34}} V_1\ten(V_2\ten (V_3\ten V_4))\xrightarrow{1\ten\Phi_{2,3,4}^{-1}}V_1\ten((V_2\ten V_3)\ten V_4)\]
and
\[B'_{1234}:(V_1\ten V_2)\ten(V_3\ten V_4)\to (V_1\ten(V_2\ten V_3))\ten V_4\]\\
%the composition
%\[(V_1\ten V_2)\ten(V_3\ten V_4)\xrightarrow{\Phi_{12,3,4}^{-1}}((V_1\ten V_2)\ten V_3)\ten V_4 \xrightarrow{\Phi_{1,2,3}} (V_1\ten(V_2\ten V_3))\ten V_4\]

\noindent
For any $v\in \Gamma(V),w\in \Gamma(W)$, define $J_{VW}(v\ten w)$ to be the
composition
\begin{multline*}
\Lm\xrightarrow{\im}\Lm\ten\Lm\xrightarrow{v\ten w}(\Npd\ten V)\ten(\Npd\ten W)\xrightarrow{A}\Npd\ten((V\ten\Npd)\ten W)\\
\xrightarrow{\beta^{-1}_{32}}\Npd\ten((\Npd\ten V)\ten W)\xrightarrow{A'}(\Npd\ten\Npd)\ten(V\ten W)\xrightarrow{\ipd\ten 1}\Npd\ten(V\ten W)
\end{multline*}
where the pair $(A,A')$ can be chosen to be $(B_{\Npd,V,\Npd,W}, B^{-1}_{\Npd,\Npd,V,W})$
or $(B'_{\Npd,V,\Npd,W}, B'^{-1}_{\Npd,\Npd,V,W})$. The map $J_{VW}(v\ten w)$ is clearly
a continuous $\g$-morphism from $\Lm$ to $\Npd\ten(V\ten W)$, so we have a well-defined
map
\[J_{VW}:\Gamma(V)\ten\Gamma(W)\to\Gamma(V\ten W)\]

\begin{proposition}\label{pr:tensor structure}
The maps $J_{VW}$ are isomorphisms of $\gD$--modules, and define a tensor structure
on the functor $\Gamma$.
\end{proposition}

\noindent
The proof of Proposition \ref{pr:tensor structure} is given in \ref{ss:start tensor}--\ref
{ss:end tensor}.

\subsection{}\label{ss:start tensor}
%--------------

The map $J_{VW}$ is compatible with the $\gD$--action.
Indeed, $\ipd$ is a morphism of right $\gD$--modules and, 
for any $x\in\gD$,
\begin{align*}
x.J_{VW}(v\ten w)&=(S^{\rho}(x)\ten\id)(\ipd\ten\id\ten\id)\wt{A}(v\ten w)\im=\\
&=(\ipd\ten\id\ten\id)(\Delta(S(x))^{\rho})_{12}\wt{A}(v\ten w)\im=\\
&=(\ipd\ten\id\ten\id)\wt{A}((S\ten S)(\Delta(x)))^{\rho})_{13}(v\ten w)\im=J_{VW}(x.(v\ten w))
\end{align*}
where $\wt{A}=A'\beta_{32}^{-1}A$.

$J_{VW}$ is an isomorphism, since it is an isomorphism modulo $\hbar$. Indeed,
$$J_{VW}(v\ten w)\equiv (i_+^*\ten 1)(1\ten s\ten 1)(v\ten w)i_-  \mod\hbar$$
To prove that $J_{VW}$ define a tensor structure on $\Gamma$, we need to show
that, for any $V_1,V_2,V_3\in\DC{U\g}{\Phi}$ the following diagram is commutative
\[{\scriptsize
\xymatrix@R=1cm@C=1.5cm{(\Gamma(V_1)\ten \Gamma(V_2))\ten \Gamma(V_3) \ar[d]_{\Phi_{D}} \ar[r]^(.53){J_{12}\ten 1} & \Gamma(V_1\ten V_2)\ten \Gamma(V_3) \ar[r]^{J_{12,3}} & \Gamma((V_1\ten V_2)\ten V_3) \ar[d]^{\Gamma(\Phi)} \\
\Gamma(V_1)\ten(\Gamma(V_2)\ten \Gamma(V_3)) \ar[r]^(.53){1\ten J_{23}} & \Gamma(V_1)\ten \Gamma(V_2\ten V_3) \ar[r]^{J_{1,23}} & \Gamma(V_1\ten(V_2\ten V_3))}}\]
where $J_{ij}$ denotes the map $J_{V_i,V_j}$ and $J_{i j,k}$ the map $J_{V_i\ten
V_j,V_k}$.

\subsection{}
%--------------

For any $v_i\in \Gamma(V_i)$, $i=1,2,3$, the map $\Gamma(\Phi)J_{12,3}J_{12}\ten1
(v_1 \ten v_2\ten v_3)$ is given by the composition
\begin{multline*}
%\Gamma(\Phi)J_{12,3}(J_{12}\ten1)(v_1\ten  v_2 \ten v_3)\\
(1\ten\Phi)(i_+^*\ten 1^{\ten 3})A_4(1\ten\beta_{1\ten2, \Npd}\ten1)A_3((i_+^*\ten 1)\ten1^{\ten 3})(A_2\ten1\ten1)\\
\cdot(1\ten\beta_{\Npd,1}\ten 1^{\ten 3})(A_1\ten1\ten1)(v_1\ten v_2\ten v_3)(i_-\ten1)i_-
\end{multline*}
where 
\[\begin{array}{ccc}
A_1=B_{\Npd,1,\Npd,2} &\hspace{.5cm} &A_3=B_{\Npd,1\ten2,\Npd,3}\\
A_2=B^{-1}_{\Npd,\Npd,1,2} &\hspace{.5cm} &A_4=B^{-1}_{\Npd,\Npd,1\ten2,3}
\end{array}\]
illustrated by the diagram
\[
{\scriptsize
\xymatrix@R=.5cm{\Lm\ar[r]^{i_-} & \Lm\ten\Lm \ar[r]^{i_-\ten 1} & (\Lm\ten\Lm)\ten\Lm \\
\ar[r]^(.1){v_1\ten v_2\ten v_3} & ((\Npd\ten V_1)\ten(\Npd\ten V_2))\ten(\Npd\ten V_3) \ar[r]^{A_1\ten 1\ten 1} & (\Npd\ten((V_1\ten\Npd)\ten V_2))\ten(\Npd\ten V_3) \\ 
\ar[r]^(.1){1\ten\beta_{\Npd,1}\ten1^{\ten 3}} &  (\Npd\ten((\Npd\ten V_1)\ten V_2))\ten(\Npd\ten V_3) \ar[r]^{A_2\ten 1\ten 1} & ((\Npd\ten\Npd)\ten(V_1\ten V_2))\ten(\Npd\ten V_3)\\
\ar[r]^(.2){(i_+^*\ten 1)\ten1^{\ten 3}} & (\Npd\ten(V_1\ten V_2))\ten(\Npd\ten V_3) \ar[r]^{A_3} & \Npd\ten(((V_1\ten V_2)\ten\Npd)\ten V_3)  \\
\ar[r]^(.2){1\ten\beta_{1\ten 2,\Npd}\ten 1} & \Npd\ten((\Npd\ten(V_1\ten V_2))\ten V_3) \ar[r]^{A_4} & (\Npd\ten\Npd)\ten((V_1\ten V_2)\ten V_3)  \\
\ar[r]^(.2){i_+^*\ten1^{\ten 3}} & \Npd\ten((V_1\ten V_2)\ten V_3) \ar[r]^{1\ten\Phi} & \Npd\ten(V_1\ten(V_2\ten V_3))}
}
\]
By functoriality of associativity and commutativity isomorphisms, we have
$$A_3(i_+^*\ten1^{\ten 4})=(i_+^*\ten1^{\ten 4})A_5$$
where $A_5=B_{\Npd\ten\Npd,12,\Npd,3}$,
$$(1\ten\beta_{12,\Npd}\ten1)(i_+^*\ten1^{\ten 4})=(i_+^*\ten1^{\ten 4})(1^{\ten2}\ten\beta_{12,\Npd}\ten1^{\ten2})$$
and
$$A_4(i_+^*\ten1^{\ten 4})=(i_+^*\ten1^{\ten 4})A_6$$
where $A_6=B^{-1}_{\Npd\ten\Npd,\Npd,1\ten2,3}$. 
\Omit{Then
$$(1\ten\Phi_{123})(i_+^*(i_+^*\ten1)\ten1^{\ten 3})=(i_+^*(i_+^*\ten1)\ten1^{\ten 3})(1^{\ten 3}\ten\Phi_{123})$$}
Finally, we have
\begin{multline}\label{eq:one side}
\Gamma(\Phi)J_{12,3}(J_{12}\ten1)(v_1\ten v_2\ten v_3)\\
=(1^{\ten 3}\ten\Phi_{123})((i_+^*(i_+^*\ten1))\ten 1^{\ten 3})
\,A\,(v_1\ten v_2\ten v_3)(i_-\ten 1)i_-
\end{multline}
where
$$A=A_6(1^{\ten 2}\ten\beta_{1\ten2,\Npd}\ten1^{\ten 2})A_5(A_2\ten1^{\ten 2})(1\ten\beta_{\Npd,1}\ten1^{\ten 3})(A_1\ten1\ten1)$$

\subsection{}
%--------------

On the other hand, $J_{1,2\ten3}(1\ten J_{23})\Phi_D(v_1\ten v_2\ten v_3)$ corresponds
to the composition
\begin{equation*}
\begin{split}
%J_{1,2\ten3}(1\ten J_{23})\Phi_D(v_1\ten & v_2\ten v_3)=\\=
(i_+^*\ten 1^{\ten 3})A_4'&(1\ten\beta_{\Npd,1}\ten1^{\ten 2})A_3'(1^{\ten 2}\ten i_+^*\ten1^{\ten 2})(1\ten1\ten A_2')\\&(1^{\ten 3}\ten\beta_{2,\Npd}\ten 1)(1\ten1\ten A_1')\Phi_D(v_1\ten v_2\ten v_3)(1\ten i_-)i_-
\end{split}
\end{equation*}
where 
\[\begin{array}{ccc}
A_1'=B_{\Npd,2,\Npd,3} &\hspace{.5cm} &A_3'=B_{\Npd,1,\Npd,2\ten3}\\
A_2'=B^{-1}_{\Npd,\Npd,2,3} &\hspace{.5cm} &A_4'=B^{-1}_{\Npd,\Npd,1,2\ten3}
\end{array}\]
illustrated by the diagram
\[
{\scriptsize
\xymatrix@R=.5cm{\Lm\ar[r]^{i_-} & \Lm\ten\Lm \ar[r]^{1\ten i_-} & \Lm\ten(\Lm\ten\Lm) \\
\ar[r]^(.1){\Phi_{D}(v_1\ten v_2\ten v_3)} & (\Npd\ten V_1)\ten((\Npd\ten V_2)\ten(\Npd\ten V_3)) \ar[r]^{1\ten 1\ten A_1'} & (\Npd\ten V_1)\ten(\Npd\ten ((V_2\ten\Npd)\ten V_3)) \\ 
\ar[r]^(.1){1^{\ten 3}\ten\beta_{2,\Npd}\ten1} &  (\Npd\ten V_1)\ten(\Npd\ten((\Npd\ten V_2)\ten V_3)) \ar[r]^{1\ten 1\ten A_2'} & (\Npd\ten V_1)\ten((\Npd\ten\Npd)\ten (V_2\ten V_3))\\
\ar[r]^(.2){1^{\ten 2}\ten i_+^*\ten1^{\ten 2}} & (\Npd\ten V_1)\ten(\Npd\ten (V_2\ten V_3)) \ar[r]^{A_3'} & \Npd\ten((V_1\ten\Npd)\ten (V_2\ten V_3))  \\
\ar[r]^(.2){1\ten\beta_{1,\Npd}\ten 1} & \Npd\ten((\Npd\ten V_1)\ten (V_2\ten V_3)) \ar[r]^{A_4'} & (\Npd\ten\Npd)\ten(V_1\ten (V_2\ten V_3))  \\
\ar[r]^(.2){i_+^*\ten1^{\ten 3}} & \Npd\ten(V_1\ten (V_2\ten V_3))}
}
\]
By functoriality of associativity and commutativity isomorphisms, we have
$$A_3'(1^{\ten 2}\ten i_+^*\ten 1^{\ten 2})=(1^{\ten 2}\ten i_+^*\ten 1^{\ten 2})A_5'$$
where $A_5'=B_{\Npd,1,\Npd\ten\Npd,2\ten 3}$,
$$(1\ten\beta_{1,\Npd}\ten1^{\ten 2})(1^{\ten 2}\ten i_+^*\ten 1^{\ten 2})=(1\ten i_+^*\ten 1^{\ten 3})(1\ten\beta_{1,\Npd\ten\Npd}\ten 1^{\ten 2})$$
and
$$A_4'(1\ten i_+^*\ten 1^{\ten 3})=(1\ten i_+^*\ten 1^{\ten 3})A_6'$$
where $A_6'=B^{-1}_{\Npd,\Npd\ten\Npd,1,2\ten3}$. Thus,
\Omit{\begin{equation}\label{eq:another side}
\begin{split}
J_{1,2\ten3}(1\ten J_{23})\Phi_D(v_1\ten & v_2\ten v_3)=\\=(i_+^*\ten 1^{\ten 3})&((1 \ten i_+^*)\ten1^{\ten 3})A_6'(1\ten\beta_{1,\Npd\ten\Npd}\ten1^{\ten 2})A_5'(1^{\ten 2}\ten A_2')\\&(1^{\ten 3}\ten\beta_{2,\Npd}\ten 1)(1\ten1\ten A_1')\Phi_D(v_1\ten v_2\ten v_3)(1\ten i_-)i_-
\end{split}
\end{equation}}
\begin{multline}\label{eq:another side}
J_{1,23}(1\ten J_{23})\Phi_D(v_1\ten v_2\ten v_3)\\
=
(i_+^*\ten 1^{\ten 3})((1 \ten i_+^*)\ten1^{\ten 3})
\,B\, %A_6'(1\ten\beta_{1,\Npd\ten\Npd}\ten1^{\ten 2})A_5'(1^{\ten 2}\ten A_2')\\&(1^{\ten 3}\ten\beta_{2,\Npd}\ten 1)(1\ten1\ten A_1')
\Phi_D(v_1\ten v_2\ten v_3)(1\ten i_-)i_-
\end{multline}
where
$$B=A_6'(1\ten\beta_{1,\Npd\ten\Npd}\ten1^{\ten 2})A_5'(1^{\ten 2}\ten A_2')(1^{\ten 3}\ten\beta_{2,\Npd}\ten 1)(1\ten1\ten A_1')$$

\subsection{}\label{ss:end tensor}
%--------------

Comparing \eqref{eq:one side} and \eqref{eq:another side}, we see that it suffices
to show that the outer arrows of the following form a commutative diagram.
\[
{\scriptsize
\xymatrix@R=1.5cm@C=.1cm{& &\Lm \ar[dll]_{(i_-\ten 1)i_-}  \ar[drr]^{(1\ten i_-)i_-}& & \\
(\Lm\ten\Lm)\ten\Lm \ar[d]_{v_1\ten v_2\ten v_3}\ar@{=}[rrrr] & & & & \Lm\ten(\Lm\ten\Lm) \ar[d]^{\Phi_{D}(v_1\ten v_2\ten v_3)}\\
((\Npd\ten V_1)\ten(\Npd\ten V_2))\ten(\Npd\ten V_3) \ar[d]_{A} \ar[rrrr]^{\Phi}&&&&(\Npd\ten V_1)\ten((\Npd\ten V_2)\ten(\Npd\ten V_3))\ar[d]^{B}\\
((\Npd\ten\Npd)\ten\Npd)\ten((V_1\ten V_2)\ten V_3)\ar[d]_{(i_+^*(i_+^*\ten 1))\ten1^{\ten 3}} \ar[rrrr]^{\Phi\ten\Phi}&&&&(\Npd\ten(\Npd\ten\Npd))\ten(V_1\ten (V_2\ten V_3))\ar[d]^{(i_+^*(1\ten i_+^*))\ten1^{\ten 3}}\\
\Npd\ten((V_1\ten V_2)\ten V_3)\ar[rrrr]^{1\ten\Phi}&&&&\Npd\ten(V_1\ten (V_2\ten V_3))
}
}
\]

Using the pentagon and the hexagon axiom, we can show that
$$(\Phi\ten\Phi)A=B\Phi$$
We have to show that 
\[\Gamma(\Phi)J_{12,3}(J_{12}\ten1)(v_1\ten  v_2  \ten v_3)=J_{1,23}(1\ten J_{23})\Phi_D(v_1\ten v_2\ten v_3)\]
in $\Hom_{\g}(\Lm,\Npd\ten(V_1\ten(V_2\ten V_3)))$:
\begin{equation*}
\begin{split}
J_{1,23}&(\id\ten J_{23})\Phi_D(v_1\ten v_2\ten v_3)=\\
&=(\ipd(\id\ten\ipd)\ten\id^{\ten 3})B\Phi_{D}(v_1\ten v_2\ten v_3)(\id\ten i_-)\im\\
&=(\ipd(\id\ten\ipd)\ten\id^{\ten 3})B\Phi_{D}(v_1\ten v_2\ten v_3)\Phi(\im\ten\id)\im\\
&=(\ipd(\id\ten\ipd)\ten\id^{\ten 3})B\Phi\Phi_{D}(v_1\ten v_2\ten v_3)(\im\ten\id)\im\\
&=(\ipd(\id\ten\ipd)\Phi\ten\Phi)A\Phi_{D}(v_1\ten v_2\ten v_3)(\im\ten\id)\im\\
&=(\ipd(\id\ten\ipd)\Phi\ten\Phi)(S^{\ten 3}(\Phi_D)^{\rho}\ten\id^{\ten 3})A(v_1\ten v_2\ten v_3)(\im\ten\id)\im\\
&=(\ipd(\id\ten\ipd)\Phi S^{\ten 3}(\Phi_D)^{\rho}\ten\Phi)A(v_1\ten v_2\ten v_3)(\im\ten\id)\im\\
&=(\ipd(\ipd\ten\id)\ten\Phi)A(v_1\ten v_2\ten v_3)(\im\ten\id)\im\\
&=\Gamma(\Phi)J_{12,3}(J_{12}\ten\id)(v_1\ten  v_2  \ten v_3)
\end{split}
\end{equation*}
where the second and seventh equalities follow from Proposition \ref{pr:coalgebra}, the fifth
one from the definition of the $\gD$--action on the modules $\Gamma(V_i)$ and the others 
from functoriality of the associator $\Phi$.  This complete the proof of Theorem \ref{th:Gamma}.

%----------------------------------------------
%COMPUTING THE 1-JETS
%----------------------------------------------

\subsection{$\mathbf{1}$--Jets of relative twists} 
%=====================================

The following is a straightforward extension of the computation
of the 1--jet of the Etingof--Kazhdan twist given in \cite{ek-1}.
\Omit{An easy computation in \cite{ek-1} shows that the twist $J$ of the 
Etingof--Kazhdan functor satisfies, under the identification $\tau$ with
the forgetful functor,
\[\tau_{V\ten W}\circ J \circ(\tau_{V}^{-1}\ten\tau_W)^{-1}\equiv1+\frac{\hbar}{2}r\quad\mod\hbar^2\]
in $\End(V\ten W)[[\hbar]]$. A similar result holds for the relative functor.}

\begin{proposition}\label{pr:jets-rel-twist}
Under the natural identification
\[\alpha_V:\Gamma(V)\to V[[\hbar]]\]
the relative twist $J_{\Gamma}$ satisfyies
\[\alpha_{V\ten W}\circ J_{\Gamma}\circ(\alpha_V^{-1}\ten\alpha_W^{-1})\equiv 1
+\frac{\hbar}{2}(r+r_D^{21})\quad\mod\hbar^2\]
in $\End(V\ten W)[[\hbar]]$.
\end{proposition}

\begin{pf}
For $v\in V, w\in W$, let 
\[\alpha^{-1}_V(v)(1_-)=\sum f_i\ten v_i\qquad \alpha^{-1}_W(w)(1_-)=\sum g_j\ten w_j\]
in $(\Npd\ten V)^{\Lpp}$ and $(\Npd\ten W)^{\Lpp}$ respectively. Then using
\[\langle (1_+\ten 1)^{\ten2},\ol{\Omega}_{23}\sum_{i,j}f_i\ten v_i\ten g_j\ten w_j\rangle= -\ol{r}(v\ten w)\]
and
\[\langle (1_+\ten 1)^{\ten2},\Omega_{D,23}\sum_{i,j}f_i\ten v_i\ten g_j\ten w_j\rangle=-\Omega_D(v\ten w)\]
where $\Omega=\ol{\Omega}+\Omega_D$, we get
\[\alpha_{V\ten W}\circ J_{\Gamma}\circ(\alpha_V^{-1}\ten\alpha_W^{-1})(v\ten w)\equiv v\ten w + \frac{\hbar}{2}(r+r_D^{21})(v\ten w)\mod\hbar^2\]
because the definition of $J_{\Gamma}$ involves the braiding $\beta'_{XY}=\beta^{-1}_{YX}$.
\end{pf}

%In analogy with \cite{vtl-4}, we get the following
\begin{corollary}\label{cor:twist-alt}
The relative twist $J_{\Gamma}$ satisfies
\[\operatorname{Alt}_2J_{\Gamma}\equiv
\frac{\hbar}{2}\left(\frac{r-r^{21}}{2}-\frac{r_D-r_D^{21}}{2}\right)\mod\hbar^2\]
\end{corollary}

%-----------------------------------
% QUANTUM VERMA
%-----------------------------------

\section{Quantization of Verma modules}\label{s:Gammah}

This section and the next contain results about the quantization of classical Verma modules,
which are required to construct the morphism of $D$--categories between the representation
theory of ${\Ueh{\g}}$ and that of ${\DJ{\g}}$. In particular, from now on, we will assume the
existence of a finite $\IN$--grading on $\g$, which induces on $\g$ a Lie bialgebra structure and
allows us to consider the quantization of $\g$ through the Etingof--Kazhdan functor, $\qEK{\g}$.

\subsection{Quantum Verma Modules}
%-------------------------------------------------
Because of the functoriality of the quantization defined by Etingof and Kazhdan in \cite{ek-2}, in the category of Drinfeld-Yetter modules over $\EK\gm$ we can similarly define \emph{quantum Verma modules}.\\

The standard inclusions of Lie bialgebras $\gpm\subset\g\simeq\D\gm$ lift to $\EK\gpm\subset\EK\g\simeq D\EK\gm$, and we can define the induced modules of the trivial representation over $\EK\gpm$
$$M_{\pm}^{\hbar}=\Ind_{\EK\gpm}^{\EK\g}\sfk[[\hbar]]$$
Similarly, we have Hopf algebra maps $\EK\Lppm\subset\EK\g$ and $\EK\Lppm\to\EK\gD$, and we can define induced modules
\[L_-^{\hbar}=\Ind_{\EK\Lpp}^{\EK\g}\sfk[[\hbar]]\qquad N_+^{\hbar}=\Ind_{\EK\Lpm}^{\EK\g}\EK\gD\]
We want to show that the equivalence $\wt{F}:\TDY{U\gm}{\Phi}[[\hbar]]\to\DY_{\EK\gm}$ matches these modules. We start proving the statement for $\Mm,\Mpd$.

\subsection{Quantization of $M_{\pm}$}
%--------------------------------------------------
We denote by $(\Mp^{\hbar})^*$ the $\qEK{\g}$--module
\[\Hom_{\sfk}(\Ind_{\EK\gm}^{\EK\g}\sfk[[\hbar]],\sfk[[\hbar]])\]
\begin{theorem}\label{th:quantumverma}
In the category of left $\qEK{\g}$--modules,
\begin{itemize}
\item[(a)] $F(\Mm)\simeq\Mm^{\hbar}$\\
\item[(b)] $F(\Mpd)\simeq(\Mp^{\hbar})^*$
\end{itemize}
\end{theorem}
\begin{pf}
The Hopf algebra $\EK\gm$ is constructed on the space $F(\Mm)$ with unit element $u\in F(\Mm)$ defined by $u(\um)=\coup\ten\um$, where $\coup\in\Mpd$ is defined as $\coup(x1_+)=\cou(x)$ for any $x\in U\gp$. Consequently, the action of $\EK\gm$ on $u\in F(\Mm)$ is free, as multiplication with the unit element. The coaction of $\EK\gm$ on $F(\Mm)$ is defined using the $\R$-matrix associated to the braided tensor functor $F$,\ie
$$\pi^*_{\Mm}:F(\Mm)\to F(\Mm)\ten F(\Mm),\quad \pi^*(x)=\R(u\ten x)$$
where $x\in F(\Mm)$ and $\R_{VW}\in\End_{\EK\g}(F(V)\ten F(W))$ is given by $\R_{VW}=\sigma J_{WV}^{-1}F(\beta_{VW})J_{VW}$, $\{J_{V,W}\}_{V,W\in\DY_{U\gm}}$ being the tensor structure on $F$. It is easy to show that
$J(u\ten u)|_{\um}=\coup\ten\um\ten\um$, and, since $\Omega(\um\ten\um)=0$, we have
$$\R(u\ten u)=u\ten u$$
For a generic $V\in\DY_{U\gm}[[\hbar]]$, the action of $\EK\gm^*$ is defined as 
$$F(\Mm)^*\ten F(V)\to F(\Mm)^*\ten F(\Mm) \ten F(V)\to F(V)$$
This means, in particular that, for every $\phi\in I\subset\EK\gm^*$, where $I$ is the maximal ideal corresponding to $u^{\perp}$, we have $\phi.u=0$. This proves $(a)$.\\

The module $\Mpd$ satisfies the following universal property: for any $V$ in the Drinfeld category of \emph{equicontinuous} $U\g$-modules, we have
$$\Hom_{U\g}(V,\Mpd)\simeq\Hom_{U\gm}(V,\sfk)$$
Indeed, to any map of $U\g$-modules $f:V\to\Mpd$, we can associate $\hat{f}:V\to\sfk$, $\hat{f}(v)=\langle f(v),\up\rangle$. It is clear that $\hat{f}$ factors through ${V}/{\gm.V}$. The equicontinuity property is necessary to show the continuity of $\hat{f}$ with respect to the topology on $V$.\\
Since $F$ defines an equivalence of categories, we have
$$\Hom_{\EK\g}(F(V),F(\Mpd))\simeq\Hom_{U\g}(V,\Mpd)[[\hbar]]\simeq\Hom_{U\gm}(V,\sfk)[[\hbar]]$$
Using the natural isomorphism $\alpha_V:F(V)\to V[[\hbar]]$, defined by $$\alpha_V(f)=\langle f(\um), \up\ten\id\rangle$$
we obtain a map $\Hom_{U\gm}(V,\sfk)[[\hbar]]\to\Hom_{\sfk}(F(V),\sfk[[\hbar]])$. Consider now the linear isomorphism $\alpha: \EK\gm\to U\gm[[\hbar]]$ and for any $x\in U\gm$ consider the $\g$-intertwiner $\psi_x:\Mm\to\Mpd\ten\Mm$ defined by $\psi_x(\um)=\coup\ten x\um$.
It is clear that, if $f(\um)=f_{(1)}\ten f_{(2)}$ in Swedler's notation,
\begin{eqnarray*}
\alpha_V(\psi_x.f)&=&\langle (\ipd\ten\id)\Phi^{-1}(\id\ten f)(\coup\ten x.\um), \up\ten\id \rangle\\
&=&\langle \Phi^{-1}(\coup\ten\id\ten\id)(\id\ten\Delta(x))(\id\ten f_{(1)}\ten f_{(2)}), (T\ten \id)(\up\ten\up\id) \rangle\\
&=&\langle \Delta(x)(f_{(1)}\ten f_{(2)}), \up\ten\id\rangle\\
&=&\langle f_{(1)},\up\rangle x.f_{(2)}\\
&=&x.\alpha_V(f)
\end{eqnarray*}
using the fact that $(\epsilon\ten1\ten1)(\Phi)=1^{\ten 2}$ and $(\epsilon\ten 1)(T)=1$.
So, clearly, if $\phi\in\Hom_{U\gm}(V,\sfk)$, then $\phi\circ\alpha_V\in\Hom_{\EK\gm}(F(V),\sfk[[\hbar]])$. Then $F(\Mpd)$ satisfies the universal property of $\Hom_{\sfk}(\Ind_{\EK\gm}^{\EK\g}\sfk[[\hbar]],\sfk[[\hbar]])$ and $(b)$ is proved.
\end{pf}

% THE RELATIVE CASE
\subsection{Quantization of relative Verma modules}
The proof of $(b)$ shows that the linear functional $F(\Mpd)\to\sfk[[\hbar]]$ is, in fact, the trivial deformation of the functional $\Mpd\to\sfk$. These results extend to the relative case and hold for the right $\gD$--action on $\Lm,\Npd$.

\begin{theorem}\label{th:quantization-LN} In the category $\DY_{\qEK{\gm}}$
\begin{itemize}
\item[(a)] $F(\Lm)\simeq\Lmh$\\
\item[(b)] $F(\Npd)\simeq(\Np^{\hbar})^*$
\end{itemize}
Moreover, as right $\qEK{\gD}$--module
\begin{itemize}
\item[(c)] $F_D(\Lm)\simeq\Lmh$\\
\item[(d)] $F_D(\Npd)\simeq(\Np^{\hbar})^*$
\end{itemize}
\end{theorem}

The proof of $(a)$ and $(b)$ amounts to constructing the morphisms
\[\sfk[[\hbar]]\to F(\Lm)\qquad F(\Npd)\to\qEK{\gD}^*\]
equivariant under the action of $\qEK{\Lpp}$ and $\qEK{\Lpm}$ respectively.\\

A direct construction along the lines of the proof of Theorem \ref{th:quantumverma} is however not straightforward. We prove this theorem in the next section by using a description of the modules $\Lm, \Npd$ and their images through $\wt{F}$ via \PROP categories. These descriptions show that the classical intertwiners
\[\sfk\to\Lm\qquad\Npd\to U\gD^*\]
satisfy the required properties and yield canonical identifications
\[\wt{F}(\Lm)\simeq\Lmh\qquad \wt{F}(\Npd)\simeq(\Np^{\hbar})^*\]

%------------------------------
% PROPS
%------------------------------

\section{Universal relative Verma modules}\label{s:props}
%--------------------------------------------------------------------------------

In this section, we prove Theorem \ref{th:quantization-LN} by using suitable
\PROP (\emph{product-permutation}) categories compatible with the EK
universal quantization functor \cite{ek-2, EG}.

\subsection{\PROP description of the EK quantization functor}
%--------------------------------------------------------------------------------

We will briefly review the construction of Etingof--Kazhdan in the setting of \PROP categories
\cite{ek-2}.

A \PROP is a symmetric tensor category generated by one object. More precisely, a \emph{
cyclic category over $S$} is the datum of
\begin{itemize}
\item a symmetric monoidal $\sfk$--linear category $(\C,\ten)$ whose objects are non--negative
integers, such that $[n]=[1]^{\ten n}$ and the unit object is $[0]$
\item a bigraded set $S=\bigcup_{m,n\in\IZ_{\geq0}}S_{nm}$ of morphism of $\C$, with
\[S_{nm}\subset\Hom_{\C}([m],[n])\]
\end{itemize}  
such that any morphism of $\C$ can be obtained from the morphisms in $S$ and permutation 
maps in $\Hom_{\C}([m],[m])$ by compositions, tensor products or linear combinations over $\sfk$.
We denote by $\F_S$ the free cyclic category over $S$. Then there exists a unique symmetric tensor
functor $\F_S\to\C$, and the following holds (cf.~\cite{ek-2})

\begin{proposition}
Let $\C$ be any cyclic category generated by a set $S$ of morphisms. Then $\C$ has the form $\F_S/\I$, where $\I$ is a tensor ideal in $\F_S$.
\end{proposition}

Let $\N$ be a symmetric monoidal $\sfk$--linear category, and $X$ an object in $\N$. A \emph{linear 
algebraic structure of type $\C$ on $X$} is a symmetric tensor functor $\G_X:\C\to\N$ such that $\G_{X}([1])=X$.
 A linear algebraic structure of type $\C$ on $X$ is a collection of morphisms between tensor powers of $X$ 
 satisfying certain consistency relations.

We mainly consider the case of \emph{non--degenerate} cyclic categories, \ie symmetric tensor categories
 with injective maps $k[\SS_n]\to\Hom_{\C}([n],[n])$.  We first consider the Karoubian envelope of $\C$ 
 obtained by formal addition to $\C$ of the kernel of the idempotents in $k[\SS_n]$ acting on $[n]$. 
 Furthermore, we consider the closure under inductive limits. In this category, denoted $S(\C)$, every 
 object is isomorphic to a direct sum of indecomposables, corresponding to irreducible representations 
 of $\SS_n$ (cf.~\cite{ek-2,EG}). In particular, in $S(\C)$, we can consider the symmetric algebra 
\[S[1]=\bigoplus_{n\geq0} S^n[1]\] 
If $\N$ is closed under inductive limits, then any linear algebraic structure of type $\C$ extends to an 
additive symmetric tensor functor
\[\G_X:S(\C)\to\N\]

We introduce the following fundamental \PROPs.
\begin{itemize}
\item {\bf Lie bialgebras.} In this case the set $S$ consists of two elements of bidegrees $(2,1),(1,2)$, the universal commutator and cocommutator. The category $\C=\LBA$ is $\F_S/\I$, where $\I$ is generated by the classical five relations.\\
\item {\bf Hopf algebras.} In this case, the set $S$ consists of six elements of bidegrees $(2,1),(1,2),(0,1),(1,0),(1,1,),(1,1)$, the universal product, coproduct, unit, count, antipode, inverse antipode. The category $\C=\HA$ is $\F_S/\I$, where $\I$ is generated by the classical four relations. 
\end{itemize}

The quantization functor described in Section \ref{s:ek} can be described in this generality, as stated by the following (cf.~\cite[Thm.1.2]{ek-2})

\begin{theorem}
There exists a universal quantization functor $Q:\HA\to S(\LBA)$.
\end{theorem}

Let $\gm$ be the canonical Lie dialgebra $[1]$ in $\LBA$ with commutator $\mu$ and cocommutator $\delta$. Let $U\gm:=S\gm\in S(\LBA)$ be the universal enveloping algebra of $\gm$. The construction of the Etingof--Kazhdan quantization functor amounts to the introduction of a Hopf algebra structure on $U\gm$,  which coincides with the standard one modulo $\langle\delta\rangle$, and yields the Lie bialgebra structure on $\gm$ when considerd modulo $\langle\delta^2\rangle$. This Hopf algebra defines the object $Q[1]$, where $[1]$ is the generating object in $\HA$. The formulae used to defined the Hopf structure coincide with those defined in \cite[Part II]{ek-1} and described in Section \ref{s:ek}. In particular, they rely on the construction of the Verma modules 
\[\Mm:=S\gm \Mpd=\wh{S\gm}\]
realized in the category of Drinfeld--Yetter modules over $\gm$ as object of \LBA~. 

\subsection{\PROPs for split pairs of Lie bialgebras}

Let $(\gm,\gDm)$ be a \emph{split} pair of Lie bialgebras, \ie there are Lie bialgebra maps
\[\gDm\xrightarrow{i}\gm\xrightarrow{p}\gDm\]
such that $p\circ i=\id$. These maps induce an inclusion $\Db\gDm\subset\Db\gm$ and consequently 
an inclusion of Manin triple $(\gD,\gDm,\gDp)\subset(\g,\gm,\gp)$, as described in Section \ref{s:parabolic}.

\begin{definition}
We denote by \PLBA the Karoubian envelope of the multicolored {\sc Prop}, whose class of objects is 
generated by the Lie bialgebra objects $[\gm],[\gDm]$, related by the maps $i:[\gDm]\to[\gm]$ ,
$p:[\gm]\to[\gDm]$, such that $p\circ i=\id_{[\gDm]}$. 
\end{definition}

The Karoubian envelope implies that $[\Lmm]:=\ker(p)\in\PLBA$.

\begin{proposition}
The multicolored \PROP \PLBA is endowed with a pair of functors $U,L$
\[U,L:\LBA\to\PLBA\qquad U[1]:=[\gm],\quad L[1]:=[\gDm]\]
and natural transformations $i,p$, induced by the maps $i,p$ in $\PLBA$,
\[
\xymatrix{\LBA\ar@/^15pt/[rr]^U="U"\ar@/_15pt/[rr]_L="L"  & & \PLBA \ar@{=>}@<1ex> "U";"L"^p \ar@{=>}@<1ex> "L";"U"^i}
\]
such that $p\circ i=\id$. Moreover, it satisfies the following universal property: for any tensor category $\C$, 
closed under kernels of projections, with the same property as \PLBA, there exists a unique tensor functor 
$\PLBA\to\C$ such that the following diagram commutes\\ 
\[
\xymatrix{\LBA\ar@/^15pt/[rr]^U="U"\ar@/_15pt/[rr]_L="L"  \ar@/^30pt/[rrr] \ar@/_30pt/[rrr] & & \PLBA \ar@{=>}@<1ex> "U";"L"^p \ar@{=>}@<1ex> "L";"U"^i  \ar@{-->}[r] & \C}
\]
\end{proposition}

\subsection{\PROPs for split pairs of Hopf algebras}
We can analogously define suitable \PROP categories corresponding to split pairs of Hopf algebras. 
In particular, we consider the \PROP \PHA characterized by functors $U_{\hbar},L_{\hbar}$ and natural 
transformations $p_{\hbar}, i_{\hbar}$ satisfying
\[
\xymatrix{\HA\ar@/^15pt/[rr]^{U_{\hbar}}="U"\ar@/_15pt/[rr]_{L_{\hbar}}="L"  \ar@/^30pt/[rrr] \ar@/_30pt/[rrr] & & \PHA \ar@{=>}@<1ex> "U";"L"^{p_{\hbar}} \ar@{=>}@<1ex> "L";"U"^{i_{\hbar}}  \ar@{-->}[r]^{\exists!}& \C}
\]
where \HA denotes the \PROP category of Hopf algebras. These also satisfy
\[
\xymatrix@C=0.7in{\HA\ar[r]^{Q\oEK}\ar@<1ex>[d] \ar@<-1ex>[d] & \hLBA \ar@<1ex>[d] \ar@<-1ex>[d] \\ \PHA \ar[r]^{Q_{\PLBA}} & \hPLBA}
\]
where $Q_{\PLBA}$ is the extension of the Etingof--Kazhdan quantization functor to $\PLBA$, obtaine by the  universal property described above with $\C=S(\PLBA)$. 

\subsection{\PROPs for parabolic Lie subalgebras}
In order to describe the module $\Npd$ it is necessary to deal with the Lie bialgebra object $\Lpm$ or, in other words to introduce the double of $\gDm$ and the \PROP $D_{\oplus}(\LBA)$ \cite{EG}. We then introduce the multicolored \PROP as a cofiber product of \PLBA and $D_{\oplus}(\LBA)$ over \LBA. 

\begin{proposition}
The multicolored \PROP \PLBAD is endowed with canonical functors 
\[
D_{\oplus}(\LBA)\to\PLBAD\leftarrow\PLBA
\]
and satisfies the following universal property:
\[
\xymatrix@C=0.8in{\LBA \ar[r]^{\operatorname{double}} \ar[d] & D_{\oplus}(\LBA) \ar[d] \ar@/^15pt/[ddr]& \\ \PLBA\ar@/_15pt/[rrd] \ar[r] & \PLBAD \ar@{-->}[dr]^{\exists!} & \\ && \C}
\]
where $\operatorname{double}$ is the \PROP map introduced in \cite{EG}.
\end{proposition}

In \PLBAD we can consider the Lie bialgebra object $[\Lpm]$. 

\subsection{\PROPs for parabolic Hopf subalgebras} Similarly, we introduce the multicolored \PROP \PHAD, endowed with canonical functors (cf. \cite{EG})
\[D_{\ten}(HA)\to\PHAD\leftarrow\PHA\]
and satisfying an analogous universal property:
\[
\xymatrix@C=0.8in{\HA \ar[r]^{\operatorname{double}} \ar[d] & D_{\ten}(\HA) \ar[d] \ar@/^15pt/[ddr]& \\ \PHA\ar@/_15pt/[rrd] \ar[r] & \PHAD \ar@{-->}[dr]^{\exists!} & \\ && \C}
\]

Moreover, we then have a canonical functor
\[
Q_{\PLBAD}: \PHAD\to\hPLBAD
\]
obtained applying such universal property with $\C=S(\PLBAD)$ and satisfying
\[
\xymatrix{& \HA\ar[dl]^{L_{\HA}} \ar[rr]^{\operatorname{double}}  \ar[ddd]^{Q\oEK}& &  D_{\otimes}(\HA) \ar[ddd]^{Q_2} \ar[dl] \\ \PHA\ar[ddd]^{Q_{\PLBA}}\ar[rr]& & \PHAD \ar[ddd]^{Q_{\PLBAD}} & \\ &&&\\ & \hLBA \ar[dl]\ar[rr]^(.3){S(\operatorname{double})} & & S(D_{\oplus}(\LBA)) \ar[dl] \\ \hPLBA \ar[rr] & & \hPLBAD & }
\]
The commutativity of the square on the back is given by the compatibility of the quantization functor with the doubling operations, proved in \cite{EG}.

\subsection{\PROP description of $\Lm,\Npd$} 

The modules $\Lm,\Npd$ can be realized in $S(\PLBAD)$. The module $\Lm$ is constructed over the object $S\Lmm\in S(\PLBA)$. The structure of Drinfeld-Yetter module over $\gm$ is determined in the following way:
\begin{itemize}
\item the free action of the Lie algebra object $\Lmm$ is defined by the map
\[S\Lmm\ten S\Lmm\to S\Lmm\]
given by Campbell-Hausdorff series, describing on $S\Lmm$ the multiplication
in $U\Lmm$.\\
\item we define the action of $\gDm$ to be trivial on ${\bf 1}\to S\Lmm$.\\ 
\item The actions of $\Lmm,\gDm$, the relation
\[\pi\circ([,]\ten1)=\pi\circ(1\ten\pi)-\pi\circ(1\ten\pi)\circ\sigma_{12}\]
and the map $[,]:\gDm\ten\Lmm\to\Lmm$ define the action of $\gm$.\\
\item We then impose the trivial coaction on ${\bf 1}\to S\Lmm$ and the compatibility condition 
between action and coaction
\[
\quad\qquad\pi^*\circ\pi=(1\ten\pi)\sigma_{12}(1\ten\pi^*)-(1\ten\pi)(\delta\ten1)+(\mu\ten 1)(1\ten \pi^*)
\]
determines the coaction for $S\Lmm$. The action defined is compatible with $[,]:\gDm\ten\Lmm\to\Lmm$
\end{itemize}

Similarly, the module $\Npd$ can be realized on the object $\wh{S\Lpm}$, formally 
added to $S(\PLBAD)$.\\

We determine the formulae for the action and the coaction of $\gm$ by direct inspection
of the action of $\g=\gm\oplus\gp$ on $\Np$ in the category $\vect$. Namely, the 
identification $\Npd=\wh{S\Lpm}$ is clearly obtained through the invariant bilinear form $\iip{-}{-}$
and there are topological formulae expressing the action of $\g$ on $\Np$. 
Therefore we determine action and coaction on $\Npd$ in the following way:
\begin{itemize}
\item the $\gp$--action on $\Np=S\Lpp=S\gp\ten S\g_{D,-}$ is given by the free action on the first factor 
$S\gp$ expressed by Campbell--Hausdorff series.\\
\item the action of $\gm=\Lmm\oplus\g_{D,-}$ on the subspace $S\g_{D,-}\subset S\Lpp$ is given by
the trivial action of $\Lmm$ and the usual free action of $\g_{D,-}$ by multiplication.\\
\item The action of $\gm$ is then interpreted as a topological coaction of $\gp$ and the aforementioned 
compatibility condition between action and coaction allows to extend the formula for the topological $\gp$ 
coaction on the entire space $S\Lpp$.\\
\item Through the invariant bilinear form $\iip{-}{-}$, these formulae are carried over $\Npd=\wh{S\Lpm}$,
by switching, in particular, the bracket and the topological cobracket on $\gp$ with the cobracket
and the bracket in $\gm$, respectively.\\
\item The obtained formulae, describing the action and the coaction of $\gm$ on $\Npd$, are well--defined
in the category $\PLBAD$ and define the requested structure of Drinfeld-Yetter module over $\gm$.
\end{itemize}

\subsection{Proof of Theorem \ref{th:quantization-LN}}

The {relative} Verma module 
\[\Np=\ind_{\Lmm}^{\g}\sfk\simeq\ind_{\Lpm}^{\g}U\gD\] 
satisfies 
\[\Hom_{U\g}(\Np, V)\simeq\Hom_{U\Lpm}(U\gD, V)\]
for every $U\g$-module $V$. We have a canonical map of $\Lpm$-modules $\rho_D:U\gD\to\Np$ corresponding to the identity in the case $V=\Np$. We get a map of $\Lpm$-modules $\rho_D^*:\Npd\to U\gD^*$ inducing an isomorphism
\[\Hom_{U\g}(V,\Npd)\simeq\Hom_{U\Lpm}(V,U\gD^*)\]
The morphism $\rho_D^*$ can indeed be thought as 
\[
\xymatrix{U\Lpm\ten\Npd\ar@<1ex>[d]\ar@<-1ex>[d]\ar[r] & \Npd \ar[d]^{\rho_D^*} \\ U\gD\ten U\gD^*\ar[r] & U\gD^*}
\]

Assuming the existence of a suitable finite $\IN$--grading, a split pair of Lie bialgebras $(\gm,\gDm)$,  gives rise to a functor
\[\PLBAD\to\vect\]
Consider now the trivial split pair given by $(\gDm,\gDm)$. We have a natural transformation
\[
\xymatrix{\PLBAD\ar@/^15pt/[rr]^{(\gm,\gDm)}="U"\ar@/_15pt/[rr]_{(\gDm,\gDm)}="L"  & & \vect \ar@{=>}"U";"L"^p}
\]
where $p$ naturally extends to the projection $\Lpm\to\gD$.\\ 

The module $U(\gD)^*$ is indeed the module $\Npd$ with respect to the trivial pair $(\gDm,\gDm)$. Consequently, the existence of the $\Lpm$--intertwiner $\rho_D^*$ can be interpreted as a simple consequence of the existence of natural transformation $p$.\\

The quantization functor $Q_{\PLBAD}$  extends the natural transformation $p$ to
\[
\xymatrix{&&\\\PHAD\ar@/^40pt/[rrr]\ar@/_40pt/[rrr]\ar[r]&S(\PLBAD)\ar@/^15pt/[rr]^{(\gm,\gDm)}="U"\ar@/_15pt/[rr]_{(\gDm,\gDm)}="L"  & & \vect \ar@{=>}"U";"L"^{S(p)}\\ &&}
\]
and shows that 
\[
F(\Npd)\simeq (N_+^{\hbar})^*
\]

Similarly, we can consider the natural transformation $S(i)$ and the diagram
\[
\xymatrix{&&\\ \PHAD\ar@/^40pt/[rrr]\ar@/_40pt/[rrr]\ar[r]&S(\PLBAD)\ar@/^15pt/[rr]^{(\gm,\gDm)}="U"\ar@/_15pt/[rr]_{(\gDm,\gDm)}="L"  & & \vect \ar@{=>}"L";"U"^{S(i)}\\&&}
\]
implying
\[
F(\Lm)\simeq\Lmh
\]
We can make analogous consideration for the right $\gD$--action on $\Lm,\Npd$. This leads to isomorphisms of right $\qEK{\gD}$--modules
\[
\wt{F}_D(\Npd)\simeq\Npdh\qquad F_D\oEK(\Lm)\simeq\Lmh
\]

%
% MANIN CHAINS
%

\section{Chains of Manin triples}\label{s:mchain}
%========================

%In this section, we apply the previous results to chains of Manin triples.  

\subsection{Chains of length $2$}
%----------------------------------------

In Section \ref{s:Gammah}, given an inclusion of Manin triples $i_D:\gD\subseteq\g$,
we introduced the relative quantum Verma modules
\[L_-^{\hbar}=\Ind_{\EK\Lpp}^{\EK\g}\sfk[[h]]\qquad N_+^{\hbar}=\Ind_{\EK\Lpm}^{\EK\g}\EK\gD\]
These modules allow to define the functor
\[\Gamma_{\hbar}:\DC{\qEK{\g}}{}\to\DC{\qEK{\gD}}{}\]
by
\[\Gammah(\V)=\Hom_{\qEK{\g}}(\Lmh,\Npdh\ten\V)\]
\begin{lemma}\label{lem:hGamma-res}
The functor $\Gammah$ is naturally tensor isomorphic to the 
restriction functor $(\qEK{(i_D^{\hbar})})^*$.
\end{lemma}
\begin{pf}
The proof of the existence of the natural isomorphism as $\qEK{\gD}$--module 
is identical to that of Proposition \ref{pr:isom to forget}. The isomorphism respects
the tensor structures, because there are only trivial associators involved.
\end{pf}

\subsection{}
%--------------

We now prove the following
\begin{theorem}\label{thm:ek-fac}
Let ${\g}, \gD$ be Manin triples with a finite $\IZ$--grading and $i_D:\gD\subseteq\g$ an
inclusion of Manin triples compatible with the grading.
Then, there exists an algebra
isomorphism
\[\Psi:\wh{\qEK{{\g}}}\to\wh{\Ueh{{\g}}}\]
restricting to $\Psi\oEK_D$ on $\wh{\qEK{\gD}}$, where the completion is given \wrt
Drinfeld--Yetter modules.
\end{theorem}

\begin{pf}
In the previous section, we showed that the quantization of the $(\qEK{\g},\qEK{\gD})$--modules
$\Npd,\Lm$ gives
\[
\wt{F}(\Npd)\xrightarrow{\qEK{\g}}\Npdh\xleftarrow{\qEK{\gD}}F_D\oEK(\Npd)
\]
\[
\wt{F}(\Lm)\xrightarrow{\qEK{\g}}\Lmh\xleftarrow{\qEK{\gD}}F_D\oEK(\Lm)
\]
Recall that the standard natural transformations $\alpha_V:\wt{F}(V)\simeq V[[\hbar]]$, $(\alpha_D)_V: \wt{F}_D(V)\simeq V[[\hbar]]$ give isomorphisms of right $\Ueh{\gD}$--modules
\[\wt{F}(\Npd)\simeq\Npd[[\hbar]]\qquad \wt{F}(\Lm)\simeq\Lm[[\hbar]]\]
and isomorphisms of $\Ueh{\g}$--modules
\[F_D\oEK(\Npd)\simeq\Npd[[\hbar]]\qquad F_D\oEK(\Lm)\simeq\Lm[[\hbar]]\]
In particular, we get isomorphisms of right $\qEK{\gD}$--modules
\[
F_D\oEK\circ \wt{F}(\Npd)\simeq F_D\oEK(\Npd)\simeq\Npdh\qquad F_D\oEK\circ \wt{F}(\Lm)\simeq F_D\oEK(\Lm)\simeq\Lmh
\]
and isomorphisms of $\qEK{\g}$--modules
\[
F_D\oEK\circ \wt{F}(\Npd)\simeq \wt{F}(\Npd)\simeq\Npdh\qquad F_D\oEK\circ \wt{F}(\Lm)\simeq \wt{F}(\Lm)\simeq\Lmh
\]
We have a natural isomorphism through $J$:
\[
\Hom_{\qEK{\g}}(F(\Lm),F(\Npd)\ten F(V))\simeq\Hom_{\g}(\Lm,\Npd\ten V)[[\hbar]]
\]
This is indeed an isomorphism of $\Ueh{\gD}$--modules, since, for $x\in U\gD$, $\phi\in\Hom_{\qEK{\g}}(F(\Lm),F(\Npd)\ten F(V))$, we have
\[
x.\phi:=(F(x)\ten\id)\circ\phi\qquad J\circ(F(x)\ten\id))=F(x\ten\id)\circ J
\]
Quantizing both side and using the isomorphism $F_D\oEK\circ F(\Npd)\simeq\Npdh$, we obtain a natural transformation
\[
\gamma_D:\Gammah\circ \wt{F}\simeq \wt{F}_D\circ\Gamma
\]
making the following diagram commutative
\[
\xymatrix{\DC{U{\g}}{\Phi}\ar[r]^{\wt{F}} \ar[d]_{\Gamma} & \DC{\DJ{\g}}{}\ar[d]^{\Gamma^{\hbar}} \ar@{=>}[dl]_{\gamma_D} \\ \DC{\Ueh{\gD}}{\Phi_{D}} \ar[r]^{\wt{F}_D}& \DC{\DJ{\gD}}{}}
\]

Applying the construction above to the algebra of endomorphisms of the fiber functor, we get the result.
\end{pf}

\subsection{Chains of arbitrary length}
%----------------------------------------------

For any chain 
\[{\bf C}: 0=\g_0\subseteq\g_1\subseteq\cdots\subseteq\g_{n-1}\subseteq\g_n=\g\]
of inclusions of Manin triples, the natural transformations 
\[\gamma_{i,i+1}\in\sfNat_{\ten}(\Gamma^{\hbar}_{i,i+1}\circ \wt{F}_{i+1},\wt{F}_i\circ\Gamma_{i,i+1})\] 
where $0\leq i\leq n-1$, $\Gamma_{0,1}:\D_\Phi(\g_1)\to\kvect$ is the EK fiber functor,
and $F\EEK_0=\id$, yield a natural transformation 
\[\begin{split}
\gamma_{\bf C}
&=
\gamma_{0,1}\circ\cdots\circ\gamma_{n-1,n}\\
&\in\sfNat_{\ten}
(\Gamma^\hbar_{0,1}\circ\cdots\circ\Gamma^\hbar_{n-1,n}\circ \wt{F}_n,
\Gamma_{0,1}\circ\cdots\circ\Gamma_{n-1,n})\\
&\cong\sfNat_{\ten}
((i_{0,n}^*)^\hbar\circ \wt{F}_n,
\Gamma_{0,1}\circ\cdots\circ\Gamma_{n-1,n})\\
&=\sfNat_{\ten}
(\wt{F}_n,\Gamma_{0,1}\circ\cdots\circ\Gamma_{n-1,n})
\end{split}\]
where we used $\Gamma_{i,i+1}^\hbar\cong(i_{i,i+1}^*)^\hbar$, and the fact that
the composition $(i_{0,n}^*)^\hbar\circ \wt{F}_n$ is the EK fiber functor for $\g_n$,
which we denote by the same symbol as $\wt{F}_n$.

This proves the following
\begin{theorem}\hfill
\begin{enumerate}
\item
For any chain of Manin triples
\[{\bf C}: \g_0\subseteq\g_1\subseteq\cdots\subseteq\g_n\subseteq\g\]
there exists an isomorphism of algebras
\[\Psi_{\bf C}:\wh{\qEK{{\g}}}\to\wh{\Ueh{{\g}}}\] 
such that $\Psi_{\bf C}(\wh{\qEK{{\g_{i}}}})
=\wh{\Ueh{\g_{i}}}$ for any $\g_i\in{\bf C}$.
\item
Given two chains ${\bf C,C'}$, the natural transformation 
\[\Phi_{\bf C C'}:=\gamma_{\bf C}^{-1}\circ\gamma_{\bf C'}\in\sfAut{F\EEK}\]
satisfies
\[\sfAd{\Phi_{\bf CC'}}\Psi_{\bf C'}=\Psi_{\bf C}\]
\end{enumerate}
\end{theorem}

\begin{proposition}
The natural transformations $\{\Phi_{\bf CC'}\}_{\bf C, C'}$ satisfy the 
following properties
\begin{itemize}
\item[(i)] {\bf Orientation.} Given two chains $\bf C, C'$
\[\Phi_{\bf CC'}=\Phi_{\bf C'C}^{-1}\]
\item[(ii)] {\bf Transitivity.} Given the chains $\bf C, C', C''$ 
\[\Phi_{\bf CC'}\circ\Phi_{\bf C'C''}=\Phi_{\bf CC''}\]
\item[(iii)] {\bf Factorization.} Given the chains
\[{\bf C},{\bf C'}: \g_0\subseteq\g_0\subseteq\cdots\subseteq\g_n
\qquad
{\bf D}, {\bf D'}:\g_n\subseteq\cdots\subseteq\g_{n+n'}\]
\[\Phi_{({\bf C}\cup{\bf D})({\bf C'}\cup{\bf D'})}=
\Phi_{{\bf C}{\bf C'}}\circ\Phi_{{\bf D}{\bf D'}}\]
\end{itemize}
\end{proposition}

%-----------------------------------------------------------------------------
% ABELIAN SUBALGEBRAS AND CENTRAL EXTENSION
%-----------------------------------------------------------------------------

\subsection{Abelian Manin triples and central extensions} 
%===========================================

We will now consider the following special case, that generalizes the role
of Levi subalgebras for Kac--Moody algebras.

\begin{proposition}\label{pr:ab-ext}
If $\g$ admits a Manin subtriple $\ll_D$, obtained by a central extension of
$\g_D$, then the relative twists and the gauge transformations are invariant
under $\ll_D$. In particular, the Etingof--Kazhdan constructions are invariant
under abelian Manin subtriples.
\end{proposition}

\begin{pf}
For $\gD=\{0\}$, the statement reduces to prove that the Etingof--Kazhdan functor preserves
the action of an abelian Manin subtriple $\a\subset\g$ (cf.~\cite[Thm. 4.3]{ek-6}, with $\a_-=\h$). Under this assumption, the natural map 
\[\xymatrix@R=.1in{U\a_-\ar[r] & \EK{\gm}:=\wt{F}(\Mm)\\ 
a \ar@{|->}[r] &\{\psi_a: 1_-\mapsto \varepsilon\ten a1_-\}}\] 
defines an inclusion of bialgebras. For any $V[[\hbar]]\in\DC{U{\g}}{\Phi}$, the natural identification
\[\alpha_V:F(V)\to V[[\hbar]]\]
is then an isomorphism of $U\a$--modules. This gives the following commutative diagram
\[\xymatrix{
\DC{U{\g}}{\Phi} \ar[rr]^{\wt{F}} \ar[dr] \ar@/_15pt/[ddr]_{F} & & \DC{\EK{\g}}{} \ar@/^15pt/[ddl] \ar[dl] \\
 & \DC{\Ueh{\a}}{} \ar[d] & \\
  & \A & }
\]
We can observe that the tensor restriction functor fits in an analogous diagram. 
It is easy to show that the object $\Gamma_D(\Lm)$ is naturally a pointed Hopf
algebra in the category $\DC{U\gD}{\Phi_D}$. We denote by 
$\D_{\gD}(\Gamma_D(\Lm))$ the category of Drinfeld--Yetter modules over 
$\Gamma_D(\Lm)$ in the category $\DC{U\gD}{\Phi_D}$. This category is 
naturally equivalent to the category of Drinfeld--Yetter module over the Radford's 
product $\Ueh{\gDm}\#\Gamma_D(\Lm)$ and there is a natural identification
\[
\xymatrix{
\DC{U{\g}}{\Phi} \ar[dr]_{\Gamma_D}\ar[rr] && \D_{\gD}(\Gamma_D(\Lm))\ar[dl]\\
&\DC{U\gD}{\Phi_D}&
}
\]
Moreover, there is a natural inclusion of bialgebras
\[U\ll_D\subset\Ueh{\gD}\#\Gamma_D(\Lm)\]
and a natural $U\ll_D$--module identification $\Gamma_D(V)\to V[[\hbar]]$. This originates natural identifications
\[\xymatrix{
\DC{U{\g}}{\Phi} \ar[rr] \ar[dr] \ar@/_15pt/[ddr]_{\Gamma_D} & & \D_{\gD}(\Gamma_D(\Lm)) \ar@/^15pt/[ddl] \ar[dl] \\
 & \DC{U{\ll_D}}{} \ar[d] & \\
  & \DC{U\gD}{\Phi_D} & }
\]
This proves that relative twists are invariant under $\ll_D$. It is clear that the Casimir
operator $\Omega_D\in (\gD\ten\gD)^{\ll_D}$ defines a braided tensor structure on 
$\DC{\Ueh{\ll_D}}{}$ that is preserved by the restriction functor induced by the inclusion
$j_D:\gD\subset\ll_D$. Given the decomposition $\ll_D=\gD\rtimes\c_D$, the natural
map $U\c_D\to\End_{\gD}(j_D^*V)$ induces an action of $U\c_D$ on $\wt{F}_D(j_D^*V)$,
commuting with the action of $\qEK{\gD}$. Therefore, we obtain a naturally commutative
diagram
\[
\xymatrix{
\DC{U\ll_D}{\Phi_D}\ar[r]^{\wt{\wt{F}_D}} \ar[d]_{j_D^*} & \DC{\qEK{\ll_D}}{}\ar[d]\\
\DC{U\gD}{\Phi_D}\ar[r]^{\wt{F}_D} & \DC{\qEK{\gD}}{} 
}
\]
where $\wt{\wt{F}_D}$ is the tensor functor induced by the composition $\wt{F}_D\circ j_D^*$.
The natural transformation $\gamma$ automatically lifts to the level of $\ll_D$, as showed in 
the following diagram
\[
\xymatrix{
\DC{U{\g}}{\Phi}\ar[rrr]^{\wt{F}}\ar[dr] \ar[dd]_{\Gamma}&&& \DC{\qEK{\g}}{}\ar[dl]\ar[dd]\\
&\DC{U\ll_D}{\Phi_D}\ar[r]^{\wt{\wt{F}_D}} \ar[dl]_{j_D^*} & \DC{\qEK{\ll_D}}{}\ar[dr] & \\
\DC{U\gD}{\Phi_D}\ar[rrr]^{\wt{F}_D} &&& \DC{\qEK{\gD}}{} 
}
\]
\end{pf}

%--------------------
% FINAL SECTION CONCLUSION
%--------------------

\section{An equivalence of quasi--Coxeter categories}\label{s:qcstructure}
%=========================================

The following is the main result of this paper.

\begin{theorem}
Let $\g$ be a symmetrizable Kac--Moody algebra with a fixed $D_{\g}$--structure. Then the completion $\wh{\DJ{{\g}}}$ is {isomorphic} to a quasi-Coxeter quasitriangular quasibialgebra of type $D_{{\g}}$ on the quasitriangular $D_{{\g}}$--quasibialgebra 
\[(\wh{\Ueh{{\g}}},\{\wh{\Ueh{{\gD}}}\}, \Delta_0,\{\Phi\oKZ_D\},\{R\oKZ_D\})\]
where the completion is taken with respect to the integrable modules in category $\O$.
\end{theorem}

%-----------------------------
% D-STRUCTURE
%-----------------------------
\subsection{$D$--structures on Kac--Moody algebras}\label{ss:kmdstructure}

Let $\sfA=(a_{ij})_{i,j\in\bfI}$ be a complex $n\times n$ matrix and $\g=\g(\sfA)$ the corresponding generalized Kac--Moody 
algebra defined in Section \ref{s:ek}. Let $\bfJ$ be a nonempty subset of $\bfI$. Consider the submatrix of $\sfA$
defined by
\[\sfA_{\bfJ}=(a_{ij})_{i,j\in\bfJ}\] 

We recall the following proposition from \cite[Ex.1.2]{K}
\begin{proposition}
Let 
\[\Pi_{\bfJ}:=\{\alpha_j\;|\;j\in\bfJ\}\qquad\Pi^{\vee}_{\bfJ}:=\{h_j\;|\;j\in\bfJ\}\]
Let $\h_{\bfJ}^{\prime}$ be the subspace of $\h$ generated by $\Pi^{\vee}_{\bfJ}$ and 
\[\mathfrak{t}_{\bfJ}=\bigcap_{j\in\bfJ}\Ker\alpha_j=\{h\in\h\;|\;\iip{\alpha_j}{h}=0\;\forall j\in\bfJ\}\]
Let $\h^{''}_{\bfJ}$ be a supplementary subspace of $\h_{\bfJ}^{\prime}+\mathfrak{t}_{\bfJ}$ in $\h$ and let
\[\h_{\bfJ}=\h^{\prime}_{\bfJ}\oplus\h_{\bfJ}^{''}\]
Then,
\begin{itemize}
\item[(i)] $(\h_{\bfJ},\Pi_{\bfJ},\Pi_{\bfJ}^{\vee})$ is a realization of the generalized Cartan matrix $\sfA_{\bfJ}$.\\
\item[(ii)] The subalgebra $\g_{\bfJ}\subset\g$, generated by $\{e_j,f_j\}_{j\in\bfJ}$ and $\h_{\bfJ}$, is the Kac--Moody algebra associated to the realization $(\h_{\bfJ},\Pi_{\bfJ},\Pi_{\bfJ}^{\vee})$ of $\sfA_{\bfJ}$.
\end{itemize}
Set 
\[Q_{\bfJ}=\sum_{j\in\bfJ}\IZ\alpha_j\subset Q\qquad\g=\g(\sfA)=\bigoplus_{\alpha\in Q}\g_{\alpha}\]
Then,
\begin{itemize}
\item[(iii)] \[\g_{\bfJ}=\h_{\bfJ}\oplus\bigoplus_{\alpha\in Q_{\bfJ}\setminus\{0\}}\g_{\alpha}\]
\end{itemize}
Let $\sfA$ be a symmetrizable matrix with a fixed decomposition and $(-|-)$ be the standard normalized 
non--degenerate bilinear form on $\h$. 
Then,
\begin{itemize}
\item[(iv)] The restriction of $(-|-)$ to $\h_{\bfJ}$ is non--degenerate.
\end{itemize}
\end{proposition}

\begin{pf}
Since $\dim(\h_{\bfJ}'\cap\t_{\bfJ})=\dim(\z(\g_{\bfJ})=n_{\bfJ}-l_{\bfJ}$, where $n_{\bfJ}=|\bfJ|$ and $l_{\bfJ}=\rank(\sfA_{\bfJ})$,
 it follows that
\[\dim\h_{\bfJ}''=n_{\bfJ}-l_{\bfJ}\qquad \dim\h_{\bfJ}=2n_{\bfJ}-l_{\bfJ}\]
Moreover, by construction, the restriction of $\{\alpha_j\}_{j\in\bfJ}$ to $\h_{\bfJ}$ are linearly independent. Indeed, since 
$\iip{\sum c_j\alpha_j}{\t_{\bfJ}}=0$ for all $c_j\in\IC$, 
\[\iip{\sum_{j\in\bfJ}c_j\alpha_j}{\h_{\bfJ}}=0\quad\Longrightarrow\quad \iip{\sum_{j\in\bfJ}c_j\alpha_j}{\h}=0\quad\Longrightarrow\quad c_j=0\]
This proves $(i)$. The proof of $(ii)$ and $(iii)$ is clear.\\

Assume now that $\sfA$ is irreducible and symmetrizable and there exists $h\in\h_{\bfJ}$ such that
\[(h|h')=0\qquad\forall h'\in\h_{\bfJ}\] 
In particular, $(h|\alpha_{j}^{\vee})=0$ and $h\in\h_{\bfJ}'\cap\t_{\bfJ}\subset\h'_{\bfJ}$. Therefore, $h=\sum c_j\alpha_j^{\vee}$ and 
\[(\sum c_j\alpha_j^{\vee}|h')=\sum c_j(\alpha_j^{\vee}|h')=\iip{\sum c_jd_j\alpha_j}{h'}=0\]
Since the operators $\{\alpha_j\}$ are linearly independent over $\h_{\bfJ}$ and $d_j\neq0$, we have $c_j=0$ and $h=0$. 
We conclude that $(|)$ is non--degenerate on $\h_{\bfJ}$ and $(iv)$ is proved.
\end{pf}

\begin{rem}
The derived algebra $\g'_{\bfJ}=[\g_{\bfJ}, \g_{\bfJ}]$ is generated by $\{e_j,f_j, h_j\}_{j\in\bfJ}$, where $h_j=[e_j,f_j]$. Therefore, it does not depend of the choice of the subspace $\h''_{\bfJ}$. The assignment $\bfJ\mapsto\g'_{\bfJ}$ defines a structure that coincides with the one provided in \cite[3.2.2]{vtl-4}.
\end{rem}

Let now $\sfA$ be an irreducible, generalized Cartan matrix. Let $D_{\g}=D(\sfA)$ be the Dynkin diagram of $\g$, 
that is, the connected graph having $\mathbf{I}$ as vertex set and an edge between $i$ and $j$ if $a_{ij}\neq0$. For any $i\in\mathbf{I}$, 
let $\sl{2}^{i}\subset\g$ be the three--dimensional subalgebra spanned by $e_i,f_i,h_i$.\\

Any connected subdiagram $D\subseteq D_{\g}$ defines a subset $\bfJ_{D}\subset\bfI$. We would like to use the 
assignement $\bfJ\mapsto\g_{\bfJ}$ to define a $D_{\g}$--algebra structure on $\g=\g(\sfA)$.

\begin{rem}
For any subset $\bfJ$ of finite type, $\dim\h''_{\bfJ}=n_{\bfJ}-l_{\bfJ}=0$ and $\h_{\bfJ}=\h'_{\bfJ}$.
Therefore, if $\sfA$ is a generalized Cartan matrix of finite type, $\h''_{\bfJ}=\{0\}$ for any subset $\bfJ\subset\bfI$. The $D_{\g}$--algebra structure on $\g=\g(\sfA)$ is then uniquely defined by the subalgebras $\{\sl{2}^i\}_{i\in\bfI}$ and the Cartan subalgebra is defined for any subdiagram $D\subset\ D_{\g}$ by
\[\h_D=\{h_i\;|\;i\in\sfV(D)\}\]
If $\sfA$ is a generalized Cartan matrix of affine type, we obtain diagrammatic Cartan subalgebras $\h_D$, where
\[\h_D=\left\{\begin{array}{ccl}\{h_i\;|\;i\in\sfV(D)\}&\mbox{if}& D\subset D_{\g}\\ \h & \mbox{if} & D=D_{\g}\end{array}\right.\]
If $\sfA$ is an irreducible generalized Cartan matrix of hyperbolic type, \ie every submatrix 
is of finite or affine type, it is still possible to define a $D_{\g}$--algebra structure, depending upon the choice 
of the subspaces $\h''_{\bfJ}$ for $|\bfI\setminus\bfJ|=1$. 
\end{rem}

It is not always possible to define a $D_{\g}$--algebra structure for a generic matrix of order $\geq 3$. 
In order to obtain a $D_{\g}$--algebra structure on $\g=\g(\sfA)$, we have to satisfy the following condition:
\[\h_{\bfJ}\subset\mathfrak{t}_{\bfJ^{\perp}}\cap\bigcap_{\bfJ\subset\bfJ'}\h_{\bfJ'}\]
Since $\mathfrak{t}_{\bfJ^{\perp}}+\mathfrak{t}_{\bfJ}=\h$, we can always choose $\h_{\bfJ}\subseteq\mathfrak{t}_{\bfJ^{\perp}}$. 

\begin{lemma}
Assume given a $D_{\g}$--algebra structure on $\g=\g(\sfA)$. Then for any two subsets $\bfJ',\bfJ''\subset\bfI$,
\[\corank(\sfA_{\bfJ'\cap\bfJ''})\leq\corank(\sfA_{\bfJ'})+\corank(\sfA_{\bfJ''}) \]
In particular, if $\corank(\sfA_{\bfJ'})=\corank(\sfA_{\bfJ''})=0$, then $\corank(\sfA_{\bfJ'\cap\bfJ''})=0$.
\end{lemma}

\begin{pf}
The result is an immediate consequence of the estimate, given by the construction,
\[\dim(\h_{\bfJ'}\cap\h_{\bfJ''})\leq |\bfJ'\cap\bfJ''|+(\corank(\sfA_{\bfJ'})+\corank(\sfA_{\bfJ''}))\]
and the constraint
\[\h_{\bfJ'\cap\bfJ''}\subseteq\h_{\bfJ'}\cap\h_{\bfJ''}\]
\end{pf}

\begin{rem}
Indeed, it is easy to show that the symmetric irreducible Cartan matrix
\[\sfA=\left[
\begin{array}{rrrr}
2&-1&0&0\\
-1&2&-2&0\\
0&-2&2&-1\\
0&0&-1&2
\end{array}
\right]\]
does not admit any $D_{\g}$--algebra structure on $\g(\sfA)$, since $\dim\h_{23}=3$ and $\dim\h_{123}\cap\h_{234}=2$.
\end{rem}

The previous condition on the corank is not sufficient to obtain a $D_{\g}$--algebra structure on $\g(\sfA)$. Consider the 
symmetric Cartan matrix
\[\sfA=\left[\begin{array}{rrrrr}
2&-2&0&0\\
-2&2& -1 & 0\\
0& -1 & 2 & -1\\
0&0&-1& 2  
\end{array}\right]\]
$\sfA$ clearly satisfies the above condition. Nonetheless, a suitable  $\h_{12}''$, complement in $\h$ of $(\h'_{12}+\mathfrak{t}_{12})$, should satisfies:
\[\h_{12}''\subset\h_{123}=\h_{123}'\quad\mbox{and}\quad\h_{12}''\subseteq\t_{4}=\langle\h_{12}', -2\alpha_3^{\vee}+\alpha_4^{\vee}\rangle\]
that are clearly not compatible conditions. Therefore, there is no suitable structure for $\sfA$.\\

In the following, we will consider only symmetrizable Kac--Moody algebras $\g$ that admit such a structure. 
It automatically defines an analogue structure on $\EK{\g}$.

%--------------------------------
% QCQTQBA ON Uhg
%--------------------------------

\subsection{qCqtqba structure on $\DJ{\g}$}

Given a fixed $D_{\g}$--structure on the Kac--Moody algebra $\g$, the quantum 
enveloping algebra $\DJ{\g}$ is naturally endowed with a quasi--Coxeter 
quasitriangular quasibialgebra structure of type $D_{\g}$ defined by 
\begin{itemize}
\item[(i)]  {\bf $D_\g$-algebra:} for any $D\in\Csd{D_{\g}}$, let
$\gD\subset\g$ be the corresponding Kac--Moody subalgebra. The $D_g$-algebra
structure is given by the subalgebras $\{\DJ{\gD}\}$.
\item[(ii)] {\bf Quasitriangular quasibialgebra:} the universal $R$-matrices $\{R_{\hbar,
D}\}$, with trivial associators $\Phi_D=1^{\otimes 3}$ and structural twists 
$F_{\mathcal F}=1^{\otimes 2}$.
\item[(iii)] {\bf Quasi-Coxeter:} the local monodromies are the quantum Weyl
group elements $\{S_i^{\hbar}\}_{i\in\bfI}$. The Casimir associators $\Phi
_{\G\F}$ are trivial.
\end{itemize}

We transfer this qCqtqba structure on $\Ueh{\g}$. More precisely, we define an equivalence of quasi--Coxeter categories between the representation theories of $\DJ{\g}$ and $\Ueh{\g}$.

%---------------------------
%INDUCTIVE STEP
%--------------------------

\subsection{Gauge transformations for $\g(\sfA)$}
For any $D\subset D_{\g}$, the inclusion $\g_{D}\subset\g$, defined in the previous section, lifts to an inclusion of Manin triples
\[\gD\oplus\h_D\subset\g\oplus\h\]
We denote by $\wt{\g}_{D}=(\g_{D}\oplus\h_{D},\b_{D,+},\b_{D,-})$ the Manin triple attached to $\gD$, for any $D\subseteq D_{\g}$.

\begin{theorem}\label{thm:equivDY}
There exists an equivalence of braided $D_{\g}$--monoidal categories from
\[(\{(\DC{\Ueh{\wt{\g}_B}}{\Phi_B}, \ten_B, \Phi_{B}, \sigma R_B)\}, \{(\Gamma_{BB'}, J^{BB'}_{\F})\})\]
to
\[(\{(\DC{\DJ{\wt{\g}_B}}{}, \ten_B,\id, \sigma R^{\hbar}_B)\}, \{(\Gamma^{\hbar}_{BB'}, \id)\})\]
given by $(\{\wt{F}_B\}, \{\gamma^{\F}_{BB'}\})$.
\end{theorem}

\begin{pf}
The natural transformations $\gamma_{BB'}$, $B\subseteq B'\subseteq D_{\g}$ constructed in Section \ref{s:mchain}, define, by vertical composition, a natural transformation 
\[\gamma_{BB'}^{\bf C}\in\sfNat_{\ten}({\Gamma^{\hbar}_{BB'}\circ \wt{F}_{B'}},{\wt{F}_{B}\circ\Gamma_{BB'}})\]
for any chain of maximal length
\[{\bf C}: B=C_0\subset C_1\subset \cdots \subset C_r=B'\]
Any chain of maximal length defines uniquely a maximal nested set $\F_{\bf C}\in\Mns{B,B'}$, but this is not a one to one correspondence. For example, for $D=\sfA_3$, the maximal nested set 
\[\F=\{\{\alpha_1\}, \{\alpha_3\}, \{\alpha_1,\alpha_2,\alpha_3\}\}\]
corresponds to two different chains of maximal length
\[{\bf C}_1:\{\alpha_1\}\subset\{\alpha_1\}\sqcup\{\alpha_3\}\subset\sfA_3\qquad {\bf C}_2:\{\alpha_3\}\subset \{\alpha_1\}\sqcup\{\alpha_3\}\subset \sfA_3\]
In order to prove that the natural transformations $\gamma$ define a morphism of braided $D_{\g}$--monoidal categories, we need to prove that the transformation $\gamma_{BB'}^{\bf C}$ depend only on the maximal nested set corresponding to ${\bf C}$.\\ 

In particular, we have to prove that, for any $B_1\perp B_2$ in $\ID$, the construction of the fiber functor 
\[
\xymatrix{& \C_{B_1\sqcup B_2} \ar[dr]^{F_{B_1, B_1\sqcup B_2}} \ar[dl]_{F_{B_2, B_1\sqcup B_2}} \ar[dd] & \\ \C_{B_1} \ar[dr]_{F_{B_1}} & & \C_{B_2} \ar[dl]^{F_{B_2}} \\ & \C_{\emptyset} & }
\]
is independent of the choice of the chain. In our case,
\[\C_{B_1\sqcup B_2}=\DC{\Ueh{\wt{\g}_{B_1}}\ten\Ueh{\wt{\g}_{B_2}}}{}\]
and the braided tensor structure is given by product of the braided tensor structures on 
\[\C_{B_1}=\DC{\Ueh{\wt{\g}_{B_1}}}{\Phi_{B_1}}\qquad\C_{B_2}=\DC{\Ueh{\wt{\g}_{B_2}}}{\Phi_{B_2}}\]
Similarly, the tensor structure on the forgetful functor 
\[\C_{B_1\sqcup B_2}\to \C_{B_i}\qquad i=1,2\]
is obtained killing the tensor structure on $\C_{B_i}$, $i=1,2$, \ie applying the tensor structure on $\C_{B_i}\to\C_{\emptyset}$. In particular, the tensor structure on $F_{B_1}\circ F_{B_1,B_1\sqcup B_2}$ and $F_{B_2}\circ F_{B_2,B_1\sqcup B_2}$ coincide, since $[\wt{\g}_{B_1},\wt{\g}_{B_2}]=0$.\\

Analogously we have an equality of natural transformation
\[\gamma_{B_1}\circ\gamma_{B_1,B_1\sqcup B_2}=\gamma_{B_2}\circ\gamma_{B_2,B_1\sqcup B_2}\]

Therefore, for any maximal nested set $\F\in\Mns{B,B'}$, it is well defined a natural transformation
\[
\gamma_{BB'}^{\F}\in\sfNat_{\ten}({\Gamma^{\hbar}_{BB'}\circ \wt{F}_{B'}},{\wt{F}_{B}\circ\Gamma_{BB'}})
\]
so that the data $(\{\wt{F}_B\}, \{\gamma_{BB'}^{\F}\})$ define an isomorphism of $D$--categories from $\{\DC{\Ueh{\wt{\g}_B}}{\Phi_{B}}\}$ to $\{\DC{\DJ{\wt{\g}_B}}{}\}$.
\end{pf}

\subsection{Extension to Levi subalgebras.} In analogy with \cite[Thm. 9.1]{vtl-4}, we want to show that
the relative twists and the Casimir associators are weight zero elements. This corresponds to show that
the corresponding tensor functors $\Gamma$ and the natural transformations $\gamma$ lift to the level
of Levi subalgebras:
\[\g_D\subset\ll_D=\n_{D,+}\oplus\h\oplus\n_{D,-}\subset\g\]

\begin{proposition}\label{prop:extLevi}
The relative twists and the Casimir associators are weight zero elements.
\end{proposition}

\begin{pf}
For $D=\emptyset$, the statement reduces to prove that the Etingof--Kazhdan functor preserves
the $\h$--action \cite[Thm. 4.3]{ek-6}. The result is a consequence of Proposition \ref{pr:ab-ext}
applied to Levi subalgebras.
\end{pf}

\subsection{Reduction to category $\O^{\operatorname{int}}$}
%==============================================

The Etingof--Kazhdan functor gives rise, by restriction, to an equivalence of categories
\[\wt{F}:\O_{\g}[[\hbar]]\to\O_{\DJ{\g}}\]
We will show now that this equivalence can be further restricted to integrable modules in category $\O$, \ie modules in category $\O$ with a locally nilpotent action of the elements $\{e_i,f_i\}_{i\in\bfI}$ (respectively $E_i,F_i$).

\begin{proposition}\label{prop:Oint-red}
The Etingof--Kazhdan functor restricts to an equivalence of braided tensor categories 
\[\wt{F}:\O^{\operatorname{int}}_{\g}[[\hbar]]\to\O^{\operatorname{int}}_{\DJ{\g}}\]
which is isomorphic to the identity functor at the level of $\h$--graded $\sfk[[\hbar]]$--modules. 
\end{proposition}

\begin{pf}
Let $V\in\Oint_{\g}$. Then, the elements $e_i,f_i$ for $i\in\bfI$ act nilpotently on $V$. Then, by \cite{K}, for all $\lambda\in\mathsf{P}(V)$, there exist $p,q\in\IZ_{\geq0}$ such that
\[\{t\in\IZ\;|\;\lambda+t\alpha_i\in \mathsf{P}(V)\}=[-p,q]\]
Since the Cartan subalgebra $\h$ is not deformed by the quantization, the functor $\wt{F}$ preserves the weight decomposition. In $\DJ{\g}$, for any $h\in\h$ and $i\in\bfI$, we have
\[[h,E_i]=\alpha_i(h)E_i\]
Therefore the action of the $E_i$'s on $V$ is locally nilpotent. The action of the $F_i$'s is always locally nilpotent, since 
\[\mathsf{P}(V)\subset\bigcup_{s=1}^{r}\mathsf{D}(\lambda_s)\]
The result follows.
\end{pf}

\begin{corollary}\label{cor:catOequiv}
\return
\begin{itemize}
\item[(i)] There exists an equivalence of braided $D_{\g}$--monoidal categories between
\[\O:=(\{(\Oint_{\g_B}, \ten_B, \Phi_{B}, \sigma R_B)\}, \{(\Gamma_{BB'}, J^{BB'}_{\F})\})\]
and
\[\O_{\hbar}:=(\{(\Oint_{\DJ{\g_B}}\}, \ten_B,\id, \sigma R^{\hbar}_B)\}, \{(\Gamma^{\hbar}_{BB'}, \id)\})\]
\item[(ii)]There exists an isomorphism of $D_{\g}$--algebras
\[\Psi_{\F}:\wh{\DJ{{\g}}}\to\wh{\Ueh{{\g}}}\] 
such that $\Psi_{\F}(\wh{\DJ{{\g_{D_i}}}})=\wh{\Ueh{\g_{D_i}}}$ for any $D_i\in\F$, where the completion is taken with respect to the integrable modules in category $\O$.
\end{itemize}
\end{corollary}

\subsection{Quasi--Coxeter structure}
%============================

The previous equivalence of braided $D_{\g}$--monoidal categories induces on 
\[\O=(\{(\Oint_{\g_B}[[\hbar]], \ten_B, \Phi_{B},  \sigma R_B)\}, \{(\Gamma_{BB'}, J^{BB'}_{\F})\})\]
a structure of quasi--Coxeter category of tipe $D_{\g}$, given by the Casimir associators 
$\Phi_{\G\F}\in\sfNat_{\ten}(\Gamma_{\F},\Gamma_{\G})$ and the local monodromies 
$S_i\in\sfEnd{\Gamma_i}$ defined for any $\G,\F\in\Mns{B,B'}$ and 
$i\in\bfI(D)$ by
\[\wt{F}_{B}(\Phi_{\G\F})=(\gamma^{\F}_{BB'})^{-1}\circ\gamma_{BB'}^{\G}
\qquad
S_i=\Psi\oEK_{i}(S_i^{\hbar})\]
where $\Psi\oEK_i:\wh{\DJ{\sl{2}^i}}\to\wh{\Ueh{\sl{2}^i}}$ is the isomorphism induced at the 
$\sl{2}^i$ level by the Etingof--Kazhdan functor.

\begin{proposition}\label{prop:catOqC}
The equivalence of braided $D_{\g}$--monoidal categories $\O\to\O_{\hbar}$ 
induces a structure of quasi--Coxeter category on $\O$.
\end{proposition}

\begin{pf}
In order to prove the proposition, we have to prove the compatibility relations of the 
elements $\Phi_{\G\F}$, $S_i$ with the underlying structure of braided $D_{\g}$--
monoidal category on $\O$.\\

The element $S_i$'s satisfy the relation
\begin{equation*}
\Delta_{\F}(S_i)=(R_i)^{21}_{\F}\cdot(S_i\ten S_i)
\end{equation*}
since $\Psi_{\F}$ is given by an isomorphism of braided $D$--monoidal categories 
and therefore
\begin{equation*}
\Psi_{\F}((R_i^{\hbar})_{\F})=(R_i)_{\F}
\end{equation*}
Similarly, the braid relations are easily satisfied, since 
\[\sfAd{\Phi_{\G\F}}\Psi_{\F}=\Psi_{\G}\]
The elements $\Phi_{\F\G}$ defined above satisfy all the required properties:\\
\begin{itemize}
\item[(i)] {\bf Orientation} For any elementary pair $(\F,\G)$ in $\Mns{B,B'}$
\[\wt{F}_B(\Phi_{\F\G})=(\gamma_{BB'}^{\F})^{-1}\circ\gamma_{BB'}^{\G}=\big(\wt{F}_B(\Phi_{\G\F})\big)^{-1}\]
\item[(ii)] {\bf Coherence} For any $\F,\G,\H\in\Mns{B,B'}$
\begin{align*}
\wt{F}_B(\Phi_{\F\G})=(\gamma_{BB'}^{\F})^{-1}\gamma_{BB'}^{\H}\circ(\gamma_{BB'}^{\H})^{-1}&\circ\gamma_{BB'}^{\G}=\\=&\wt{F}_B(\Phi_{\F\H})\circ H_{B}(\Phi_{\H\G})
\end{align*}
This property implies the coherence.\\
\item[(iii)] {\bf Factorization.} Clear by construction.\\
\end{itemize}
Finally, the elements $\Phi_{\G\F}$ satisfy
\[\Delta(\Phi_{\G\F})\circ J_{\F}=J_{\G}\circ\Phi_{\G\F}^{\ten 2}\]
because they are given by composition of invertible natural tensor transformations. 
\end{pf}

\subsection{Normalized isomorphisms}
%============================

In the completion $\wh{\Ueh{\sl{2}^i}}$ with respect to category $\O$ integrable
modules, there are preferred element $S_{i,C}$ 
\[S_{i,C}=\wt{s_i}\exp(\frac{\hbar}{2}C_i)\]
where
\[\wt{s_i}=\exp(e_i)\exp(-f_i)\exp(e_i)\qquad C_i=\frac{(\alpha_i,\alpha_i)}{2}(e_if_i+f_ie_i+\frac{1}{2}h_i^2)\]

\begin{proposition}\label{prop:norm-isom}
There exists an equivalence of quasi--Coxeter categories of type $D_{\g}$ between
\[\O:=(\{(\Oint_{\g_B}, \ten_B, \Phi_{B}, \sigma R_B)\}, \{(\Gamma_{BB'}, J^{BB'}_{\F})\}, \{\Phi_{\G\F}\}, \{S_{i,C}\})\]
and
\[\O_{\hbar}:=(\{(\Oint_{\DJ{\g_B}}\}, \ten_B,\id, \sigma R^{\hbar}_B)\}, \{(\Gamma^{\hbar}_{BB'}, \id)\}, \{\id\}, \{S^{\hbar}_i\})\]
\end{proposition}

\begin{pf}
Using the result of Proposition \ref{prop:catOqC}, it is enough to prove that the natural transformation
$\gamma_i$
\[
\xymatrix{
\Oint_i \ar[rr]^{\wt{F}_i}\ar[dr]_{\sf f}^{}="A"& &\ar[dl]^{{\sf f}_{\hbar}} \ar@{=>}"A"^{\gamma_i} \Oint_{i,\hbar}\\ & \A & }
\]
can be modified in such a way that the induced isomorphism at 
the level of endomorphism algebras
$\wh{\DJ{\sl{2}^i}} \to \wh{\Ueh{\sl{2}^i}}$ maps $S_i^{\hbar}$ to $S_{i,C}$.
The natural transformation used in Corollary \ref{cor:catOequiv} induces the
Etingof--Kazhdan isomorphism
\[\Psi_i\oEK:\wh{\DJ{\sl{2}^i}} \to \wh{\Ueh{\sl{2}^i}}\]
which is the identity mod $\hbar$ and the identity on the Cartan subalgebra.
As above, we denote by $S_i$ the element $\Psi_i\oEK(S_i^{\hbar})$.
Then $S_i\equiv\wt{s_i}\mod\hbar$ and, by \cite[Proposition 8.1, Lemma 8.4]{vtl-4},
we have
\[S_i^2=S_{i,C}^2 \qquad S_i=\sfAd{x}(S_{i,C})\]
on the integrable modules in category $\O$, for $x=(S_{i,C}\cdot S_i^{-1})^{\frac{1}{2}}$. Therefore, the modified isomorphism
\[\Psi_i:=\sfAd{x}\circ\Psi_i\oEK\]
maps $S_i^{\hbar}$ to $S_{i,C}$. Moreover, $\Psi_i$ correspond with the natural transformation 
given by the composition of $\gamma_i$ with $x\in\wh{\Ueh{\sl{2}^i}}=\sfEnd{\sf f}$
\[
\xymatrix{
\Oint_i \ar@{=}[r] \ar[d]^{\sf f}="A"& \Oint_i \ar[r]^{\wt{F}} \ar[d]^{\sf f}="B" & \Oint_{i,\hbar} \ar[d]^{{\sf f}_{\hbar}}="C" \ar@{=>}"C";"B"_{\gamma_i} \ar@{=>}"B";"A"_{x}\\
\A \ar@{=}[r]& \A \ar@{=}[r] & \A
}
\]
The result follows substituting $\gamma_i$ with $x\circ\gamma_i$ in Proposition \ref{prop:catOqC}.
\end{pf}

\subsection{The main theorem}
%=======================

We now state in more details the main theorem of the paper and summarize the proof 
outlined in the previous results.

\begin{theorem}\label{thm:main}
Let $\g$ be a symmetrizable Kac--Moody algebra with a fixed
$D_{\g}$--structure and $\DJ{\g}$ the corresponding Drinfeld--Jimbo
quantum group with the analogous $D_{\g}$--structure. For any choice 
of a Lie associator $\Phi$, there exists an equivalence of quasi--Coxeter 
categories between
\[\O:=(\{(\Oint_{\g_B}, \ten_B, \Phi_{B}, \sigma R_B)\}, \{(\Gamma_{BB'}, J^{BB'}_{\F})\}, \{\Phi_{\G\F}\}, \{S_{i,C}\})\]
and
\[\O_{\hbar}:=(\{(\Oint_{\DJ{\g_B}}\}, \ten_B,\id, \sigma R^{\hbar}_B)\}, \{(\Gamma^{\hbar}_{BB'}, \id)\}, \{\id\}, \{S^{\hbar}_i\})\]
where $\ten_B$ denotes the standard tensor product in $\Oint_{\g_B}$ and
\begin{eqnarray*}
S_{i,C}&=&\wt{s_i}\exp(\frac{h}{2}\cdot C_i)\\
\Phi_B&=& 1 \mod\hbar^2 \\
R_B&=&\exp(\frac{\hbar}{2}\Omega_D)\\
\operatorname{Alt}_2 J_{\F}^{BB'}&=& \frac{\hbar}{2}\big(\frac{r_{B'}-r_{B'}^{21}}{2}-\frac{r_{B}-r_{B}^{21}}{2}\big)
\end{eqnarray*}
and $\Phi_{\G\F}$, $J_{\F}^{BB'}$ are weight zero elements.
\Omit{
\begin{itemize}
\item[(i)] for any choice of a Lie associator $\Phi$, there exists an equivalence of quasi--
Coxeter categories between
\[\O:=(\{(\Oint_{\g_B}, \ten_B, \Phi_{B}, \sigma R_B)\}, \{(\Gamma_{BB'}, J^{BB'}_{\F})\}, \{\Phi_{\G\F}\}, \{S_{i,C}\})\]
and
\[\O_{\hbar}:=(\{(\Oint_{\DJ{\g_B}}\}, \ten_B,\id, \sigma R^{\hbar}_B)\}, \{(\Gamma^{\hbar}_{BB'}, \id)\}, \{\id\}, \{S^{\hbar}_i\})\]
where $\ten_B$ denotes the standard tensor product in $\Oint_{\g_B}$ and
\begin{eqnarray*}
S_{i,C}&=&\wt{s_i}\exp(\frac{h}{2}\cdot C_i)\\
\Phi_B&=& 1 \mod\hbar^2 \\
R_B&=&\exp(\frac{\hbar}{2}\Omega_D)\\
\operatorname{Alt}_2 J_{\F}^{BB'}&=& \hbar\big(\frac{r_{B'}-r_{B'}^{21}}{2}-\frac{r_{B}-r_{B}^{21}}{2}\big)
\end{eqnarray*}
and $\Phi_{\G\F}$, $J_{\F}^{BB'}$ are weight zero elements.
\item[(ii)] for $\Phi=\Phi^{\operatorname{KZ}}$, the algebra $\wh{\DJ{{\g}}}$ is equivalent to a quasi--Coxeter quasitriangular quasibialgebra of type $D_{{\g}}$ of the form
\[(\{\wh{\Ueh{{\gD}}}\}, \Delta_0,\{\Phi\oKZ_D\},\{R\oKZ_D\}, \{J_{\F}^{BB'}\}, \{\Phi_{\G\F}\}, \{S_{i,C}\})\]
where the completion is taken with respect to the integrable modules in category $\O$.
\end{itemize}
}
\end{theorem}

\begin{pf}
The existence of an equivalence is a consequence of the constructions of Section \ref{s:mchain} and proved in Theorem \ref{thm:equivDY} and Proposition \ref{prop:catOqC}, \ref{prop:norm-isom}, 
concerning the local monodromies $S_{i,C}$.\\

The properties of associators $\Phi_B$ and $R$--matrices $R_{B}$ are direct consequences of 
the construction in Section \ref{s:ek},\ref{s:Gamma}. The relation satisfied by the relative twists $J_{\F}^{BB'}$ is proven by a simple application of Proposition \ref{pr:jets-rel-twist} and Corollary 
\ref{cor:twist-alt}. It is easy to check that the $1$--jet of the twist $J_{\F}^{BB'}$ differs from the 
$1$--jet of the twist $J^{BB'}$ (as defined in Section \ref{s:Gamma}) by a symmetric element
that cancels out computing the alternator. Therefore, Corollary \ref{cor:twist-alt} holds for 
$J^{BB'}_{\F}$ as well.\\

Finally, as previously explained, the weight zero property of the relative twists $J_{\F}^{BB'}$
and the Casimir associators $\Phi_{\G\F}$ is proved in Proposition \ref{pr:ab-ext},
\ref{prop:extLevi}. This complete the proof of Theorem \ref{thm:main}. 
\end{pf}

%---------------------------------------------
% BIBLIOGRAPHY
%--------------------------------------------

\end{document}